\documentclass[10pt]{amsart}

\usepackage{amssymb} 
\usepackage[latin1]{inputenc}
\usepackage[matrix,arrow,curve]{xy}
\usepackage[mathscr]{eucal}
\usepackage{epic,eepic}
\usepackage{bbm}
\usepackage{multirow}
\usepackage{graphicx}
% ...
\makeindex

\newtheorem{itheorem}{Theorem}
\renewcommand{\theitheorem}{\roman{itheorem}}
\newtheorem{icor}{Corollary}

\newtheorem{theorem}{Theorem}[subsection]

\newtheorem{lemma}[theorem]{Lemma}
\newtheorem{prop}[theorem]{Proposition}

\newtheorem{cor}[theorem]{Corollary}

\theoremstyle{definition}

\newtheorem{idefinition}[icor]{Definition}

\newtheorem{iexample}[icor]{Example}

\newtheorem{definition}[theorem]{Definition}
\newtheorem{remark}[theorem]{Remark}
\newtheorem{example}[theorem]{Example}

\newenvironment{pf}
{\medskip\noindent {\it Proof --- \ }}
{\hfill\nobreak $\Box$ \par\bigbreak}

\newtheoremstyle{numero}
{1ex}
{1ex}
{\it}
{1cm}
{\scshape}
{}
{4pt}
{}
\theoremstyle{numero}
\newtheorem{num}{}[subsection]

\newtheorem{sousnum}{}[num]

\newcommand{\Rc}{{\mathcal R}}

\newcommand{\F}{ \mathbb F}

\newcommand{\Q}{{ \mathbb Q } }
\newcommand{\Z}{{ \mathbb Z  }}
\newcommand{\N}{{ \mathbb N  }}

\newcommand{\m}{{\mathfrak m}}

\renewcommand{\ker}{{\rm{Ker}\,}}
\newcommand{\Gal}{{\mathrm{Gal}\,}}

\newcommand{\M}{{\mathcal M}}

\newcommand{\Aut}{{\text{Aut}\,}}
\newcommand{\End}{{\text{End}}}
\newcommand{\Id}{{\text{Id}}}

\newcommand{\U}{{\text{U}}}

\newcommand{\Gl}{{\text {GL}}}

\newcommand{\GL}{{\text {GL}}}
\newcommand{\PGL}{{\text {PGL}}}
\newcommand{\SL}{{\text {SL}}}

\newcommand{\PSL}{{\text {PSL}}}

\newcommand{\spec}{{\text{Spec\,}}}

\newcommand{\CC}{{\mathcal{C}}}

\newcommand{\G}{{\mathfrak g}}

\newcommand{\Ind}{{\text{Ind}}}
\newcommand{\ad}{{\text{ad}}}

\newcommand{\Frob}{{\text{Frob\,}}}

\newcommand{\tr}{{\rm{tr\,}}}

\newcommand{\rhob}{{\bar \rho}}

\newcommand{\rad}{{\mathrm{rad}}}

\newcommand{\univ}{{\text{univ}}}
\newcommand{\cycl}{{\text{cycl}}}

\renewcommand{\Gal}{{\rm Gal}}

\newcommand{\chib}{\bar \chi}

\renewcommand{\G}{G}

\newcommand{\Kab}{\text{Kab}}

\newcommand{\Fc}{{\mathcal F}}

\newcommand{\SR}{SR}
\newcommand{\Gb}{{\overline{G}}}
\newcommand{\Gbb}{\overline{\Gb}}

\newcommand{\uker}{{\underline{\rm{Ker}}}\,}

\newcommand{\smat}[1]{ \left( \begin{smallmatrix} #1 \end{smallmatrix}  \right)}
\newcommand{\mat}[1]{\left( \begin{matrix} #1 \end{matrix} \right)} 
\newcommand{\mato}[1]{\left( \begin{matrix} #1 \end{matrix}\right)^{\hspace{-0.1cm} 0}} 
\newcommand{\smato}[1]{\left( \begin{smallmatrix} #1 \end{smallmatrix}\right)^{\hspace{-0.05cm} 0}} 
\newcommand{\ess}{\text{ess}}
\newcommand{\spe}{\text{spe}}

\newcommand{\PGb}{\Gb}
\begin{document}

\title[Images of pseudorepresentations]{Image of pseudorepresentations and coefficients of modular forms modulo $p$}

\renewcommand{\theitheorem}{\Roman{itheorem}}

\begin{abstract} We describe the image of general families of two-dimensional representations over compact semi-local rings.
Applying this description to the family carried by the universal Hecke algebra acting on the space of modular forms of level $N$ modulo a prime $p$, we prove new results about 
the coefficients of modular forms mod $p$. If $f=\sum_{n=0}^\infty a_n q^n$ is such a form, for which we can assume without loss of generality that $a_n=0$ if $(n,Np)>1$, calling $\delta(f)$ the density of the set of primes $\ell$ such that $a_\ell \neq 0$,
we prove that $\delta(f)>0$ provided that $f$ is not zero (and if $p=2$, not a multiple of $\Delta$). More importantly, we prove, when $p>2$,
a {\it uniform} version of this result, namely that 
there exists a constant $c>0$ depending only on $N$ and $p$ such that $\delta(f)>c$ for all forms $f$ except for those in an explicit subspace of infinite codimension of the space of all modular forms mod $p$ of level $N$. Forms in this subspace, called {\it special} modular forms mod $p$,
are proved to be closely related to certain classes of modular forms mod $p$ previously studied by the author, Nicolas and Serre, called cyclotomic and CM modular forms mod $p$. 
\end{abstract}

\baselineskip16pt
\subjclass[2000]{11R}

\author[J.~Bella\"iche]{Jo\"el Bella\"iche}
\thanks{Joël Bellaïche was supported by NSF grants DMS 1405993.}
\email{jbellaic@brandeis.edu}
\address{Mathematics Department\\Brandeis University\\
415 South Street\\Waltham, MA 02454-9110\\U.S.A}
\maketitle

\setcounter{tocdepth}{1}
\tableofcontents

\section{Introduction} 

This article has two parts. In the first, we describe  the image of general families of two-dimensional representations of a pro-finite group. In the second, we use these descriptions  to study the behavior of the coefficients at primes of modular forms modulo an odd prime $p$, focussing especially on results which are {\it uniform} in the modular form.

\subsection{Image of family of representations}

\label{introimage}
Let $\Pi$ be a profinite group, $A$ a compact  local ring of maximal ideal $\m$. The residue field $\F=A/\m$ is thus a finite field, and we assume throughout \S\ref{introimage}  that its characteristic $p$ is different from $2$. 

The families we are interested in are families of two-dimensional representations of $\Pi$ carried by $A$. As past work using family of Galois representations has made clear, it is important for many applications to consider not only families of representations that can be described as a representation of $\Pi$ on a rank-two free $A$-module, but more generally two-dimensional pseudo-representations of $\Pi$ over $A$. Hence we consider a family defined as a continuous two-dimensional 
pseudo-representation\footnote{We use Chenevier's notion \cite{chenevier} of pseudo-representations, which is the most general and the most elegant, though since we assume $p>2$ for most of this paper, Chenevier's notion is equivalent to Rouquier's one.} $(t,d)$ of $\Pi$ over $A$. 

We put certain restrictions to the family we consider. First, the residual representation of the family may be irreducible or 
the sum of two characters. In the latter case, we assume that those two characters are distinct, and also that $\Pi$ satisfies the $p$-finiteness condition of Mazur. Second, we assume that as a topological $W(\F$)-algebra, $A$ is generated by $t(\Pi)$. Third, we assume that $d$ is constant, that is for every $g \in \Pi$, $d(g)$ is the Teichmüller lift of $\bar d(g)$. The last two are not serious restrictions: the second assertion can always be made true by replacing $A$ by its sub-algebra generated by $t(\Pi)$, the third by twisting $(t,d)$ by a suitable character.

Though it is not always true that $(t,d)$ comes from a representation $\rho : \Pi \rightarrow \GL_2(A)$, there always exists a Generalized Matrix Algebra (or GMA, see \cite[\S1]{BC} or below, \S\ref{subGMA}) $R$ over $A$ and a representation $\rho : \Pi \rightarrow R^\ast$ with trace $t$ and determinant $d$. We may assume that $R$ is faithful (see below~\ref{subGMA}), and generated as an $A$-module by $\rho(\Pi)$, and if we do, $R$ and $\rho$  are unique up to unique isomorphism, $R$ has a natural topology and the representation $\rho$ is continuous.

We set $G:=\rho(\Pi)$ and call this closed subgroup of $R^\ast$ the {\it image} of our family $(t,d)$. The aim is to describe as precisely as possible the group $G$. We shall handle this group using a slight generalization (from the case $R=M_2(A)$ to the case of arbitrary GMAs) 
of the remarkable theory of Lie Algebras of Pink (see \S\ref{sectionpinkpseudo}). 
This theory attaches to every closed subgroup $\Gamma$ of $\SR^1:= \{x \in R^\ast,\ \det x=1,\ x \equiv \Id \pmod{\rad R}\}$ a closed Lie subring $L=L(\Gamma)$ of $(\rad R)^0= \{x \in \rad R, \tr x = 0\}$. Contrarily to the classical theory of Lie algebras, the subgroup $\Gamma$ is not uniquely determined by $L=L(\Gamma)$. However, its closed derived subgroup $\Gamma_2$ is, as well as all the further terms of its descending central series, so that the knowledge of $L$ gives us a good, if partial, grasp on what $\Gamma$ is. We  apply this theory to the subgroup $\Gamma = G \cap SR^1$, which has finite index in $G$. 

\par \bigskip
We obtain a complete description of the Lie ring $L$ after extending the scalars from $\Z_p$ to $W(\F)$, the ring of Witt vectors of the finite filed $\F$. 
Note that $W(\F)/\Z_p$ is only a small extension, finite and unramified, which is harmless in the applications to modular forms (we do not extend the scalars to $A$, which would be much more destructive). The description of $W(\F)L$ we obtain depends, unsurprisingly, of the nature of the projective image of the representation $\rhob$. There are five cases to consider, according to the projective image being {\it exceptional} (that is, either isomorphic to $A_4$, $S_4$, $A_5$) or {\it large} (that is isomorphic to $\PGL_2(\F_q)$ or $\PSL_2(\F_q)$ for some subfield $\F_q$ of $\F$), {\it dihedral of order $>4$}, {\it dihedral of order 4}, {\it cyclic of order $>2$}, {\it cyclic of order $2$}. 

Rather than giving all the results, which the reader will find in Theorems~\ref{structureLZ2Z}, \ref{structureLcyclic}, \ref{structureLZ2ZZ2Z}. \ref{structureLdihedral} and \ref{structureLlargeimage}, let us just illustrate them by giving two examples:
\begin{itemize} \item in the {\bf large} or {\bf exceptional} projective image case, we prove that there exists a closed $W(\F)$-submodule $I_1$ of $A$ such that $I_1^2 \subset I_1$ and $W(\F) L=\mato{I_1 & I_1 \\ I_1 & I_1}$ is the set of matrices of trace $0$ with coefficients in $I_1$;
\item in the {\bf cyclic of order $>2$} projective image case, we can write the GMA $R=\mat{A & B \\ C & A}$ with $B,C$ two $A$-modules with a bilinear map $B \times C \rightarrow A$ denoted as multiplication, and we prove that there exists a $W(\F)$-module $I_1$ such that $BC \subset I_1 \subset A$ satisfying $I_1^3 \subset I_1$ and $W(\F) L=\mato{I_1 & B \\ C & I_1}$.
\end{itemize}

Moreover, we prove in each case that the description of $W(\F) L$ we obtain is optimal, in the sense that any $W(\F)$-Lie algebra satisfying the given description can be obtained from a family of representations of the type considered. In other words, nothing more can be said on $W(\F) L$.

In many cases (for instance when $\F = \F_p$ or when the projective image of $\rhob$ is large, or when this image is cyclic of order $n$ such that $\gcd(n,p-1)>2$, etc.) we obtain, better than a description of $W(\F) L$, a description of $L$ which we again prove to be optimal. We refer the reader to the Theorems cited above for the precise statements.
\par \bigskip
Recently there has been a surge in activity concerning the study of the image of families of Galois representations, represented by papers by Hida \cite{hida}, Lang \cite{lang}, and Conti-Iovita-Tilouine \cite{CIT}.  In these articles, the authors study the image of families of Galois representations attached to Hida or Coleman families of modular forms. Among the five possibilities concerning the projective image of $\rhob$ enumerated above, these authors only consider two, namely the cases when the projective image of $\rhob$ is large/exceptional or dihedral of order $>4$. Their main result is that except if all forms in the family is CM, and under various supplementary assumptions, the image $G$ of the family is large, in the following sense: there is an explicit subring $A_0$ of $A$ such that the family of representations is virtually defined over $A_0$ (i.e. is defined over $A_0$ after restricting  it to an open subgroup $\Pi_0$ of the Galois group, which is explicit in their work), and the image $G_0$ of $\Pi_0$ contains a non-trivial {\it congruence subgroup} of $\SL_2(A_0)$.  (Actually, the result of Conti-Iovita-Tilouine is slightly weaker, as it only proves this for $G_0$ replaced by its Zariski closure).

In Section~\ref{sectioncongruencelarge} (which is not used in the rest of the paper), we prove a similar result  in the case where $\rhob$ is  large or exceptional, dihedral of order $>4$ and cyclic of order $>2$. In the two remaining cases (cyclic of order $2$, and dihedral of order $4$), we show in section~\ref{sectionexample} that no result of this type is to be expected. Our result is  more general than the ones mentioned above in that it works for almost arbitrary families of representations of an arbitrary profinite group $\Pi$, instead of only specific families of representations of the absolute Galois group of $\Q$ (though it fails to deal with a few representations that Lang's result is able to deal with). Dually, our methods are much more elementary, in that we use only basic group theory and Pink's theory of Lie algebras, rather than the theory of classical and $p$-adic modular forms, the structure of the Galois group, and advanced Hodge-Tate theory as in the afore-mentioned articles.

\subsection{Coefficients of modular forms}

\subsubsection{Individual density result}

Let $N \geq 1$ be an integer, $p$ any prime, $k \in \Z/(p-1)\Z$. For $\F$ a finite extension of $\F_p$, 
we shall denote by $M_k(N,\F)$  the algebra of modular forms of level $\Gamma_0(N)$, weight $k$, with coefficients in $\F$, in the sense of Swinnerton-Dyer. If $f= \sum_{n=0}^\infty a_n q^n$ is an element of $M(N,\F)$, then the set $\{ \ell \text{ prime }, a_\ell \neq 0\}$ is Frobenian, as was known already to Serre in the seventies (cf. \cite{SerreEM}), and therefore has a density, which is a rational number between $0$ and $1$. We shall denote this number by $\delta(f)$, and refer to it as {\it the density of} $f$.

  
Let $\Fc_k(N,\F)$ be the subspace of $M_k(N,\F)$ of forms $f=\sum a_n q^n$ such that $a_n \neq 0 \Rightarrow (n,Np)=1$. Equivalently, $\Fc_k(N,\F)$ is the intersection of the kernels of the operators $U_\ell$ for $\ell$ prime, $\ell \mid Np$, defined by $U_\ell (\sum a_n q^n) = \sum a_{n\ell} q^n$ (those operators leave $M_k(N,\F)$ stable, see \cite{jochnowitz}.) When studying $\delta(f)$, there is no loss of generality in supposing $f \in \Fc_k(N,\F)$, because for any $f=\sum_{n=0}^\infty a_n q^n \in M_k(N,\F)$,
the $q$-series $$f' = \sum_{n=0,(n,Np)=1}^\infty a_n q^n$$ belongs to $\Fc_k(N^2,\F)$ and obviously satisfies $\delta(f')=\delta(f)$. We shall henceforth restrict our attention to the subspace $\Fc_k(N,\F)$ of $M_k(N,\F)$. 

\begin{iexample} We let $\Delta \in \F_p[[q]]$ be the product $q\prod_{n\geq 1}(1-q^n)^{24}$. It is  the
reduction mod $p$ of the $q$-expansion of the unique normalized cuspidal eigenform of weight 12 and level $1$, and $\Delta = \sum_{n \geq 1} \tau(n) q^n$
where $\tau$ is the reduction mod $p$ of the usual Ramanujan $\tau$-function.  One has $\Delta \in M_{12}(N,\F_p)$. Let us denote by $\Delta'$ (depending implicitly of $p$ and $N$) the $q$-series $\sum_{n \geq 1, \ (n,Np)=1} \tau(n)q^n$, which belongs to $\Fc(N,\F)$. For $p=2$, $N=1$ one has $\Delta=\Delta' = \sum_{n \text{ odd }} q^{n^2}$

\end{iexample}

\begin{itheorem} \label{densnonzero} Let $\F$ be a finite extension of $\F_p$, $k \in \Z/(p-1)\Z$ and $f \in \Fc_k(N,\F)$. Assume that $f \neq 0$ (resp. $f \not \in \F \Delta'$ if $p=2$). Then $\delta(f) > 0$.
\end{itheorem}
The theorem will be proved in \S\ref{proof1}.
\begin{icor} \label{cornonzero}  Let $f=\sum a_n q^n, g=\sum b_n q^n \in \Fc_k(N,\F)$. Assume that $a_\ell = b_\ell$ for all primes $\ell$ except for a set of density $0$ (and that $a_1=b_1$ if $p=2$). Then $f=g$.
\end{icor}
\begin{pf} Since $\delta(f-g)=0$, Theorem~\ref{densnonzero} implies $f-g=0$ if $p>2$, and $f-g \in \F \Delta'$ if $p=2$. In this case, since $a_1(f-g)=0$, $f-g=0$ as well.
 \end{pf}

%

\subsubsection{Uniformity?} 

We now turn to the main subject of this paper, the question of {\it uniformity} in the lower bound of Theorem~\ref{densnonzero}: when $f$ varies in the infinite-dimensional space $\Fc_k(N,\F)$, we know that $ \delta(f)>0 $, but is it possible that  $\delta(f)$ goes to $0$, or will $\delta(f)$ stay bounded away from $0$, at least when $f$ is supposed to stay in some large subset of $\Fc_k(N,\F)$? This question of uniformity is not only natural, but of crucial importance if we hope to obtain new results for coefficients of weakly holomorphic modular forms of half-integral weight, such as the inverse Dedekind $\eta$-function, $\eta^{-1}$, whose coefficients are the value of the partition function $p(n)$. Indeed, those weakly holomorphic modular forms are in an appropriate sense limits of classical modular forms. We plan to go back to these applications in a subsequent paper.

\begin{example} In the case $N=2$, $p=1$, the vector space $\Fc=\Fc_0(1,\F_p)$ has $(\Delta^n)_{n=1,3,5,7,\dots}$ as a basis.
It was proved by the author (letter to Nicolas and Serre, July 2012) that for $p=2$, and any integer $r\geq 1$,
\begin{num} \label{densdelta} $\delta(\Delta^{2^{r}+1})=2^{-\left\lfloor \frac{r-1}{2}\right\rfloor -2}$, \ \ \ $\delta(\Delta^{(2^{2r+1}+1)/3}) = 2^{-r-1}$ \end{num}
\noindent
Hence those forms (except perhaps a finite number of them) must be excluded if we want a positive lower bound for $\delta(f)$. For others odd powers of $\Delta$, experimental computations done with sage and certain partial results strongly suggest a different and striking pattern: it seems
 that $ \delta(\Delta^n)=1/8$ for all $n > 1$ not of the form $2^r+1$ or $\frac{2^{(2r+1)+1}}{3}$.
 \end{example} 
 
Though in this paper we are forced to exclude the case $p=2$ (both because Pink's theory requires $p>2$ and because our GMA methods require a multiplicity free hypothesis which is not satisfied if $p=2$), the example above, together with analog computations done by Medvedovski in the case $p=3$, showed that to obtain a uniform lower bound $\delta(f)>c>0$, it is necessary to exclude some exceptional forms $f$, and at the same time suggested that such a lower bound was otherwise possible. Indeed we prove:

\begin{itheorem} \label{densbounded} (cf. \S\ref{proofdensbounded}.) Let us assume that $p>2$. There exists a canonical subspace $\Fc_{k,\spe}(N,\F)$ of $\Fc_k(N,\F)$, of infinite codimension, and a constant $c>0$ (depending only on $N,\F$) such that for every modular form $f \in \Fc_k(N,\F) - \Fc_{k,\spe}(N,\F)$, one has $\delta(f)>c$.
\end{itheorem}

The constant $c$ is effective (we can take $c = \frac{p-1}{pn}$ where $n$ is the product of the orders of the image of all representations 
$\rhob\in \Rc(k,N,\F)$, see below).

The definition of the subspace $\Fc_\spe(N,\F)$, which we call the {\it subspace of special forms of $\Fc$}, is given in~\ref{defformespeciales}.  This definition uses the image of the natural Galois pseudo-representation over the semi-local Hecke algebra $A$
acting of $\Fc$, as well as Pink's Lie algebra of that image. 

This subspace is called the subspace of {\it special forms}. To analyze this subspace in more detail, we need to introduce some notations and recall some elementary facts.

\subsubsection{Decomposition of $\Fc_k(N,\F)$}

For simplicity we shall often drop the level $N$, the weight $k$ (which are fixed during all the discussion) and the finite field $\F$ from the notation and write $\Fc$ for $\Fc_k(N,\F)$, $\Fc_\spe$ for $\Fc_{k,\spe}(N,\F)$.

The space $\Fc$ is endowed with an action of the Hecke operators $T_\ell$  for $\ell \nmid Np$. After replacing $\F$ by a large enough
finite extension, we may assume (cf.~\cite{jochnowitz}) that all eigenvalues of these operators are in $\F$.
Let $A_k(\F)$ for $k \in \Z/(p-1)\Z$ be the closed $\F$-subalgebra of $\End_{\F}(\Fc_k(\F))$ generated by the Hecke operators $T_\ell$ for $\ell$ not dividing $Np$. The sequences 
$(\lambda_\ell)_{\ell\nmid Np}$ with $\lambda_\ell \in \F$ which are
 systems of eigenvalues for the operators $T_\ell$ of a common eigenvector in $\Fc$ are in bijection, by a theorem of Deligne,
 with a certain set $\Rc= \Rc(k,N,\F)$ of semi-simple continuous Galois representations $\rhob: G_{\Q,Np} \rightarrow \Gl_2(\F)$: the correspondence is given by $\lambda_\ell = \tr \rhob(\Frob_\ell)$ for all $\ell \nmid Np$. This set $\Rc(k,N,\F)$ can 
 be described as the set of all semi-simple representations $\rhob: G_{\Q,Np} \rightarrow \Gl_2(\F)$ of determinant $\omega_p^{k-1}$ and Serre's level $N$. This is the content of Serre's conjecture, now a theorem of Khare and Wintenberger.

If $\rhob$ corresponds to a system of eigenvalues $(\lambda_\ell)$, we shall denote by $\Fc_\rhob = \Fc_\rhob(N,\F)$ the {\it generalized eigenspace} in $\Fc$ for the $T_\ell$ ($\ell \nmid Np$) with eigenvalues $\lambda_\ell$, that is the set of forms $f \in \Fc$ such that $\forall \ell \nmid Np,\ \exists n \in \N,\ (T_\ell-\lambda_\ell)^n f=0$.

We thus have a decomposition
\begin{eqnarray} \label{decFck0}  \Fc = \bigoplus_{\rhob \in \Rc} \Fc_\rhob \end{eqnarray}
of $\Fc$ into generalized eigenspaces.

\subsubsection{Special modular forms in $\Fc_\rhob$}

We define $\Fc_{\rhob,\spe}$ as the space of modular forms in $\Fc_\rhob$ that are special, that is $\Fc_{\rhob,\spe} = \Fc_\rhob \cap \Fc_\spe$. The following result refines the statement that $\Fc_\spe$ is of infinite codimension given in Theorem~\ref{densbounded}.

\begin{itheorem} \label{ithmspecialinfinite} (cf. \S\ref{proofithmspecialinfinite})
 Let $\rhob$ be any representation  in $\Rc$. Assume  that $p>2$, and if $p=3$, assume also that $\rhob$ is a twist of $1 \oplus \omega_3$, where $\omega_3$ is the cyclotomic character. The space $\Fc_{\rhob,\spe}$  has infinite codimension in $\Fc_\rhob$. \end{itheorem}

\subsubsection{Special modular forms, CM forms, abelian forms}

For many representations $\rhob$, we are able to give a much more precise description of $\Fc_{\rhob,\spe}$. 

\begin{idefinition} Let  $f = \sum a_n q^n \in \Fc$. Let $K$ be a quadratic extension of $\Q$. 
We shall say that $f$ is {\it cyclotomic} (resp. $K$-abelian) if there exists a finite cyclotomic extension $L/\Q$ (resp. an abelian extension $L/K$, Galois over $\Q$) such that 
such that for $\ell$ prime not dividing $Np$, the coefficient $a_\ell$ of $f$ depends only on $\ell$ through $\Frob_{\ell,L/\Q}$.  \end{idefinition} 
 Thus, a form $f$ is cyclotomic if there exists $M \geq 1$ such that $a_\ell$ depends only on $\ell \pmod{M}$.
\begin{example} In the case $p=2$, $N=1$, it was proved by Nicolas and Serre (\cite{NS2}) that the  forms $\Delta^n$ for $n=2^r+1$ and $n=(2^{2r+1}+1)/3$ appearing in~\ref{densdelta} are $K$-abelian, and it was proved by the author (letter to Serre and Nicolas, October 2013)
that only for those odd values of $n$ were $\Delta^n$ $K$-abelian (for $K=\Q(i)$ or $K=\Q(i \sqrt{2})$) but not cyclotomic. The forms $\Delta^n$ are known to be cyclotomic for $n=13,5,7,19,21$ and conjectured not to be so for other values of $n$. There also exists forms which are $K$-abelian or cyclotomic not of the form $\Delta^n$: they have been classified and their density $\delta$ has been computed, and often goes to zero along infinite sequences of such forms.    
\end{example}

Once again, though we exclude the case $p=2$, this example suggested a close relation between the $K$-abelian and cyclotomic forms on the one hand, and the so-called special modular forms which we need to exclude in Theorem~\ref{densbounded}, in the other hand. Indeed, we prove

\begin{itheorem} \label{ithmspecial} (cf. Cor. \ref{corspered}, Cor. \ref{corspedih} and \S\ref{spelarge}). We assume $p>2$

\begin{itemize} \item If $\rhob$ has large projective image, the space of special modular forms $\Fc_{\rhob,\spe}$ is finite-dimensional.
\item If $\rhob$ has a dihedral projective image which is of order $n$ with $n>4$, $4 \mid n$, then the space of special modular forms $\Fc_{\rhob,\spe}$  contains as a finite codimension subspace the space of $K$-CM forms, where $K$ is the quadratic extension of $\Q$ fixed by the unique quotient of order $2$ of the projective image of $\rhob$.
\item If $\rhob$ has cyclic projective image which is not of order $2$, then the space of special modular forms $\Fc_{\rhob,\spe}$ is exactly the space of cyclotomic modular forms.
\end{itemize}
Moreover, in all the cases considered above, the space $\Fc_{\rhob,\spe}$ is stable by all Hecke operators.
\end{itheorem}
By contrast, in the remaining two degenerate cases where the projective image of $\rhob$ is $\Z/2\Z$ or $\Z/2\Z \times \Z/2\Z$, the space of special modular forms $\Fc_{\rhob,\spe}$ is not in general stable by all the Hecke-operators, and while it may be proved to contain all\footnote{In the case of projective image $\Z/2\Z \times \Z/2\Z$, there exists non-zero $K$-abelian forms, but no cyclotomic forms, in $\Fc_\rhob$, for exactly three quadratic fields $K$.  In the case when the projective image is of order $2$, there exists non-zero $K$-CM forms for exactly one quadratic field, plus non-zero cyclotomic modular forms. For more about 
cyclotomic forms and CM forms, see \S\ref{sectionspecial}.} cyclotomic and CM-forms in $\Fc_{\rhob}$, I do not know at this point how much larger $\Fc_{\rhob,\spe}$ is.

\subsubsection{A rough outline of the proofs}

To prove Theorems~\ref{densnonzero},~\ref{densbounded},~\ref{ithmspecialinfinite} and ~\ref{ithmspecial}, we consider the Hecke algebra $A$ acting on the space of modular forms $\Fc$ mod $p$. This is by construction a compact semi-local Hecke algebra, which carries a natural pseudo-representation $(t,d)$ of the Galois group $G_{\Q,Np}$.
A crucial ingredient is the description of the image $G$ of this pseudo-representation, or at least, of its Pink's Lie algebra, 
a special case of the general results described in~\ref{introimage} and proved in section~\ref{sectionLieStudy}.

A form $f$ in $\Fc$ defines an open and closed subset $N_f$ of the compact group $G$ (namely $N_f=\{g \in G, a_1(\tr(g) f) \neq 0\}$)   
such that $\mu_G(N_f)=\delta(f)$ (as is shown by a simple application of Chebotarev, see~\S\ref{proof1}), where $\mu_G$ is the probability Haar measure on $G$.  Theorem~\ref{densnonzero} is thus reduced to checking that $N_f$ is not empty (except when $f=0$, or in the case $p=2$, when $f$ is proportional to $\Delta'$), which is not hard (see~\S\ref{proof1}). 

To prove the other theorems we need to understand how $\mu_G(N_f)$ varies with $f$. Since we have more control on the finite index subgroup $\Gamma$ of $G$ that on $G$ itself, we cut $N_f$ into parts related to $\Gamma$-cosets. To be precise, if $X$ is a set of representatives in $G$ of $G/\Gamma$, so that $G = \coprod_{x \in X} x \Gamma$, we cut $N_f$ into pieces $N_{f,x}:= x^{-1} N_f \cap \Gamma$, so that $\mu_G(N_f) = \sum_{x \in X} \mu_G(N_{f,x})$ and our problem is to understand for a given $x$, how $\mu_G(N_{f,x})$ varies with $f$. 

Since $N_{f,x}$ is a subset of $\Gamma$, we can transport the question to the Lie Algebra $L$ of $\Gamma$, that is study instead $\mu_L(M_{f,x})$ where $M_{f,x}=\Theta(N_{f,x}) \subset L$, $\Theta$ being the `logarithm' in Pink's Lie theory. (Here I ignore, for simplicity, the fact that $\Theta$ is not always a measure-preserving bijection between $\Gamma$ and $L$. This is remedied by replacing $\Gamma$ and $L$ by $\Gamma_2$ and $L_2$, their derived subgroup and derived Lie algebra respectively. However, this changes is source of important, and essential, complications. See Remark~\ref{LL2} for a more detailed discussion of this fine point). 

The Lie algebra $L$ is an infinite-dimensional vector space over $\F_p$, and it turns out that $M_{f,x}$ is the complement in $L$ of an algebraic hypersurface of $L$ (here I am assuming $\F=\F_p$ for simplicity), 
that is a subset a polynomial on $L$ involving finitely many variables. Thus, $\mu_L(M_{f,x})$ is the proportion of points that does not lie on an  hypersurface in a finite-dimensional space over $\F_p$. Unfortunately the dimension of the ambient space as well as the degree of that hypersurface depend on $f$, and the estimates given by the Weil's conjectures proved by Deligne are not sufficient to get the desired lower bound for $\mu_L(M_{f,x})$ in general. 

However, when we choose for $x$ the image $c$ of a complex conjugation in $G$, we can show that the equation defining $M_{f,c}$ is, after a measure-preserving change of variables, affine. This is the main point of the proof of Theorem~\ref{densbounded}, and is dealt with in a more general settings in \S\ref{keymeasure}. If we denote by $M'_{f,c}$ the transform of $M_{f,c}$ by this change of variable, $\mu_L(M_{f,c})=\mu_L(M'_{f,c})$ and $M'_{f,c}$ is either empty, or an hyperplane of $L$, or $L$. In the last two cases, $\mu_L(M_{f,c}) \geq 1/p$, which gives us the desired lower bound. We need to determine for which forms $f$ we have $M_{f,c}$ empty. This is done in section~\ref{sectionessential}, relying on the explicit description of $L$ given in section~\ref{sectionLieStudy} which leads us to the notion of the essential subgroup $A_\ess$ of $A$, studied in \S\ref{sectionessential} and to the definition of {\it special modular forms} (cf. \S\ref{defformespeciales}), the forms $f$ which are orthogonal to $A_\ess$, and which happen to be the same as those for which $M_{f,c}$ is empty. This proves that forms $f$ which are non-special, the quantity $\delta(f)$ is bounded below by a positive constant independent of $f$.

To prove that the special forms are rare (cf. \S\ref{Proofthmspecialinfinite}), we need to show that $A_\ess$ is big,
and a crucial ingredient, that we borrow from recent previous works of the author, Khare, Deo, Medvedovski, inspired by Nicolas and Serre, is that each local component of $A$ is noetherian and of Krull dimension at least 2 (except when $p=2,3$, where we only know that some components have dimension at least $2$).

\par \bigskip \par \bigskip

The author is grateful to G. Chenevier, A. Conti, S. Deo, J. Lang, A. Medvedovski, P. Monsky, J.-L. Nicolas, J.-P. Serre, J. Tilouine for many useful and interesting discussions.

\section{Pseudo-representations and GMA} 
\label{sectionpseudorep}

\subsection{Reminder and complements on pseudo-representations of dimension 2}

\subsubsection{Pseudo-representations of a group}

For the general definition of a  {\it pseudo-representation} of a group $\Pi$ with values in a commutative ring $A$, we refer the reader to \cite{chenevier}. In dimension $2$, which is the only case we shall need, it is not long to recall the equivalent definition proposed in {\it loc. cit.}, Lemma 1.9: a (two-dimensional) {\it pseudo-representation} of $\Pi$ with values in $A$ is a pair of maps $t: \Pi \rightarrow A$, $d: \Pi \rightarrow A$, such that
\begin{num} \label{td1} $d$ is a group homomorphism from $\Pi$ to $A^\ast$. \end{num}
\begin{num} \label{td2} $t$ is a central function from $\Pi$ to $A$. \end{num} 
\begin{num} \label{td3} $t(1)=2$. \end{num}
\begin{num} \label{td4}  $t(xy)+d(y) t(xy^{-1}) = t(x) t(y)$ for all $x,y \in \Pi$. \end{num}
If $\Pi$ is a topological group, $A$ a topological ring, one says that the pseudo-representation $(t,d)$ is {\it continuous} if $t$ and $d$ are.
If $2$ is invertible in $A$, $d$ can be recovered from $t$ by the formula $d(x)=\frac{t(x)^2-t(x^2)}{2}$. If $\rho$ is any representation
$\Pi \rightarrow \GL_2(A)$, then it is easy to check that $(\tr \rho, \det \rho)$ is a pseudo-representation of dimension $2$.

The {\it kernel} of $(t,d)$ is defined by
$$\ker (t,d):=\{y \in \Pi, d(y)=1 \text{ and }  \forall x \in \Pi,\ t(xy)=t(x)\}.$$ By \ref{td1} and \ref{td2}, this is a normal subgroup of $\Pi$, closed if $(t,d)$ is continuous. We observe that if $2$ is invertible in $A$, we can omit the condition on $d$ in the definition of $\ker (t,d)$ as it follows from the condition on $t$.
Both the maps $t$ and $d$ factors through the quotient group $\Pi/\ker (t,d)$, and they define a  pseudo-representation of dimension $2$ of $\Pi/\ker(t,d)$ with values in $A$ whose kernel is trivial.

\subsubsection{Pseudo-representations of an algebra}

Let $R$ be an $A$-algebra (non-necessarily commutative), and let $(T,D)$ be a pair of maps $R \rightarrow A$. We say that $(T,D)$ is a {\it pseudo-representation} of dimension $2$ of $R$ with values in $A$, if
\begin{num} \label{TD1} $D(1)=1$, $D$ is multiplicative (i.e. $D(xy)=D(x)D(y)$ for $x,y \in R$) and homogeneous of degree $2$ (i.e. $D(ax)=a^2D(x)$ for $a \in A$, $x \in R)$. \end{num}
\begin{num} \label{TD2} $T$ is $A$-linear and $T(xy)=T(yx)$ for all $x,y \in R$. \end{num}
\begin{num} \label{TD3} $T(1)=2$. \end{num}
\begin{num}  \label{TD4} $D(x+y)=D(x)+D(y)+T(x)T(y)-T(xy)$ for all $x,y \in R$. \end{num}

 \begin{lemma}\label{tdTD} If $R=A[\Pi]$, the map $(T,D)\mapsto (T_{|\Pi},D_{|\Pi})$ is a bijection between the set of all pseudo-representations of dimension $2$ of $R$ and the sets of all pseudo-representations of dimension $2$ of $\Pi$.
 \end{lemma} 
\begin{pf} The proof below is closely inspired by \cite{chenevier}.

 If $(T,D)$ satisfies \ref{TD1} to \ref{TD4}, it is clear that $(T_{|\Pi},D_{|\Pi})$ satisfies \ref{td1} to \ref{td3}.
 Set $f(x,y):=T(x)T(y)-T(xy)$ for $x,y \in R$, so that \ref{TD4} becomes
\begin{num}  \label{TD5} $D(x+y)=D(x)+D(y)+f(x,y)$ for all $x,y \in R$. \end{num}
For $x,y,z \in R$ one has $D((x+y)z)=D(xz)+D(yz)+f(xz,yz)$ but also, since $D$ is multiplicative $D((x+y)z) = D(x+y) D(z) = D(xz) + D(yz)+ f(x,y) D(z)$, hence
\begin{num}  \label{TD6} $f(xz,yz)=f(x,y)D(z)$ for all $x,y,z \in R$. \end{num}
If $y$ is invertible in $R$, of inverse $z=y^{-1}$, applying \ref{TD6} gives
$f(x,y) D(y^{-1}) = f(xy^{-1},1)$. Since for every $x$, $T(x)=f(x,1)$ by \ref{TD3}, one obtains 
$T(xy^{-1})=f(x,y) D(y^{-1}) = T(xy) D(y)^{-1} - T(x)T(y) D(y)^{-1}$,
that is
\begin{num}  \label{TD7}  $T(xy)+D(y) T(xy^{-1}) = T(x) T(y)$ for all $x \in R$, $y \in R^\ast$. \end{num}
In particular, the restrictions of $T$ and $D$ to $\Pi$ satisfy \ref{td4}, hence $ (T_{|G},D_{|G})$ is a pseudo-representation of $G$ of dimension $2$.
\par \bigskip

Conversely, if $(t,d)$ is a pseudo-representation of dimension 2 of $\Pi$ with values in $A$, let us denote by $T$ the unique $A$-linear map $A[\Pi] \rightarrow A$ which coincides with $t$ on $\Pi$ and by $f$ the symmetric bilinear form on $A[\Pi]$ defined by $$f(x,y):=T(x)T(y)-T(xy).$$ For $x \in G$, one has $f(x,x)=T(x)^2-T(x^2) =2 d(x)$ by \ref{td4} and \ref{td3}. Therefore, there exists a unique quadratic form $D:A[\Pi] \rightarrow A$
such that 
\begin{num} \label{condDf} $D(x+y) -D(x)-D(y) = f(x,y)$ for all $x,y \in R$,\footnote{This condition~\ref{condDf} is expressed by saying that $f(x,y)$ is the {\it polarization} of the quadratic form $D$.}
and \end{num} \begin{num} \label{condDd} $D(g)=d(g)$ for all $g \in \Pi$. \end{num}

Thus we have defined functions $T,D$ from $A[\Pi]$ to $A$ that extends $t$ and $d$, and that satisfies \ref{TD2} to \ref{TD4}, as well as $D(1)=1$ and $D$ homogeneous of degree $2$. We now proceed to show that $D$ is multiplicative.

From \ref{td4} one gets $f(x,y)=t(xy^{-1}) d(y)$ for $x,y \in \Pi$ hence 
\begin{num} \label{td5} $f(zx,zy)=f(xz,yz) = f(x,y) d(z)$ for $x,y,z \in \Pi$. \end{num}
This relation holds more generally for $x,y,z \in A[\Pi]$ by linearity. 
For $z \in G$, the quadratic forms on $A[\Pi]$ given by $x \mapsto D(xz)$ and $x \mapsto D(x) D(z)$
have the same polarization (namely $f(x,y)d(z)$, using~\ref{td5}), and agrees on the basis $\Pi$ on $A[\Pi]$.
They are therefore equal:
\begin{num} \label{Dsemimult} $D(xz)=D(x) D(z)$ for $x \in A[\Pi]$, $z \in G$. \end{num}
Again, the quadratic forms $z \mapsto D(xz)$ and $z \mapsto D(x)D(z)$ have the same polarization by \ref{td5},
and they agree on $\Pi$ by \ref{Dsemimult}, hence they are equal.
Therefore $(T,D)$ is a pseudo-representation of $R$ with values in $A$, and the map $(t,d) \mapsto (T,D)$ is an inverse of the restriction map considered in the statement. 
\end{pf}

There is a notion of kernel for a pseudo-representation $(T,D)$ of an algebra $R$:
$$\uker (T,D) = \{y \in R,\  D(y)=0 \text{ and } T(yx)=0\  \forall x \in R\}.$$
We say that $(T,D)$ is {\it faithful} if $\uker(T,D)=0$. It is easy to see that $\uker(T,D)$ is a two-sided ideal of $R$, and that $(T,D)$ factors through $R/\uker(T,D)$ and defines a  faithful pseudorepresentation of that algebra with values in $A$.

If $(T,D)$ is a pseudo-representation of $A[\Pi]$, and $(t,d)$ is the pseudo-representation of $\pi$ obtained by restriction,  then the relation between the $\ker(t,d)$ and $\uker(T,D)$ is as follows:
\begin{lemma} \label{keruker} For $g \in \Pi$, one has $g \in \ker (t,d)$ if and only if $g-1 \in \uker(T,D)$.
\end{lemma}
\begin{pf}
If $g \in \Pi$, by linearity of trace $t(gh) = t(h)$ for all $h \in \Pi$ if and only if $T(gy)=T(y)$ for all $y$ in $R=A[\Pi]$. If the latter condition holds, then in particular $t(g)=2$, and under this condition
$d(g)=1$ and $D(g-1)=0$ are equivalent since $D(g-1)=d(g)-t(g)+1$.
\end{pf}
However, in general $\uker (T,D)$ is strictly larger than the two-sided ideal generated by the elements $g-1$, $g \in \ker (t,d)$. If $(T,D)$ is faithful then $\ker (t,d)=\{1\}$, but the converse is false in general.

We say that a pseudo-representation $(T,D)$ of $R$ is {\it Cayley-Hamilton} if for every $x \in R$, one has $x^2-T(x)x+D(x)=0$.	A
faithful pseudo-representation is Cayley-Hamilton, but the converse is false in general.

\subsection{Generalized Matrix Algebras}

\label{subGMA}
 
The notion of Generalized Matrix Algebra (GMA) is defined and studied in detail in \cite[\S1.3]{BC}. Here we will content ourselves with an {\it ad hoc} definition which is equivalent to  the notion called GMA of type $(1,1)$ in the terminology of {\it loc. cit.}

 Let $A$ be a commutative ring. Suppose given two $A$-modules $B$ and $C$, and a morphism of $A$-modules $m: B \otimes_A C \rightarrow A$ such that
 \begin{num} \label{asso} for all $b,b' \in B$ and $c,c' \in C$, $m(b,c)b'=m(b',c)b$ and $m(b,c') c = m(b,c) c'$. 
 \end{num} 
 With this data we define a not necessarily commutative $A$-algebra $R$, $R=A \oplus B \oplus C \oplus A$ as an $A$-module, endowed with the multiplication 
$$(a,b,c,d) \times (a',b',c',d') = (aa'+m(b,c'), ab'+d'b, a'c+dc', dd'+m(b',c)),$$
for $a,a',d,d' \in A$, $b,b' \in B$, $c,c' \in C$: the distributivity of multiplication over addition is obvious, the unity for multiplication is $(1,0,0,1)$, and the associativity of multiplication is easily checked using
 \ref{asso}.
We call $(A,B,C,m,R)$, or by abuse $R$, a {\it generalized matrix algebra}.
A morphism of GMAs from $(A,B,C,m,R)$ to $(A',B',C',m',R')$ is the data $(f_A,f_B,f_C)$ of a morphism of rings $f_A : A \rightarrow A'$ and two morphisms of $A'$-modules $f_B : B \otimes_A A' \rightarrow B'$ and $f_C : C \otimes_A A' \rightarrow C'$ such that $f_A(m(b,c)) = m'(f_B(b),f_C(c))$ for every $b \in B$, $c \in C$. A morphism of GMAs induces a morphism of $A'$-algebras $f_R: R \otimes_A A' \rightarrow R'$. When $A=A'$ and $f_A=\Id_A$, we say that this morphism is {\it over $A$}, or an {\it $A$-morphism}. A {\it sub-GMA} of $(A,B,C,m,R)$ is a GMA $(A',B',C',m',R')$ where $A' \subset A$, $B' \subset B$, $C' \subset C$ such that the these three inclusions maps define a morphism of GMAs. An {\it $A$-sub-GMA} is a sub-GMA where $A'=A$.

This meaning of these definitions becomes clearer if we decide to represent $(a,b,c,d)$ as a matrix $\smat{a & b \\ c & d}$, and to simply write $bc$ or $cb$ for $m(b,c)$, for then multiplication in $R$ is computed as multiplication of ordinary matrices.

\begin{lemma} \label{basicGMA1}
If in a GMA $R$, $BC=A$, then there are isomorphisms of $A$-modules form $B$ and $C$ onto $A$ so that $m$ corresponds to the multiplication $A \times A \rightarrow A$. In other words, there is an isomorphism over $A$ of GMAs $R \simeq M_2(A)$.
\end{lemma}
\begin{pf} Let $b \in B$ and $c \in C$ such that $m(b,c) = 1$; by \ref{asso}
one gets for $b' \in B$ that $b' = m (b,c) b' = m(b',c) b$ which shows that $B$ is generated by $b$; moreover if for $a \in A$, $ab=0$, then $m(ab,c)=am(b,c)=a=0$. which shows that $(b)$ is a basis of $B$. Similarly $(c)$ is a basis of $C$ and if we identifies $B$ and $C$ with $A$ using those basis, then $m$ becomes the multiplication of $A$ because $m(ab,a'c)=aa'm(b,c)=1$.
\end{pf}

We define the {\it trace} map $\tr: R \rightarrow A$ as $\tr \smat{a & b \\ c & d} = a+d$ and the {\it determinant} map $\det:R \rightarrow A$ by
$\det \smat{a & b \\ c & d} = ad-bc$. It is clear that as in the case of usual matrix algebras, one has $\tr(rr')=\tr(r'r)$, $\det(rr')=\det(r)\det(r')$ and, if $p>2$,
$\det(r)=\frac{\tr(r)^2-\tr(r^2)}{2}$. 

It is easily checked that the pair of maps $(\tr,\det) : R \rightarrow A$ is a pseudo-representation  of dimension $2$ of $R$ with values in $A$. We say that the $GMA$ $R$ is {\it faithful} (resp. {\it Cayley-Hamilton}) 
if $(\tr,\det)$ is. It is easily seen that $T$ is faithful if and only if the map $m: B \otimes_A C \rightarrow A$ being {\it non-degenerate}, meaning that the only $b \in B$ such that $m(b,c)=0$ for all $c\in C$ is $b=0$, and the only $c \in C$ such that $m(b,c)=0$ for all $b \in B$ is $c=0$.

\begin{lemma} \label{basicGMA2}
Assume that $A$ is a domain, with fraction field $K$, and that $R=\mat{A & B \\ C & D}$ is a faithful GMA over $A$.
Then there exists embedding of 
$A$-modules of $B$ and $C$ onto $K$, such that if $B$ and $C$ are identified with their image in $K$, $m: B \times C \rightarrow A$ is given by the multiplication of $K$.
\end{lemma}
\begin{pf} 
Since $m:B \otimes C \rightarrow A$ is non-degenerate, $B$ and $C$ have no torsion.

Fix $b_0 \in B-\{0\}, c_0 \in C-\{0\}$ such that $m(b_0,c_0) \neq 0$. Define a morphism of $A$-modules $i: B \rightarrow K$ by setting $i(b)=m(b,c_0)/m(b_0,c_0)$. If $i(b)=0$, then $m(b,c_0)=0$ so $m(b,c_0) b_0 = m(b_0,c_0) b = 0$, and $b=0$ since $B$ has no torsion; thus $i$ is injective. Define $j: C \rightarrow K$ by setting $j(c)=m(b_0,c)$, which embeds $C$ into $K$, and one easily checks that $m(b,c)=i(b)j(c)$.
\end{pf}

\begin{lemma} \label{basicGMA3} 
Assume that $A$ is a domain, with fraction field $K$, and that $R=\mat{A & B \\ C & D}$ is a faithful GMA over $A$, and that $BC \neq 0$. Then $R \otimes_A K$ is isomorphic, as a GMA over $K$, to $M_2(K)$.
\end{lemma}
This follows from the preceding lemma.

\subsection{Topological GMAs}

If $A$ is a topological ring, a {\it topological GMA} is a GMA $R$ over $A$ provided with a topology that makes it a topological $A$-algebra.
More concretely, if $R=\mat{A & B \\ C & A}$ is a GMA, making $R$ a topological GMA amounts to giving a topology on $B$ and $C$ that makes them topological $A$-modules, and make the multiplication $m: B \times C \rightarrow A$ continuous.

For instance, if $A$ is a noetherian local ring which is complete for the topology defined by its maximal ideal, and if $R$ is finite as an 
$A$-module, then $R$ provided with its finite $A$-module topology is a topological GMA. 

We observe that for any topological ring $A$, $R=M_2(A)$ has a unique structure of topological GMA, namely the one given by the product topology on $M_2(A) = A^4$.

\subsection{Pseudo-representations and GMA-valued representations}

Let $A$ be a complete local ring with maximal ideal $\m$ and residue field $\F$. Let $\Pi$ be a group, $(t,d): \Pi \rightarrow A$ a pseudo-representation. 

The reduction $\bar t, \bar d$ modulo $\m$ of $t$, $d$ form a pseudo-representation of dimension $2$ of $G$ with values in $\F$.  We make the following definition:
\begin{definition} \label{multfree} We say that $(t,d)$ is {\it residually multiplicity-free} if there exists a (necessarily unique up to isomorphism) semi-simple representation 
$\rhob: \Pi \rightarrow \GL_2(\F)$ such that $\tr \rhob = \bar t$, $\det \rhob = \bar d$, and this representation 
is the direct sum of distinct absolutely irreducible representations. \end{definition}

By a theorem of Chenevier, there always exists a a finite extension $\F'$ of $\F$ and a $\rhob: \Pi  \rightarrow \GL_2 (\F)$ such that $\tr \rhob = \bar t$, $\det \rho = \bar d$. Hence saying that $(t,d)$ is residually multiplicity-free amounts to saying that we can take $\F'=\F$ and that the residual representation $\rhob$ is  either absolutely irreducible, or the direct sum of two distinct characters. 

Following Mazur \cite[page 246]{mazurModular}, we say that a pro-finite group $\Pi$ satisfies the {\it $p$-finiteness condition} if  for every open subgroup $H$ of $\Pi$, the largest pro-$p$ quotient
$H_p$ of $H$ is topologically of finite type. This condition is always satisfied for a profinite group $\Pi$ which is topologically of finite type, and it is also known to hold for a Galois group $\Pi=\G_{\Q,S}$ where $S$ is a finite set of places ({\it loc. cit.}).

\begin{prop} \label{exGMA} 
Assume that $(t,d)$ is residually multiplicity-free.
\begin{itemize}
\item[(i)] There exists a  faithful GMA $(A,B,C,m,R)$ over $A$, and a morphism of groups
$\rho: \Pi \rightarrow R^\ast$ such that 
\begin{num} \label{trrhotdetrhod} on $\Pi$, $\tr \rho = t$ and $\det \rho=d$, \end{num}
\begin{num} \label{arhoPi} $A \rho(\Pi) = R$. \end{num}
\item[(ii)] If $(\rho,R)$ and $(\rho',R')$  are as in (i), then there exists a unique isomorphism of $A$-algebras $\Psi: R \rightarrow R'$
such that $\Psi \circ \rho = \rho'$.
\item[(iii)] Given an element $g_0 \in \Pi$ such that $\rhob(g_0)$ has two distinct eigenvalues $\lambda_0,\mu_0$ in $\F$ (such an element always exists under \ref{multfree}), there exists a faithful GMA $(A,B,C,m,R)$  over $A$, and a morphism of $A$-algebras $\rho: \Pi \rightarrow R^\ast$ satisfying~\ref{trrhotdetrhod} and~\ref{arhoPi}, and such that 
\begin{num} $\rho(g_0)$ is diagonal and $\rho(g_0) \equiv   \mat{\lambda_0 & 0 \\ 0 & \mu_0} \pmod{\m}$.
\end{num}
\item[(iv)]  If $g_0 \in \Pi$, $(\rho,R)$ and $(\rho',R')$  are as in (iii), the unique isomorphism of $A$-algebras, $\Psi: R \rightarrow R'$ such that $\Psi \circ \rho = \rho'$ (cf. (ii)) is an $A$-isomorphism of GMAs.
\item[(v)]  If $\rhob$ is irreducible, then $R$ is isomorphic to $M_2(A)$ as a GMA over $A$.
If $\rhob$ is reducible, then one has $BC\subset \m$.
\item[(vi)] If $(\rho,R)$ is as in (i),  then $\ker \rho = \ker(t,d)$, and denoting by $\tilde \rho: A[\Pi] \rightarrow R$ the morphism of $A$-algebras extending $\rho$, 
one has $\ker \tilde \rho = \uker (T,D)$.
\item[(vii)] If $A$ is noetherian, if $\Pi$ is a profinite group satisfying the $p$-finiteness condition, and if $(t,d)$ is continuous then for $(\rho,R)$ as in (i), $R$ is of 
finite type as an $A$-module and if $R$  is given its unique topology of $A$-algebras, the morphism $\rho: G \rightarrow R^\ast$ is continuous.
\end{itemize}
 \end{prop}
 \begin{pf} 
 Let $(T,D)$ be the pseudo-representation of $A[\Pi]$ with values in $A$ extending $(t,d)$, as in Lemma~\ref{tdTD}.
Let $R$ be the quotient of $A[\Pi]$ by $\uker(T,D)$, let $\tilde \rho$ be the natural projection $\tilde \rho: A[\Pi] \rightarrow R$ and let $\rho$ be the restriction of $\tilde \rho$ to $\Pi$. 
Let $g_0$ be an element  of $\Pi$ as in (iii), let $\Pi_0$ be the subgroup generated by $g_0$ in $\Pi$ and  let $R_0 \subset R$ be the $A$-subalgebra $A\rho(\Pi_0)$. As proved in \cite[\S1.4]{BC}, the algebras $R$ and $R_0$ are integral over $A$.
 By \cite[\S1.4]{BC} and the hypothesis made on $g_0$, 
if $J_0$ denotes the Jacobson radical of $R_0$, then there is an isomorphism of $\F$-algebra $R_0/J_0 \simeq \F\rhob(\Pi)$. The algebra $\F \rhob(\Pi_0)$ is isomorphic to $\F \times \F$ and we can fix such an isomorphism that sends $\rhob(g_0)$ to $ \smat{\lambda_0 & 0 \\ 0 & \mu_0}$. The two obvious idempotents $(1,0)$ and $(0,1)$ of $\F \times \F$ can be lifted to idempotents $e_1$ and $e_2$ of $R_0$ such that $e_1 e_2 = 0$, $e_1+e_2=1$. This makes $R_0$ and $R$ GMAs with the properties stated in (i) and (iii). The uniqueness
statement (ii) is clear, since if $(\rho,R)$ is as in (i), $R$ has to be a quotient of $A[\Pi]$ through which $(T,D)$ factors, hence of the form $A[\Pi]/I$ with $I$ a two-sided ideal contained in $\uker (T,D)$, but since $R$ is faithful we must have $I=\uker (T,D)$. The uniqueness statement (iv) is equally easy, since a morphism $\Psi$ as in (iv), which exists and is unique by (ii), preserves the diagonal matrix $\rho(g_0)$ which has diagonal terms that are distinct modulo $\m$, hence preserves
the idempotents $e_1$ and $e_2$ and is a morphism of GMA. Finally (v) in the irreducible case is a well-known result of Rouquier and Nyssen extended by Chenevier (\cite[Theorem 2.22]{chenevier})
to the case of general pseudo-representation, and (v) in the reducible case follows from \cite[Theorem 1.4.4]{BC}.

Let us prove (vi). Since $\tilde \rho: A[\Pi] \rightarrow R$ is surjective, one has $\uker (T,D) = \tilde \rho^{-1} \uker(tr_R,\det_R)$.
Since $R$ is faithful, it follows that $\uker (T,D) = \ker \tilde \rho$. Using Lemma~\ref{keruker}, thus implies that $\ker(t,d)=\ker \rho$.

For (vii), let $A[[\Pi]]$ be the completed group algebra of the pro-finite group $\Pi$. Chenevier proves in \cite[\S4]{chenevierv1} that $t$ and $d$ can be extended into a continuous pseudo-representation $(\tilde T,\tilde D)$ of $A[[\Pi]]$ of dimension 2 with values in $A$. The restriction of $(\tilde T, \tilde D)$ to the sub-algebra $A[\Pi]$ is therefore the pseudo-representation $(T,D)$ 
of $A[\Pi]$ corresponding to $(t,d)$. Form the definition of the linear kernel, one has $\uker (T,D) = \uker(\tilde T, \tilde D) \cap A[G]$. Hence $R=A[\Pi]/\uker(T,D)$ is isomorphic to an $A$-sub-algebra of $A[[G]]/\uker(\tilde T, \tilde D)$. The latter is a finite type $A$-module by \cite[Lemma 4.5]{chenevierv1}. Since $A$ is noetherian, 
$R$ is of finite type as an $A$-module. 

Let us prove now that $\rho$ is continuous. Choose a finite family of elements $g_1,\dots,g_m$ of $\Pi$ such that the $\rho(g_i)$ generate $R$. Consider the map $R \rightarrow A^n$, $x \mapsto \tr(x \rho(g_i))$. Since $R$ is faithful, this map is an injection, and by the elementary properties of the natural topology of finite $A$-modules, it induces an homeomorphism of $R$ onto its image. It therefore suffices to prove that the map $g \mapsto \tr (\rho(g) \rho(g_i))$ is continuous for $i=1,\dots,m$, but this is clear since that map is just $t(g g_i)$. 
\end{pf}

\begin{definition} \label{tdrepresenstationdef} Any representation $\rho: \Pi \rightarrow R^\ast$ satisfying the property (i) of the above proposition
will be called a {\it $(t,d)$-representation}. If in addition $\rho$ satisfies condition (iii), we shall say that $\rho$ is {\it adapted to $(g_0,\lambda_0,\mu_0)$.}
\end{definition}

\begin{remark} Without the assumption of $p$-finiteness on $\Pi$, the  assertion (vii) of the preceding theorem is false.
For a counter-example, let $A=\F_p[\epsilon]$ with $\epsilon^2=1$, $V$ an infinite-dimensional $\F_p$-vector space seen as an $A$-module through the map $A \rightarrow \F_p, \epsilon \mapsto 0$,
$b : V \times V \rightarrow \F$ a non-degenerate $\F$-bilinear form, and $m: V \times V \rightarrow A$ defined as $\epsilon b$. Then $m$ satisfies condition \ref{asso}, hence there is a GMA $(A,V,V,m,R)$ which moreover is faithful. Define $\Pi=R^\ast$, and consider the restriction $(t,d)$ of $(\tr,\det)$ to $\Pi$. This is a  
pseudo-representation of dimension 2, and $A[\Pi]/\uker(T,D)=R$ but $R$ is not finite as an $A$-module.  
\end{remark}

\begin{lemma} \label{imagesubGma} Let $R$ be a GMA over $A$ and $\rho: \Pi \rightarrow R^\ast$ a representation of a group $\Pi$.
Assume that there exists an element $g_0 \in \Pi$ such that $\rho(g_0)$ is diagonal, with diagonal terms distinct modulo $\m$.
Then $A \rho(\Pi)$ is a sub-$A$-GMA of $R$.

Furthermore, if $R=M_2(A)$ and $\rho \mod \m$ is absolutely irreducible, then $A \rho(\Pi)=R$.
\end{lemma}
\begin{pf}
If $\rho(g_0)=\smat{\lambda_0 & 0 \\ 0 & \mu_0}$, then the matrix $e_1:=\smat{1 & 0 \\ 0 & 0}=\frac{\rho(g_0)-\mu_0 \rho(1) }{\lambda_0 - \mu_0}$ belongs to $A\rho(\Pi)$, and similarly the matrix $e_2:=1-e_1 = \smat{0 & 0 \\ 0 & 1}$. Then $e_1 A\rho(G) e_1$ is an $A$-submodule of $A$ that contains $1$, so is $A$, and similarly for $e_2 A \rho(\Pi) e_2$. Define $B':= e_1 A \rho(\Pi) e_2$, a submodule of $B$, and $C':= e_1 A \rho(\Pi) e_2$, a submodule of $C$. Then $A \rho(G)=\mat{A & B' \\ C' & A}$ an $A$-sub-GMA of $R$.

For the {\it furthermore}, suppose by contradiction that either $B'$ or $C'$ is a proper sub-module of $B=C=A$. Then $B'C'$ is a proper ideal of $A$, so is contained in $\m$, which shows that $\tr \rho \pmod{\m}$ is the sum of two characters ($\Pi \mapsto \F, g \mapsto e_1 \rho(g) e_1 \pmod{\m}$ for the first, the same with $e_2$ for the second), contradicting the hypothesis.
\end{pf}

\section{Reminder of representation theory}

\subsection{The classification of representations $\rhob$}

\label{classification}

Let  $\Pi$ be a group,  $\F$ a finite field, $\rhob: \Pi \rightarrow \GL_2(\F)$ a representation which is either absolutely irreducible or the sum of two distinct characters. Let us set $\Gb=\rhob(\Pi) \subset \GL_2(\F)$ and $\Gbb$ the {\it projective image} of $\rhob$, that is the image of $\Gb$ in $\PGL_2(\F)$. The well-known classification of such representations according to their projective image is as follows.

\par \bigskip

{\footnotesize
\begin{tabular}{ | p{2cm} |  p{3cm}			 |  p{0.9cm} |  p{3.8cm} 							| p{2.5cm} |  }
\hline
			Name 	& $\Gbb$ is isomorphic to & Subcase    & Description of $\rhob$  & Description of $\ad^0 \rhob$ \\
\hline \hline

\multirow{2}{*}{Cyclic} & $\Z/n\Z$ & $n=2$ & $\chib \oplus \chib'$, with $\chib^2=\chib'^2$ &  
  $(\chib/\chib')^2 \oplus 1$  \\
\cline{3-5}
& & $n>2$ & $\chib \oplus \chib'$, with $\chib^2 \neq \chib'^2$ & $\chib/\chib' \oplus 1 \oplus \chib'/\chib$ 
 \\
\hline

\multirow{2}{*}{Dihedral} & \multirow{2}{*}{$D_n$} & $n >2$  & 
irreducible, isomorphic to $\Ind_{\Pi_1}^{\Pi} \psi_1$ for a unique index $2$ subgroup $\Pi_1$ of $\Pi$ & $\epsilon_1 \oplus \Ind_{\Pi_1}^\Pi \tau $, with   $\Ind_{\Pi_1}^\Pi \tau $ irreducible\\

 \cline{3-5} 
 
& & $n=2$  &  irreducible,  isomorphic to $\Ind_{\Pi_1}^\Pi \psi_i$ for three index two subgroups $\Pi_1$, $\Pi_2$ and $\Pi_3$ & 
$\epsilon_1 \oplus \epsilon_2  \oplus \epsilon_3$ \\ 
\hline

Large image & $\PGL_2(\F_q)$  or  $\PSL_2(\F_q)$ & &   \multirow{2}{*}{irreducible} & \multirow{2}{*}{irreducible}  \\
\cline{1-2}  Exceptional & $A_4$, $S_4$ or $A_5$   &  & & \\
\hline 
\end{tabular} }

In the table above, $\chi$ and $\chi'$ are two distinct characters of $\Pi$, $\psi_i$ is a non-trivial character of $\Pi_i$ for $i=1,2,3$,
and $\epsilon_i$ is the character of $\Pi$ of kernel $\Pi_i$ for $i=1,2,3$, and $\tau$ is a non-trivial character of $\Pi_1$. The group $D_n$ is the dihedral group of order $2n$.

\subsection{A group cohomology computation}

\begin{prop} \label{H1adj0} In the large image and exceptional case, if $V$ is adjoint representation of the natural representation of $\bar G$,
one has $H^1(\bar G, V)=0$.
\end{prop}
\begin{pf} The representation $V$ of $\Gb$ factors through $\Gbb$. Let $Z$ be the kernel of $\Gb \mapsto \Gbb$. The inflation-restriction exact sequence is
$$0 \rightarrow H^1(\Gbb,V) \rightarrow H^1(\Gb,V) \rightarrow H^1(Z,V)$$
and since $Z$ is of order prime to $p$, and $V$ is of order a power of $p$, the last term is 0. It therefore suffices to prove that 
$H^1(G,\bar V)=0$.

For  $\Gb=\PGL_2(\F_q)$ or $\Gb=\PSL_2(\F_p)$, this follows from Matthias Wendt's answer to question 178025 of mathoverflow.

If $\Gb$ is isomorphic to $A_4$ or $S_4$, then the result is clear if $p \geq 5$. If $p=3$, we argue as follows:
Let $K_4$ be the Klein subgroup of $A_4$. One has an exact sequence $0 \rightarrow H^1(A_4/K_4,V^{K_4}) \rightarrow H^1(A_4, V) \rightarrow H^1(K_4,V)$; since $V$ is still irreducible as a representation of $K_4$, 
$V^{K_4}=0$; Since $K_4$ has order prime to $3$, $H^1(K_4,V)=0$. Hence $H^1(A_4,V)=0$. For $S_4$ we use the sequence
$H^1(S_4/A_4,V^{A_4}) \rightarrow H^1(S_4,V) \rightarrow H^1(A_4,V)$ where the first and last term are $0$.

If $\Gb$ is isomorphic to $A_5$, the result is clear if $p \geq 7$. For $p=5$, $\Gb$ is conjugate to $\PSL_2(\F_5)$, a case which has already been dealt with. For $p=3$, let us consider $A_4$ as the subgroup of $A_5$ fixing one letter, and note that since 
$A_4$ has index 5 which is prime to $|V|=27$, it suffices to prove that $H^1(A_4,V)=0$, which has already being done.
\end{pf}

\section{Pink's Lie theory for GMAs}

\label{sectionpinkpseudo}

\subsection{Assumptions concerning the base ring $A$}

In all this section, $p$ is a prime.  We suppose given
\begin{num} \label{baseringA} a topological ring $A$ which is compact and semi-local.
\end{num}
By definition, $A$ semi-local means that $A$ is a finite product $\prod_{i=1}^r A_i$, where the $A_i$ are local rings.
We provide each of the ring $A_i$ with its quotient topology from the topology of $A$. The $A_i$ are compact rings, and are local.
We shall call $\m_i$ the maximal ideal of $A_i$ and $\F_i=A_i/\m_i$ its residue field. By an abuse of language which hopefully will not induce confusion, we shall also call  $\m_i$ the corresponding maximal ideal in $A$, namely $\prod_{j \neq i} A_j \times \m_i$, so that we can write $A/\m_i=\F_i$, and $(\m_i)$, $i=1,\dots,r$ are the complete list of maximal ideals of $A$.

In general, the compact topology on $A_i$ is not the $\m_i$-adic topology. However:
\begin{lemma} \label{lemmanoetherian} \begin{itemize}
\item[(i)] The topological ring $A$ (and its factors $A_i$) is pro-finite (i.e. the open  co-finite ideals $J$ form a basis of neighborhood of $0$)
\item[(ii)] The fields $\F_i$ are finite and the ideals $\m_i$ are open in $A_i$.
\item[(iii)] Each ring $A_i$ is $\m_i$-adically complete, and its $\m_i$-adic topology is finer that its given topology.
\item[(iv)] One has $A_i$ noetherian if and only if $\m_i^2$ is open in $A_i$. In this case, the $\m_i$-adic topology on $A_i$ coincide with its given topology.
\end{itemize}
\end{lemma}
\begin{pf} Assertion (i) is  \cite[Prop. 5.1.2]{RZ}. If we write $A_i = \projlim A_i/J$ with $J$ running among open cofinite ideals of $A_i$, 
then each $A_i/J$ is local with maximal ideal $\m_i/J$ and residue field $\F_i$. In particular $\F_i$ is finite. Moreover $\m_i = \projlim \m_i/J$:
the inclusion $\m_i \subset \projlim \m_i/J$ is obvious, while if $x \in A_i$ is not in $\m_i$, $x$ is invertible, so its image in any $A_i/J$ is not in $\m_i/J$. Therefore we see that $\m_i$ is closed in $A_i$, and since it is cofinite, it is also open. This proves (ii). For $J$ any open co-finite ideal of $A_i$, 
$A_i/J$ is finite local, hence  Artinian, and there is an $n$ such that $(\m_i/J)^n = 0$ in $A_i/J$, that is $\m_i^n \subset J$ in $A_i$. Hence the family 
$\m_i^n$ is cofinal to the family of co-finite open ideals, and $A_i$ is $\m_i$-adically complete. Therefore, every open set
for the given topology contains an ideal $\m_i^n$ hence is also open for the $\m_i$-adic topology. This proves (iii). Finally, note that $\m_i^2$ is open
if and only if $\m_i/\m_i^2$ is finite, i.e. by Nakayama if and only if $\m_i$ is of finite type, i.e. if and only if $A_i$ is noetherian. In this case, all the $\m_i^n$ are cofinite, hence $A_i$ is compact for the $\m_i$-adic topology. The identity map $A_i \rightarrow A_i$ where the source is given the $\m_i$-adic topology, and the target its original topology, which is continuous by (iii), is therefore closed, hence an homeomorphism. This proves (iv).
\end{pf}

The Jacobson radical $\rad A$ of $A$ will be denoted by $\m$. We have $\m=\prod_{i=1}^r \m_i = \cap_{i=1}^r \m_i$. 
It follows from the lemma that $A$ is $\m$-complete and profinite for the $\m$-adic topology.

From now on and throughout this section, we make the following assumption:
\begin{num} \label{podd} The prime $p$ is odd. \end{num}
Since $p>2$, if $x$ is an element of $1+\m$, there exists by Hensel's lemma a unique $y \in 1+\m$ such that $y^2=x$. We shall henceforth denote that element by $\sqrt{x}$. We observe that the map $x \mapsto \sqrt{x}$ is continuous.

\subsection{A slightly generalized setting for Pink's theory}

\label{subpinksetting}

Pink's theory is concerned with certain closed subgroups of $\GL_2(A)$, the multiplicative group of invertible elements in the matrix algebra $M_2(A)$.
To allow for more generality, we shall consider closed subgroups of the multiplicative group of units of a {\it generalized matrix algebra.}

To fix notation for the rest of this section, 
\begin{num} \label{assumptionR} Let $R=\mat{A & B \\ C & A}$ be a topological GMA over $A$, which is compact and Cayley-Hamilton (cf. \S\ref{subGMA}). 
\end{num}

We denote by $R^\ast$ the multiplicative group of invertible elements in $R$. Clearly, it is also the set of elements $r$ of $R$ such that $\det r \in A^\ast$. It follows that $R^\ast$ is both open and closed in $R$, and, provided with the subspace topology, is a compact topological group. We denote by $SR^\ast$ the set of elements in $R^\ast$ with determinant $1$. Obviously this is a closed normal subgroup of $R^\ast$.

We shall denote by $\rad R$ the Jacobson radical of the algebra $R$. It is a closed hence compact additive subgroup of $R$.  We shall denote by $R^1$ the subgroup $1+ \rad R$. It is a closed normal subgroup of $R^\ast$.

We call $SR^1$ the intersection of $SR$ and $R^1$ in $R$. 
Obviously $SR^1$ is a closed normal subgroup of $R^\ast$.

\begin{remark} \label{descriptionradical}
To fix ideas, we shall now give an explicit description of the various rings and groups introduced above, in the case where $A$ is local. In this case there are two possibilities regarding the ideal $BC=m(B,C)$ of the ring $A$. Either $BC=A$, or $BC \subset \m$.

When $BC=A$, then by Lemma~\ref{basicGMA1}, $R$ is isomorphic as GMA to $M_2(A)$, so we can as well assume that $R = M_2(A)$ as a topological GMA. Its radical $\rad R$ is $M_2(\m)=\m M_2(A)$ and the quotient $R/\rad R$ is the simple algebra $M_2(\F)$.
The group $R^1$ is the multiplicative group of matrices in $M_2(A)$ which are congruent to $\Id$ modulo $\m M_2(A)$. The group $SR^1$ is the subgroup of those whose determinant is $1$. Note that in the literature, those groups $R^1$ and $SR^1$ 
are often denoted $\GL_2^1(A)$ and $\SL_2^1(A)$ respectively. In this case we are in the situation considered by Pink.

When $BC \subset \m$, the radical $\rad R$ is $\mat{\m & B \\ C & \m}$ and the quotient $R/\rad R$ is the semi-simple algebra of diagonal matrices $\mat{ \F & 0 \\ 0 & \F}$. The group $SR^1$ is the group of matrices $\mat{a & b \\ c & d}$ in $R$ such that $a \equiv d \equiv 1 \pmod{\m}$ and $ad-bc=1$. 

In the general case, if $A$ is a finite product of local rings $A_i$, then $R$ naturally decomposes as a product of GMA $R_i$ and the radical $\rad R$ as a product of $\rad R_i$, for each of which one of the two description above holds. 
\end{remark}

\begin{lemma} \label{lemmadescriptionradical}
If $m \in \rad R$, $\tr m, \tr m^2, \det m \in \m$.
\end{lemma}
\begin{pf} We may assume that $A$ is local, in which case we use the description of $R$ and $\rad R$ given in the preceding remark.
 If  $R=M_2(A)$, $m \in M_2(\m)$ and the result is clear.
If $BC \subset A$, then $m \in \mat{\m & B \\ C & \m}$ so $\tr(m) \in \m$ and $\tr(m^2) \in \m^2 + BC \subset \m$,
and finally $\det(m)=(\tr(m)^2-\tr(m^2))/2  \in \m$.
\end{pf}

\noindent {\large \bf Notation:} In the rest of this paper, we shall use freely the following notation: if $S$ is a set of matrices, $S^0$ is the set of matrices of trace zero in $S$.
 If $X,Y$ are two closed additive subgroups of $R$, we shall denote by $[X,Y]$ (resp. $X\cdot Y$ or $XY$) the {\it closure} of the subgroup generated by all commutators $[x,y]$ (resp. $xy$) for $x \in X, \ y \in Y$.

\begin{remark} We observe that $(\rad R)^0$, provided with the Lie bracket   $[r,r']=rr'-r'r$, is a Lie algebra over $A$.
Concretely, $(\rad R)^0=\mato{\m & \m \\ \m & \m}$ when $R=M_2(A)$ and $(\rad R)^0=\mato{\m & B \\ C & \m}$ when $BC \subset \m$. \end{remark}

\subsection{Pink's Theta map}

Following Pink, let $\Theta$ be the continuous $A$-linear map $R \rightarrow R$, $r \mapsto r - \frac{\tr r}{2} \Id$. Pink states eleven formulas involving $\Theta$ and $\tr$. We state the analog in our more general situation of the formulas we need:

\begin{num} \label{pinkf2} If $x,y \in R$, $[\Theta(x),\Theta(y)]=\Theta(xy)-\Theta(yx).$ \end{num}
\begin{num} \label{pinkf9} If $x \in SR$, $y \in R$, one has $\tr(x) \Theta(y) =  \Theta(xy)+\Theta(x^{-1}y)$. \end{num}
\begin{num} \label{pinkf3} If $x,y \in R$, one has $2 \Theta(xy) =  [\Theta(x),\Theta(y)] +\tr(x) \Theta(y)+\tr(y) \Theta(x)$. \end{num}
\begin{num} \label{pinkf4} If $x,y \in R$, $\tr(\Theta(x) \Theta(y) )= \tr(xy) - \tr(x) \tr(y)/2$. \end{num}
\begin{num} \label{pinkf8} If $x \in SR$,  one has $ \Theta(x^{-1}) = -\Theta(x)$. \end{num}
\begin{num} \label{pinkf10} If $x,y,u,v \in (\rad R)^0$,  one has $ 4 \tr(xy) [u,v] = [y, [x, [u,v]]] + [x,[y,[u,v]]] +  [[x,v],[y,u]] + [[y,v],[x,u]]$. \end{num}

These formulas are proved by easy computations left to the reader, using the facts that in $R$, $\tr(xy)=\tr(yx)$ and that for any $x \in R$, the Cayley-Hamilton identity holds, namely $x^2-\tr(x) x + \det(x)=0$, with $\det(x)=(\tr(x)^2-\tr(x^2))/2$. (Also useful is the formula
$xy-\tr(x) y -\tr(y) x + \tr(x)\tr(y)-\tr(xy) = 0$ for $x,y \in R$, which is obtained by bi-linearizing the Cayley-Hamilton identity).

 Alternatively, we can use Proposition 1.3.13 of \cite{BC} which implies that every Cayley-Hamilton GMA $R$ can be embedded in a trace-preserving way into $M_2(A')$ for $A'$ some commutative ring containing $A$. This reduces the formulas to prove to the case of $M_2(A')$. 
 In this case these formulas are stated in \cite{pink}, though their proofs are also left to the reader.

\begin{lemma} The map $\Theta$ induces a homeomorphism from $\SR^1$ onto $(\rad R)^0$. Its inverse is given by
\begin{num} \label{Thetainvan} $ \Theta^{-1} \smat{a & b \\ c & -a }= \smat{a + \sqrt{1+bc+a^2} & b \\ c & -a+ \sqrt{1+bc+a^2}}$ \end{num}
\noindent
or equivalently 
\begin{num} \label{Thetainvan2} $\Theta^{-1} m = m + \sqrt{1+\tr(m^2)/2} \Id$. \end{num}
\noindent
Moreover one has 
\begin{num} \label{Thetainvan3} $\tr ( \Theta^{-1} m) = 2 + \sum_{n \geq 1} 2^{1-n} {n \choose 1/2} \tr(m^2)^n$.\end{num}
\end{lemma}
\begin{pf} It is clear that $\Theta$ sends $\SR^1$ into $(\rad R)^0$.  If $m$ in $(\rad R)^0$, $x \in SR^1$ and 
$\Theta(x) = m$ then one has $x = m + \lambda \Id$ for some $\lambda \in 1+\m$ and using that $\det x=1$, one gets $\lambda^2 = 1+\tr(m^2)/2$. Since $\tr m^2 \in \m$ by Lemma~\ref{lemmadescriptionradical}, this equation defines a unique $\lambda$, which shows that for every $m \in (\rad R)^0$, there exists a unique $x$ such that $\Theta(x)=m$, and proves the formula for $\Theta^{-1}$. Formula~\ref{Thetainvan3} follows using Newton's Taylor expansion for $\sqrt{1+t}$.

\end{pf}

\subsection{The Lie algebra $L$ attached to a subgroup of $SR^1$} 
 
\label{subpinktheory}

The object of Pink's theory is to understand the structure of the closed subgroups of $SL_2^1(A)$, using Lie-theoretic methods.  
Our objective here is to expand Pink's constructions and results to the case of subgroups of $SR^1$, where $R$ is a GMA over $A$ as above.
We shall offer from this sub-section \S\ref{subpinktheory} to \S\ref{subdescending} a self-contained presentation, where arguments, whose details follow closely those of \cite{pink} are re-organized and somewhat simplified.

Let $\Gamma$ be a closed subgroup of $SR^1$. Following Pink we define 
a closed subgroup $L$ of $(\rad R)^0$ as the closure of the additive subgroup of $(\rad R)^0$ generated by $\theta(\Gamma)$.

Obviously, $\Gamma \subset \Theta^{-1}(L)$ but we may not have equality. Observe that the subgroup $L$ is not in general an $A$-submodule of $(\rad R)^0$. 

\begin{theorem}[Pink] \label{thmLlie} One has $[L,L] \subset L$, that is $L$ is a Lie subring of $(\rad R)^0$.
\end{theorem}
\begin{pf} It suffices to show that if $x,y \in \Gamma$, $[\Theta(x),\Theta(y)] \in L$, that is $\Theta(xy) - \Theta(yx) \in L$ by~\ref{pinkf2}.
Since $xy$ and $yx$ are in $\Gamma$, this is clear. 
\end{pf}

\begin{definition} We call $L=L(\Gamma)$ the {\it Pink's Lie algebra} of $\Gamma$. \end{definition}

\begin{lemma}[Pink] For $\gamma \in \Gamma$, one has $\tr(\gamma) L \subset L$.
\end{lemma}
\begin{pf} This follows immediately from \ref{pinkf9}. \end{pf}

\subsection{The pseudo-ring $P$ attached to a closed subgroup $\Gamma$ of $SR^1$.}

For $\Gamma$ and $L$ as in the preceding section, we define $$P=P(\Gamma)=\tr(L^2).$$ 
This is a closed additive subgroup of $A$. (Note that our $P$ is denoted by $C$ in \cite{pink}).
\begin{theorem}[Pink] \label{thmpinkdirect} One has $PL \subset L$.
\end{theorem}
\begin{pf} By definition, $P$ is the closure of the additive subgroup generated by the $\tr(\Theta(x) \Theta(y))$ for $x,y \in \Gamma$.
By \ref{pinkf4}, one has $\tr(\Theta(x) \Theta(y)) = \tr(xy) - \tr(x) \tr(y)/2 \in \tr(\Gamma) + \tr(\Gamma)^2$. Thus $P \subset \tr(\Gamma) + \tr(\Gamma)^2$, and the theorem follows from the preceding lemma.
\end{pf}

\begin{cor} \label{corPpseudo} The subgroup $P$ of $A$ is stable by multiplication; in other words, it is a pseudo-subring of $A$. Moreover
$P$ is the smallest closed pseudo-subring of $A$ containing $\tr(\gamma)-2$ for all $\gamma \in \Gamma$.
\end{cor}
\begin{pf} \label{corPLL} Since $P L \subset L$, one has $P^2 = P \tr(L\cdot L) = \tr(P L \cdot L) \subset \tr(L \cdot L) = P$, hence $P$ is a pseudo-subring. Let us call by $Q$ the subgroup of $A$ generated by $\tr(\gamma)-2$, $\gamma \in \Gamma$.
Let us first show that $Q \subset P$. If $\gamma \in \Gamma$, $m=\tr(\gamma)$, one has $\tr(\gamma)= -2 + \sum_{n \geq 1} 2^{1-n} {n \choose 1/2} \tr(m^2)^n$ by \ref{Thetainvan3}. Since $\tr(m^2) \in P$ and $P$ is stable by multiplication, $\tr(m^2)^n \in P$ for all $n$ and since $P$ is closed, $Q \subset P$. On the other hand, as seen in the proof of the preceding theorem, $P$ is the closed subgroup of $A$ generated by the elements $\tr(xy)-\tr(x) \tr(y)/2$ for $x,y \in \Gamma$, that is by the elements $\tr(xy)-2 - (\tr(x)-2)(\tr(y)-2)/2 - (\tr(x)-2) - (\tr(y)-2) \in Q+Q^2$ 
Thus $Q \subset P \subset Q+Q^2$, and since $P$ is a closed pseudo-ring, it follows that the closed pseudo-subring of $A$ generated by $Q$ 
is $P$. \end{pf}

\subsection{Pink's converse theorem}

\begin{theorem}[Pink] \label{pinkconverse} Let $L$ be a Lie subring of $(\rad R)^0$. Set $P=\tr(L \cdot L)$. 
If $P L \subset L$, then $H := \Theta^{-1}(L)$ is a closed subgroup of $SR^1$, and $\Theta$ is a homeomorphism of $H$ onto $L$. In particular $L = L(H)$, and $P= P(H)$.
\end{theorem}
\begin{pf} If $P L \subset L$, then one sees as in the proof of Cor.~\ref{corPLL} that 
$P$ is a pseudo-subring and $\tr(h)-2 \subset P$ for every $h \in H$. Thus $\tr(H) L \subset L$
 
If $x,y \in H$, $2 \Theta(xy) = [\Theta(x),\Theta(y)] + \tr(x) \Theta(y) + \tr(y) \Theta(y)$ by \ref{pinkf3}. The first term is in $L$ because $L$ is a Lie subring, the last two are also in $L$ since $\tr(H) L \subset L$. Therefore, $xy \in H$. Also by \ref{pinkf8}, $\Theta(x^{-1})=-\Theta(x)$ so $x^{-1} \in H$. This shows that $H$ is a subgroup of $SR^1$, obviously closed.
\end{pf}

\subsection{Descending the central sequence}

\label{subdescending}
Let $\Gamma$ be a closed subgroup of $SR^1$, $L=L(\Gamma)$ its Pink's Lie algebra, $P=P(\Gamma)=\tr(L \cdot L)$ the attached pseudo-ring. We define: 
\begin{itemize}
\item for $n \geq 1$, closed Lie subrings $L_n$ of $(\rad R)^0$, defined by recurrence as follows: $L_1=L$, $L_{n+1}=[L_{n},L]$;
\item for $n \geq 1$, closed subsets $H_n = \Theta^{-1} L_n$ of $\SR^1$.
\item for $n \geq 1$, closed subgroups $\Gamma_n$ of $\SR^1$ defined by recurrence as follows: $\Gamma_1=\Gamma$, $\Gamma_{n+1}=(\Gamma_{n},\Gamma)$ (closed commutators subgroup) for $n \geq 1$
\end{itemize}

\begin{prop}[Pink] Let $n,m \geq 1$.
\begin{itemize}
\item[(i)]  If $n \geq 1$, $L_{n+1} \subset L_{n}$.
\item[(ii)] If $n,m \geq 1$, $[L_n,L_m] \subset L_{n+m}$.
\item[(iii)] If $n \geq 1$, for $h \in H_n$, $\tr(h)-2 \in P$.
\item[(iv)] If $n \geq 1$, $H_n$ is a closed subgroup of $SR^1$ and $\Theta: H_n \rightarrow L_n$ is an homeomorphism.
\item[(v)] If $n \geq 2$, $P L_n \subset L_{n+2}$.
\item[(vi)] If $n \geq 2$, $\Theta$ induces a bicontinuous isomorphism of groups $H_n/H_{n+1} \simeq   L_n/L_{n+1}$.

\end{itemize} 
\end{prop}
\begin{pf} Assertions (i) and (ii) follows easily by induction from Theorem~\ref{thmLlie}. 
For (iii), write $m=\Theta(h) \in L_n \subset L$. Then by \ref{Thetainvan3}, $\tr(h)-2 = \sum_{n \geq 1} 2^{1-n} {n \choose 1/2} \tr(m^2)^n \in P$.

For (iv), from $P L \subset L$ one proves by induction that $P L_n \subset L_n$. One therefore has $\tr(L_n \cdot L_n) L_n \subset \tr(L \cdot L) L_n = P L_n \subset L_n$. It then follows from Theorem~\ref{pinkconverse} applied to the Lie subring $L_n$ that $H_n = \Theta^{-1}(L_n)$ is a closed subgroup of $L_n$, and that $\Theta$ is a homeomorphism of $H_n$ onto $L_n$. 

Formula (v) follows from~\ref{pinkf10} for $n=2$ and then by induction for all $n \geq 2$.

For (vi), if $x,y \in H_n$, then by \ref{pinkf3}, $$\Theta(xy)-\Theta(x) -\Theta(y) = \frac{1}{2}([\Theta(x),\Theta(y)] - (\tr(x)-2)\Theta(y) - (\tr(y)-2)\Theta(x)).$$
Hence $\Theta(xy) - \Theta(x) - \Theta(y) \in L_{2n} + L_{n+2} \subset L_{n+1}$ by (i), (ii), (iii) and (v) (applicable since $n \geq 2$).
This shows that $\Theta$ induces a group morphism from $H_n$ to $L_n/L_{n+1}$. This morphism is surjective by (iv), and 
the kernel of this morphism is clearly $H_{n+1}$, hence (vi).
\end{pf}

The most important theorem of Pink's theory is Theorem~\ref{thmpink1} below,  which shows that for $n \geq 2$, the terms $\Gamma_n$ of the descending central sequence
of $\Gamma$ are determined by their Lie algebra $L_n$, hence by $L$.

First, we need a lemma: 
\begin{lemma}[Pink] Let $n \geq 2$. For $x \in H_1, y \in H_{n-1}$ one has
$$\Theta(xyx^{-1}y^{-1}) \equiv [\Theta(x),\Theta(y)] \pmod{L_{n+1}}.$$
In particular, $xyx^{-1}y^{-1} \in H_n$.
\end{lemma}
\begin{pf}
One writes $$2\Theta(xyx^{-1} y^{-1}) = 2\Theta( [x,y] x^{-1} y^{-1}) = [\Theta([x,y]),\Theta(x^{-1}y^{-1})] + \tr(x^{-1} y^{-1})  \Theta([x,y])$$ by ~\ref{pinkf3}. Since obviously $[\Theta(x),\Theta(y)]=[x,y] = \Theta([x,y])$, this can be written
$$2\Theta(xyx^{-1} y^{-1})  =  [[\Theta(x),\Theta(y)],\Theta(x^{-1}y^{-1})] + \tr(x^{-1} y^{-1})[\Theta(x),\Theta(y)]$$
Now $\Theta(x) \in L_1$, $\Theta(y) \in L_{n-1}$, so $[\Theta(x), \Theta(y)] \in L_n$ and $\Theta(x^{-1} y^{-1}) \in L_1$.
Thus, the first term of the RHS is in $L_{n+1}$. As for the second term, $\tr(x^{-1} y^{-1}) -2 \in P$, and since $P L_n \subset L_{n+2}$, one gets that the second term is $2 [\Theta(x),\Theta(y)] \pmod{L_{n+1}}$ and the lemma follows.
\end{pf} 
\begin{theorem}[Pink] \label{thmpink1} For $n\geq 2$, one has $ \Gamma_n = H_n = \Theta^{-1}(L_n)$. Hence $\Theta$ realizes an homeomorphism of $\Gamma_n$ on $L_n$ for $n \geq 2$.
\end{theorem}
\begin{pf} We follow approximately Pink's method.

By definition $\Gamma_1=\Gamma \subset H_1$. We prove by induction that $\Gamma_n \subset H_n$ for all $n$. Assuming $\Gamma_{n-1} \subset H_{n-1}$, we get for $x \in \Gamma$, $y \in \Gamma_{n-1}$, $\Theta(x) \in L_1$, $\Theta(y) \in L_{n-1}$, hence by the commutator relation $\Theta(xyx^{-1}y^{-1}) \in [L,L_{n-1}] + L_{n+1} \subset L_n$, and
$xyx^{-1} y^{-1} \in H_n$. Since $H_n$ is a closed subgroup of $SR^1$, and $\Gamma_n$ is the closed subgroup generates by the $xyx^{-1}y^{-1}$ as above, one gets $\Gamma_n \subset H_n$.

Let $\Delta_n$ be the closed subgroup of $(\rad R)^0$ generated by $\Theta(\Gamma_n)$. We claim by induction that $\Delta_n + L_{n+1}= L_n$ for all $n \geq 1$. This is true for  $n=1$ because by definition $\Delta_1=L_1$. For $n \geq 2$, since $\Gamma_n$ is the subgroup generated by $xyx^{-1}y^{-1}$ for $x \in \Gamma$, $y \in \Gamma_{n-1}$, and $\Theta$ is a morphism from $\Gamma_n$ to $L_n/L_{n+1}$, $\Delta_n + L_{n+1}$ is the closed subgroup of $(\rad R)^0$ generated by $L_{n+1}$ and the elements $\Theta(xyx^{-1}y^{-1})$, that is, by the lemma, the elements $[\Theta(x),\Theta(y)]$. Since the closed subgroups generated by those elements is $[L_1,L_{n-1}]=L_n$, we get that $\Delta_n + L_{n+1} = L_n$.

For $n \geq 2$, since $\Theta$ is a morphism from $H_n$ onto $L_n/L_{n+1}$, $\Theta(\Gamma_n) + L_{n+1}$ is already a closed subgroup of $L_n$, hence it is $\Delta_n + L_{n+1} = L_n$. We thus have shown, for all $n \geq 2$
$$\Theta(\Gamma_n) + L_{n+1} = L_n.$$
Applying this formula for $n$ replaced by $n+1$ gives a description of $L_{n+1}$ that we can plug in the LHS of the formula, getting
$\Theta(\Gamma_n)+ L_{n+2} = L_n$, and by induction on $m$, $\Theta(\Gamma_n) + L_{n+m} = L_n$ for all $m \geq 1$.
Since $\cap_m L_{n+m} = 0$ and $\Theta(\Gamma_n)$ is closed, one gets $\Theta(\Gamma_n) = L_n$, hence $\Gamma_n=H_n$ and the theorem. 
\end{pf}

Thus, the knowledge of the Lie algebra $L$ of $\Gamma$ determines the derived subgroup $\Gamma_2$ of $\Gamma$. There is an other result of Pink, limiting the possibilities for the quotient $\Gamma/\Gamma_2$:

\begin{theorem}[Pink] \label{thmpink2} The composition law $\ast$ on $L/L_2$ defined by
$$x \ast y = x (\sqrt{1+\tr(y^2)/2}) + y(\sqrt{1+\tr(x^2)/2})$$
makes $L/L_2$ a commutative group. The map $\Theta$ induces a bicontinuous morphism of groups $H_1/H_2 \rightarrow (L/L_2,\ast)$.
The image $\Delta$ of $\Gamma/H_2=\Gamma/\Gamma_2$ in $L/L_2$, which is obviously a subgroup of $L/L_2$ for the law $\ast$, 
topologically generates $L/L_2$ for the law $+$.
\end{theorem}
Since we shall only use this theorem  in the case where $R=M_2(A)$ (see Prop.~\ref{exampleliesimple}), we just refer to 
\cite[Prop. 2.6]{pink} for the proof.

\subsection{Complements to Pink's theory}

\subsubsection{Functoriality w.r.t. surjective morphism of rings}

\label{functorialitypink}

Let $J$ be an ideal of $A$. The ring $A/J$ is still a compact semi-local topological ring, of radical $\m/(\m \cap J)$, with residue fields a subset of the set of residue fields of $A$, hence all finite of characteristic $p>2$. In other words, $A/J$ satisfies \ref{baseringA} and \ref{podd}.

The $A/J$-algebra $R_J = R/JR = \mat{A/J & B/JB \\ C/JC & A/J}$ is a GMA which is obviously of finite type as an $A/J$-module,
and also Cayley-Hamilton. We denote by $\pi_J$ the surjective morphism of algebras $R \rightarrow R/JR$. This morphism induces a morphism of multiplicative groups $\pi_J: R^\ast \rightarrow R_J^\ast$ which is still surjective because an element of a GMA is invertible if and only if its determinant is.
 It also induces a surjection $R^1 \rightarrow R_J^1$ and a morphism $SR^1 \rightarrow SR_J^1$, which we again denote by $\pi_J$. Also $\pi_J$ induces a map $\pi_J: (\rad R)^0 \rightarrow (\rad R_J)^0$.

If $\Gamma$ is a closed subgroup of $SR^1$, let us denote by $\Gamma_J$ the closed subgroup $\pi_J(\Gamma)$. Then we can apply Pink's theory
to $\Gamma_J$ and define sub-Lie-algebras $L_n(\Gamma_J)$ of $(\rad R_J)^0$. The functoriality mentioned in the title is the fact that
\begin{num} \label{functo1}  for all $n \geq 1$, $\pi_J(L_n(\Gamma)) = L_n(\Gamma_J)$. \end{num}
\noindent This is easy to see for $n=1$ from the definition for $L_1$, and then by induction on $n$ for any $n$.

\subsubsection{Multiplication by $\tr(\Gamma)$}

\begin{lemma} \label{multtrgamma} For every $\gamma \in \Gamma$, and every $n \geq 1$, one has $\tr(\gamma) L_n = L_n$. \end{lemma}
\begin{pf}
It suffices to prove the first assertion for $n=1$, because then, one has $L_{n+1}=[L_1,L_n]=[\tr(\gamma) L_1, L_n] = \tr(\gamma) [L_1,L_n] = \tr(\gamma) L_{n+1}]$. For $n=1$ we already know that $\tr(\gamma)L \subset L$, so we just need to show that $\tr(\gamma)^{-1} L \subset L$.

Note that $\tr(\gamma)\equiv 2 \pmod{\m}$. Let $m = \Theta(\gamma)$. Then $\gamma = \Theta^{-1}(m)$
so that by \ref{Thetainvan3}, $\tr \gamma = 2+\sum_{n \geq 1} 2^{1-n} {n \choose 1/2} \tr(m^2)^n$ and $\tr(\gamma)^{-1} = 2^{-1} + \sum_{n \geq 1}  b_n \tr(m^2)^n$ for some coefficients $b_n \in \F_p$ that we need not compute. Since $\tr(m^2) \in P(\Gamma)$, $\tr(m^2)^n L \subset L$ hence $\tr(\gamma)^{-1} L \subset L$ which completes the proof of the first assertion. \end{pf}


\subsubsection{A simple class of examples}

Let $I$ be a closed pseudo-subring of $A$ contained in $\m$, that is a closed additive subgroup of $\m$, stable by multiplication. Let $R=M_2(A)$ be the standard GMA. Then 
$L=\mato{ I & I \\ I & I}$ is a $\Z_p$-Lie sub-algebra of $\mato{ \m & \m \\ \m & \m} = (\rad R)^0$.
We will determine the closed subgroups $\Gamma$ of $\SL_2^1(R)$ that have $L$ as Pink's Lie algebra; actually there is only one such subgroup:

\begin{prop} \label{exampleliesimple} 
Let $\Gamma$ be a closed subgroup of $\SL_2^1(R)$ such that $L(\Gamma)=L=\mato{ I & I \\ I & I}$. Then $\Gamma = \Theta^{-1}(L)$
and $\Theta$ realizes an homeomorphism from $\Gamma$ onto $L$. More generally $\Gamma_n = \Theta^{-1}\left(\mato{ I^n & I^n \\ I^n & I^n} \right)$ for every $n \geq 1$.
\end{prop}  
\begin{pf} If we set $X=\smat{0 & 1 \\ 0 & 0}$, $Y = \smat{0 & 0 \\ -1 & 0}$, $J=\smat{1 & 0 \\ 0 & -1}$, then the usual commutation relations
are $[qX,q'Y]=qq'J$, $[qJ,q'X]=2qq'J$ and $[qJ,q'Y]=-2qq'X$,  for any $q,q' \in I$. The additive subgroup generated by these elements, $L_2=[L,L]$, is thus $\mato{ I^2 & I^2 \\ I^2 & I^2}$. Similarly one proves by induction that $L_n=\mato{ I^n & I^n \\ I^n & I^n}$ for any $n$.

For $x \in L$, the power series defining $\sqrt{1+\tr{x^2}/2}-1$ has all its terms in $I^2$, hence is in $I$ since $I$ is closed under multiplication and topologically.
Thus for $x,y \in L$, $y \sqrt{1+\tr{x^2}}  - y \in I L \subset L_2$ and it follows that $x \ast y \equiv x+ y\pmod{L_2}$ (using the notation of Theorem~\ref{thmpink2}.).
The subgroup $\Theta(\Gamma) \pmod{L_2}$ of $(L/L_2,\ast)$ is thus also a subgroup for the additive law $+$, and therefore, by Theorem~\ref{thmpink2}, is such that its topological closure is $L/L_2$. Since it is already closed, one has $\Theta(\Gamma) \equiv L \pmod{L_2}$.
Since $\Theta(\Gamma)$ contains $L_2$, we obtain 
$\Theta(\Gamma)=L$. The proposition easily follows.
\end{pf}

\subsubsection{Haar measures}

For any compact group $\Delta$, we denote by $\mu_{\Delta}$ the Haar measure on $\Delta$ normalized so as to have a total mass $1$.

\begin{lemma} \label{haarmeasure} Let $H$ and $H'$ be two compact groups, $(H_n)_{n \geq n_0}$ (resp. $(H'_n)_{n \geq n_0}$) a decreasing sequence of closed normal subgroups in $H$ (resp. in $H'$) such that $H_{n_0}=H$ and $\cap_n H_n = \{1\}$ (resp. $H'_{n_0} = H'$ and $\cap_n H'_n = \{1\}$). Let $f$ be an homeomorphism from $H$ to $H'$ (not necessarily a group homomorphism) such that for every $h$ in 
$H$, $f(h H_n) = f(h) H'_n$. We assume that
\begin{itemize} \item[(i)] either the induced map $\bar f_n: H_n/H_{n+1} \rightarrow H'_n/H'_{n+1}$ is a morphism of groups,
\item[(ii)] or the $H_n$ are open in $H$.
\end{itemize}
Then  $f$ sends the Haar measure $\mu_H$ to the Haar measure $\mu_{H'}$.
\end{lemma}
\begin{pf} By assumption, $\bar f_n: H_n/H_{n+1} \rightarrow H'_n/H'_{n+1}$ is either an isomorphism of groups, or a bijection between finite groups, hence in both cases sends the normalized Haar
measure of $H_n/H_{n+1}$ on the normalized Haar measure of $H'_n/H'_{n+1}$. Using this, and an induction over $n$ and Fubini, one sees
that the map $\bar f: H/H_n \rightarrow H'/H'_n$ preserves Haar measures.

To prove the lemma, it suffices to prove that $\mu_H(U) = \mu_{H'}(f(U))$ for any open set $U$ in $H$. Since $H$ is compact, $U$ contains $H_n$ for some $n$, and $f$ induces
a bijection $\bar f$ from the finite group $H/H_n$ to the finite group $H'/H'_n$. If $\bar U$ is the image of $U$ in $H/H_n$, we are reduced to prove that
$\mu_{H/H_n}(U) = \mu_{H'/H_n}(\bar f(\bar U))$, which we have already done.\end{pf}

\begin{prop} \label{ThetapreservesHaar} In the situation of Theorem~\ref{thmpink1}, the homeomorphism $\Theta: \Gamma_2 \rightarrow L_2$ sends the Haar measure $\mu_{\Gamma_2}$ to the Haar measure $\mu_{L_2}.$
\end{prop}
\begin{pf} We apply the preceding lemma to $f=\Theta$, $H=\Gamma_2$, $H'=L_2$, $n_0=2$, $H_n=\Gamma_n$, $H'_n=L_n$.
\end{pf}   

Let us note for later use another application of Lemma~\ref{haarmeasure}.
\begin{lemma} \label{lemmahomeo} Let $V$ be a closed additive subgroup of $R$, $\sigma: V \rightarrow V$ a map satisfying the following property:
\begin{eqnarray*} \forall v,v' \in V, n \in \N, \ v-v' \in \m^n R \Longrightarrow \sigma(v)-\sigma(v') \in \m^{n+1} R \end{eqnarray*}
Let $\Psi: V \rightarrow V$ be the map $\Psi(v) = v + \sigma(v)$. Then $\Psi$ is an homeomorphism of $V$ onto $V$ and sends the Haar measure $\mu_V$ to itself.
\end{lemma}
\begin{pf} If $v \neq v' \in V$, let $n$ be an integer such that $v-v' \in \m^n R$ but $v-v' \not \in \m^{n+1}R$. Then $\Psi(v)-\Psi(v')=(v-v') + (\sigma(v)-\sigma(v'))$ and since $\sigma(v)-\sigma(v') \in \m^{n+1} R$, $\Psi(v)-\Psi(v')$ is not in $\m^{n+1}R$ and in particular $\Psi(v)  \neq \Psi(v')$. Hence $\Psi$ is injective. If $v' \in V$, consider the map
$h: V \rightarrow V, y \mapsto v'-\sigma(y)$. The hypothesis made on $\sigma$ implies that this map has a fixed point in $V$, so there exists $v$ such that $v'-\sigma(v)=v$, or $\Psi(v)=v'$. Hence $\Psi$ is surjective. As $\Psi$ is obviously continuous, and closed since $V$ is compact, it is a homeomorphism. To show that $\Psi$ preserves the Haar measure, we apply Lemma~\ref{haarmeasure} with $H=H'=V$, $H_n=H'_n =V \cap (\m^n R)$: for any $n$, the group $H_n$ is open in $V$ since $\m^n R$ is open
in $R$ and the hypothesis implies that $\Psi(v+H_n)=\Psi(v)+H_n$.
\end{pf}

\subsection{Decomposition of Lie algebras}

In this subsection, $R$ is a GMA over $A$ satisfying the conditions of \S\ref{subpinksetting}.

\subsubsection{Decomposable Lie algebras}

\label{decomposable}

Let $L$ be a closed subspace of $(\rad R)^0$.

\begin{num} \label{Ldiag} We shall say that $L$ is {\it decomposable} if, for any $\mat{a & b \\ c & -a} \in L$, one has
$\mat{a & 0 \\ 0 & -a} \in L$ and $\mat{ 0 & b \\ c & 0} \in L$. 
 \end{num}

 We shall denote by $\Delta$ and $\nabla$ the additive groups of diagonal matrices and anti-diagonal matrices in $L$.  Thus, $L$ is decomposable if and only if 
\begin{num} $L = \Delta \oplus \nabla.$\end{num}
\noindent
Since by definition matrices in $L$ have trace $0$ and diagonal terms are in the radical $\m$ of $A$, we see that $\Delta$ has the form \begin{num} \label{defI1} $\Delta = I_1 J$, with $I_1$ a unique additive closed subgroup of $\m$, \end{num}
\noindent where $J$ denotes as usual the matrix $\smat{1&0\\0&-1}$.
 We take \ref{defI1} as the definition of $I_1$. We thus have, if $L$ is decomposable
$$L = I_1 J \oplus \nabla.$$
Let us set $P=\tr(L^2)$.
\begin{lemma} \label{descPdec} One has $P = \tr(\Delta^2) + \tr (\nabla^2) = I_1^2 + \tr(\nabla^2)$. 
\end{lemma}
\begin{pf} 
 If $m,m' \in L$, we can write 
$m = \delta + \epsilon$, $m' = \delta' + \epsilon'$ with $\delta,\delta' \in \Delta$ and $\epsilon,\epsilon' \in \nabla$.
Then $\tr(m m') = \tr (\delta \delta') + \tr (\epsilon \epsilon')$ since the matrices $\delta \epsilon'$ and $\delta' \epsilon$ are anti-diagonal. Thus, $P \subset \tr(\Delta^2) + \tr (\nabla^2)$, and since the other inclusion is clear, this implies the result.
\end{pf}

\begin{prop} \label{structuredecomposable} Let $L = I_1 J \oplus \nabla \subset (\rad R)^0$ be a decomposable space. The following are equivalent:

\begin{num} There exists a closed subgroup $\Gamma$ of $SR^1$ such that $L$ is the Lie algebra of $\Gamma$. \end{num}
\begin{num} \label{condexGammaL} One has:
\begin{sousnum} \label{nablanabla} $[\nabla,\nabla] \subset I_1 J$,\end{sousnum}
\begin{sousnum} \label{I1Jnabla} $I_1[ J,\nabla] \subset \nabla$, \end{sousnum}
\begin{sousnum} \label{trnabla2I1} $\tr(\nabla^2) I_1 \subset I_1$, \end{sousnum}
\begin{sousnum} \label{trnabla2nabla} $\tr(\nabla^2) \nabla \subset \nabla$, \end{sousnum}
\begin{sousnum} \label{I13} $I_1^3 \subset I_1$, \end{sousnum} 
\end{num}
\end{prop}
\begin{pf} The two first conditions~\ref{nablanabla} and~\ref{I1Jnabla} are equivalent to $L$ being stable by Lie bracket.
Since $P=I_1^2 + \tr(\nabla)^2$, the condition $PL \subset L$ is equivalent to the conjunction of~\ref{trnabla2I1}, \ref{trnabla2nabla}, \ref{I13} and $I_1^2 \nabla \subset \nabla$. But this condition follows from \ref{I1Jnabla}: applied twice, this property gives  $I_1^2 [J, [J,\nabla]] \subset \nabla$, that is $I_1^2 \nabla \subset \nabla$. Therefore the five conditions~\ref{condexGammaL} together are equivalent to $L$ being a Lie subring of $(\rad R)^0$ and $PL \subset L$. The proposition thus follows from Theorems~\ref{thmpinkdirect} and~\ref{pinkconverse}.
\end{pf}

When $L$ is decomposable, we set:
\begin{num} \label{defB1C1} $B_1 := \{b \in B, \exists c \in C, \mat{0 & b \\ c & 0} \in \nabla\},$
\end{num} \begin{num}$ C_1 := \{c \in C, \exists b \in B, 
\mat{0 & b \\ c & 0} \in \nabla\}$.\end{num}
We have obviously $\nabla \subset \mat{0 & B_1 \\ C_1 & 0}$ but the inclusion may be strict.

\subsubsection{Strongly decomposable Lie algebra}

Let $L$ be a closed subspace of $(\rad R)^0$.

\begin{num} \label{Lstrongdec} We shall say that $L$ is {\it strongly decomposable} if, for any$\mat{a & b \\ c & -a} \in L$, one has $\mat{a & 0 \\ 0 & -a} \in L$, $\mat{ 0 & b \\ 0 & 0} \in L $ and $\mat{ 0 & 0 \\ c & 0} \in L$. 
 \end{num}

If we define $B_1, C_1$ and $I_1$ as above \ref{defB1C1}, one can reformulate \ref{Lstrongdec} as
 
 \begin{num} $L =  \mato{I_1 & B_1 \\ C_1 & I_1}$. \end{num}
 If $P=\tr(L^2)$, then we see that
 \begin{num} \label{Pstrongdec} $P=I_1^2 + B_1 C_1$ \end{num}
 
 \begin{prop} \label{conditionstrongdec} Let $L =  \mato{I_1 & B_1 \\ C_1 & I_1} \subset (\rad R)^0$ be a closed subgroup. The following are equivalent:

\begin{num} \label{exGammaL}There exists a closed subgroup $\Gamma$ of $SR^1$ such that $L$ is the Lie algebra of $\Gamma$. \end{num}
\begin{num} \label{condexGammaLstrong} One has:
\begin{sousnum} \label{nablanablastrong} $B_1 C_1 \subset I_1$,\end{sousnum}
\begin{sousnum} \label{I1Jnablastrong} $I_1 B_1 \subset B_1$ and $I_1 C_1 \subset C_1$ , \end{sousnum}
\begin{sousnum} \label{I13strong}$I_1^3 \subset I_1$, \end{sousnum} 
\end{num}
\end{prop}
\begin{pf} 
If $L$ is strongly decomposable, it is in particular decomposable, and we use the notation of \S\ref{decomposable}:
$L=I_1 J \oplus \nabla$ with  $\nabla=\mat{0 & B_1 \\ C_1 & 0}$. One thus has $[\nabla, \nabla] = B_1 C_1 J$ and  $I_1 [J,\nabla] = \mat{0 & I_1 B_1 \\ I_1 C_1 & 0}$, so \ref{nablanablastrong} is equivalent to \ref{nablanabla} and \ref{I1Jnablastrong} is equivalent to \ref{I1Jnabla}.

Since $\tr(\nabla^2) = B_1 C_1$, \ref{trnabla2I1} reads $B_1 C_1 I_1 \subset I_1$, which is a consequence of the above.
Similarly, \ref{trnabla2nabla} read $B_1 C_1 B_1 \subset B_1$ and $B_1 C_1 C_1 \subset C_1$, both of which follow from the above. Thus we see that \ref{condexGammaLstrong} is equivalent to \ref{condexGammaL} and the proposition follows.
\end{pf}

\section{Admissible pseudo-representations}

\subsection{Hypotheses on the base ring $A$}

In all this section, we let $\F$ be a finite field of characteristic $p$, and we denote by $W(\F)$ the ring of Witt vectors of $\F$.
We suppose given 
\begin{num} \label{baseringA2} A topological $W(\F)$-algebra $A$ which is compact and semi-local, and such that the maps $W(\F) \rightarrow A/\m_i$, where $\m_i$, $i=1,\dots,r$ are the maximal ideals of $A$, are surjective.
\end{num}
Thus $A$ satisfies the condition \ref{baseringA} with the small additional requirements that $A$ is a topological $W(\F)$-algebra and that
the maps $W(\F) \rightarrow A/\m_i$ are surjective, which implies that the residue fields $\F_i$, $i=1,\dots,r$, are all equal at $\F$. We use the same notations as in the preceding section: $A=\prod_{i=1}^r A_i$ with the $A_i$'s local,  and we write (by abuse) $\m_i$ for the maximal ideal of $A_i$.

We shall denote by $s : \F \rightarrow A$ the map obtained by taking the Teichmuller lift in $W(\F)$ of an element of $\F$ and seeing it as an element  of $A$ through the structural map $W(\F) \rightarrow A$.
The map $s$ is a set-theoretical section of the residue map $A \rightarrow A/\m = \F$, and preserve multiplication but not addition.
The elements of $A$ that belong to $s(\F)$ will be called {\it constants}.

\subsection{Admissible pseudo-deformations}

\label{admissiblepd}
We now proceed to define an admissible pseudo-deformation over $A$. It is a $4$-tuple $(\Pi,\rhob,t,d)$ where
\begin{num} \label{Pipfinite} $\Pi$ is a profinite group which satisfies Mazur's {\it finiteness $p$-condition} (i.e. the maximal pro-$p$-quotient of every open subgroup of $\Pi$ is topologically finitely generated.)\footnote{Actually, this hypothesis is used only when at least one of the $\rhob_i$ is reducible, to ensure that $R$ is provided with a natural topology below, or in Corollary~\ref{admissiblenoetherian} which itself is only needed in section \S\ref{sectioncongruencelarge}. In other words, all the results concerning the dihedral, large, and exceptional cases in the next section does not need this hypothesis.} 
\end{num}
\begin{num} \label{tmultfree} $\rhob=(\rhob_i)_{i=1}^r$ is a family of isomorphism classes of continuous representations $\rhob_i : \Pi \rightarrow \GL_2(\F)$, each of them being either   absolutely irreducible or the sum of two distinct characters.\end{num} \noindent
\begin{num} \label{deform} $(t,d)$ is a continuous pseudo-representation of $\Pi$ over $A$ such that for $i=1,\dots,r$
we have $\tr \rhob \equiv t_i \pmod{\m_i}$ and $\det \rhob_i \equiv d \pmod{\m_i}$. \end{num}
\begin{num} \label{dconstant} We have $d(g) \in s(\F)$ for all $g \in \Pi$ \end{num}
\begin{num} \label{tgenerate}  As a topological $W(\F)$-algebra, $A$ is generated by $t(\Pi)$. \end{num}

The condition \ref{dconstant} expresses the fact that this pseudo-representation has {\it constant determinant}. 
Even if we do not assume it, there is always a twist of $(t,d)$ which has constant determinant, namely the twist by the character $g \mapsto \sqrt{d(g)^{-1} s(d(g))}$.

If we denote by $(t_i,d_i)$ the composition of $(t,d)$ with $A \rightarrow A_i$, the condition \ref{deform} says that $(t_i,d_i)$ is a {\it deformation} over $A$ of the pseudo-representation $(\tr \rhob_i, \det \rhob_i)$ attached to $\rhob_i$, or as it is customary to say, a {\it pseudo-deformation} of $\rhob_i$ over $A$. 

If $A \rightarrow A'$ is a surjective map, then $A'$ with its quotient topology satisfies~\ref{baseringA2}, and we can write $A'=\prod_{j \in J} A'_j$, where $J$ is a subset of $\{1,\dots,r\}$ and the map $A \rightarrow A'$ is the product of surjective maps $A_j \rightarrow A'_j$
for $j \in J$. If we denote by $(t',d')$ the composition of $(t,d)$ with the map $A \rightarrow A'$, then it is clear that $(\Pi,(\rhob_j)_{j \in J},t',d')$ is an admissible pseudo-deformation over the ring $A'$. In particular, for every $i=1,\dots,r$, $(\Pi,\rhob_i,t_i,d_i)$ is an admissible pseudo-deformation over the local ring $A_i$.

\subsection{Equivalent formulations for~\ref{tgenerate}}

Following \cite{LenstraDeSmit}, let $\CC$ be the category of topological $W(\F)$-algebra $B$ that are compact and local, and such that the map $W(\F) \rightarrow B/\m_B$ is surjective, where $\m_B$ is the maximal ideal of $B$. Given a topological group $\Pi$ and a continuous representation $\rhob: \Pi \rightarrow \GL_2(\F)$,
we consider the functor $\Fc_\rhob$ from $\CC$ to the categories of sets, such that $\Fc_\rhob(B)$ is the set of continuous pseudo-representation $(t,d): \Pi \rightarrow B$ such that $t \equiv \tr \rhob \pmod{\m_B}$,  $d \equiv \det \rhob \pmod{\m_B}$, and $d(g) \in s(\F)$ for all $g \in \Pi$. By \cite{LenstraDeSmit}, this functor is representable by a ring $A_{\rhob,\univ}$.
 
 Let $(\Pi,\rhob,t,d)$ be a pseudo-representation over $A=\prod_i A_i$ satisfying \ref{tmultfree}, \ref{deform} and \ref{dconstant}, and let $i \in \{1,\dots,r\}$. Thus $(\Pi,t_i,d_i)$ defines an element of $\Fc_{\rhob_i}(A_i)$ hence a map $A_{\rhob_i,\univ} \rightarrow A_i$.
 
\begin{prop} \label{tgenerateuniv} $(\Pi,\rhob,t,d)$ satisfies~\ref{tgenerate} (i.e. is admissible) if and only if for $i=1,\dots,r$, the morphisms $A_{\rhob_i,\univ} \rightarrow A$ are surjective.
\end{prop}
This is clear. 
\begin{cor} \label{admissiblenoetherian} If $(\Pi,\rhob,t,d)$ is an admissible pseudo-deformation over $A$, then $A$ is noetherian.
\end{cor}
\begin{pf} Since $\Pi$ satisfies Mazur's $p$-finiteness condition, we now that $A_{\rhob_i,\univ}$ is noetherian by \cite{LenstraDeSmit} in the case $\rhob_i$ absolutely irreducible, by \cite{pseudodef} in the case $\rhob_i$ reducible and $p>2$ and by \cite{chenevier} in the case $p=2$. Thus $A_i$ is noetherian for all $i$, and $A$ is noetherian.
\end{pf}
 
\begin{prop} Assume $p>2$. In the definition of an admissible pseudo-representation, condition~\ref{tgenerate} can be replaced by the apparently weaker condition
\begin{num}  \label{tgeneratemodule}  As a topological $W(\F)$-module, $A$ is generated by $t(\Pi)$. \end{num}
\end{prop}
Indeed, the $W(\F)$-module generated by $t(\Pi)$ is already a $W(\F)$-algebra, for it contains $t(1)=2$, hence $1$ since $p>2$, and it is stable by multiplication: if $x,y \in \Pi$, $t(x) t(y) = t(xy) + d(y) t(xy^{-1})$, and $d(y) \in W(\F)$ by \ref{dconstant}.

\subsection{$(t,d)$-representations attached to an admissible pseudo-deformation and their image}

 If $(\Pi,\rhob,t,d)$ is an admissible pseudo-deformation, then for every $i \in \{1,\dots,r\}$, there exists, by Theorem~\ref{exGMA}, a unique up to unique isomorphism $A_i$-GMA $R_i$ and a $(t_i,d_i)$-representation $\rho_i: \Pi \rightarrow R_i^\ast$.  Let us remind that that means that there exist a faithful GMA $R_i=\mat{ A_i & B_i \\ C_i & A_i}$ and a representation $\rho_i: \Pi \rightarrow \GL_2(A_i)$ of trace $t_i$ and determinant $d_i$, and that given another GMA $R'_i$ and representation $\rho'_i$ satisfying the same conditions, there exists a unique isomorphism of $A$-algebras $f: R_i \rightarrow R'_i$
such that $f \circ \rho_i = \rho'_i$. We note that by Corollary~\ref{admissiblenoetherian} and Theorem~\ref{exGMA}, the ring $A_i$ is noetherian, the algebra $R_i$ is finite-type as an $A_i$-module, and when $R_i$ is provided with its natural topology, the representation $\rho_i$ is continuous.

Setting $R=\prod_{i=1}^r R_i$ and seeing this ring as an $A=\prod_{i=1}^r A_i$-algebra (component-wise),  we get a continuous representation $\rho : \Pi \rightarrow R^\ast$ of trace $t$ and determinant $d$ which is unique up to unique isomorphism. We call this representation a $(t,d)$-representation.

Given such a representation $\rho$, we set
\begin{num} \label{defG} $G = \rho(\Pi)$ \end{num}
\begin{num} \label{defGamma} $\Gamma = G \cap SR^1$, \end{num}
\noindent  where $SR^1$ is defined as in~\S\ref{subpinksetting}.
Note that $G$ is a closed subgroup of $R^\ast$ and $\Gamma$ a closed subgroup of $SR^1$.

We denote by  $\Gb$ the image of $G$ by the map $R^\ast \rightarrow (R/\rad R)^\ast$
\begin{lemma} The sequence
\begin{num} \label{exactGammaG} $1 \rightarrow \Gamma \rightarrow G \rightarrow \Gb \rightarrow 1$\end{num}
\noindent is exact. In particular, $\Gamma$ is a finite index normal subgroup in $G$.
\end{lemma}
\begin{pf}  Though $\Gamma$ is defined as $G \cap SR^1$, we claim that $\Gamma=G\cap R^1$.  Indeed, let $g \in G \cap R^1$ and write $g = \rho(x)$ for $x \in \Pi$. Then $\det g =s( \overline{\det(g)})$ by~\ref{dconstant}. Since $g \in R^1$, $\det g \in 1+\m \subset A^\ast$ and $\overline{\det(g)}=1$. Thus $\det(g)=s(1)=1$ and $g \in \Gamma$.

 Since the kernel of $G \rightarrow \Gb$ is $G \cap R^1$, the result follows.
\end{pf}

We also define
 \begin{num} \label{defGi} $G_i = \rho_i(\Pi)$ \end{num}
 The group $G_i$ is the image of $G$ by the map $R^\ast \rightarrow R_i^\ast$. The surjective maps $G \rightarrow G_i$ for $i=1,\dots,r$ define a map $G \rightarrow \prod_{i=1}^r G_i$ which is always injective, but not necessarily surjective.

We observe that the choice of a representation $\rho_i$ specifies a single representation $\tilde \rho_i: \Pi \rightarrow \GL_2(\F)$ in the isomorphism class $\rhob_i$, as follows: consider the composition  $\tilde \rho_i: \Pi \stackrel{\rho}\rightarrow R_i^\ast \rightarrow (R_i/\rad R_i)^\ast$. We know 
that $R_i /\rad R_i$ is $M_2(\F)$ if $\rhob_i$ is absolutely irreducible and $\mat{\F & 0 \\ 0 & \F}$ otherwise, so $\tilde \rho_i$ can be considered in both cases as a semi-simple representation of $G$. The trace and determinant of $\tilde \rho_i$ are reduction mod $\m_i$ of those of $\rho_i$, hence are identical to those of $\rhob_i$. Therefore, $\tilde \rho_i$ is a representation in the equivalence class $\rhob_i$. By a slight abuse of notations, when a representation $\rho_i$ is fixed, we shall denote by $\rhob_i$ its reduction $\tilde \rho_i$.

\par \bigskip 

\section{Lie-theoretic study of admissible pseudo-deformations}

\label{sectionLieStudy}

\subsection{Hypothesis on the base ring $A$.}

In this section, we let $\F$ be a finite field of characteristic $p>2$, and we consider
\begin{num} \label{baseringA3} A topological ring $A$ which is compact and local, with residue field $\F$. \end{num}
Such a ring $A$ is automatically a topological $W(\F)$-algebra, and the map $W(\F) \rightarrow A \rightarrow A/\m=\F$ is the residue map of $W(\F)$, hence surjective. Hence our hypothesis implies \ref{baseringA2}, and actually is equivalent to it combined with the supplementary assertion that $A$ is local (and $p>2$).

Our aim is to study the image $G$ of $\rho$, with a special attention to its subgroup $\Gamma$. The group $G$ depends on the chosen $(t,d)$-representation $\rho: \Pi \rightarrow R$, but only up to unique isomorphism. We can choose to work with any 
$(t,d)$-representation $\rho: \Pi \rightarrow R^\ast$ that simplifies our analysis. 
According to~\ref{tmultfree}, there is an element $g_0 \in \Pi$ such that $\rhob(g_0)$ has two distinct eigenvalues in $\F$, $\lambda_0$ and $\mu_0$. Actually, there are in general many of them. Given such an element $g_0$ as well as an ordering $(\lambda_0,\mu_0)$ of 
the eigenvalues of $\rhob(g_0)$, there exists  a $(t,d)$-representation 
$\rho : \Pi \rightarrow R^\ast$ adapted to $(g_0, \lambda_0, \mu_0)$. Let us remind that that means that $\rho(g_0)$ is a diagonal matrix which reduces modulo $\m$ to $\mat{\lambda_0 & 0 \\ 0 & \mu_0}$. We shall see that working with $(t,d)$-representations $\rho$ which are adapted to a well-chosen element $g_0$ is often the right choice.

In order to study the group $G$, and its subgroup $\Gamma$, we shall make use of the generalization of Pink's theory exposed in the preceding section. Note that the GMA $R$ is Cayley-Hamilton, since it is faithful, and that $\Gamma$ is a closed subgroup of $SR^1$, so this theory 
applies and attach to $\Gamma$ a Lie subring $L=L(\Gamma)$ of $(\rad R)^0$. To $L$ is attached a pseuso-ring $P = \tr(L^2)$
such that $PL \subset L$, and the full descending spectral sequence $L_1=L$, $L_2=[L,L]$, etc. 

\subsection{Finding constant elements in $G$}

Given a faithful GMA $R$ over $A$, the multiplicative section $s: \F \rightarrow A$ induces a set-theoretic section of the map $R \rightarrow R/\rad R$. This section, still denoted by $s: R/\rad R \rightarrow R$, sends a matrix $\smat{a & b \\ c & d}$ to $\smat{s(a) & s(b) \\ s(c) & s(d)}$ in the case $R = M_2(A)$ and $\smat{a & 0 \\ 0 & d}$ to $\smat{s(a) &0\\ 0 & s(d)}$ in the case $R =\mat{A & B \\ C & D}$ with $BC \subset\m$. We shall call a matrix of $R$ which lies in $s(R/\rad R)$ {\it constant}. 

Note that the section $s$ is multiplicative in the second case, but is not in the first, because multiplications of matrices involve addition of the coefficients in general, and $s$ does not preserve addition. However, when $m, m'$ are two matrices in $R/\rad R$ which are either diagonal or anti-diagonal, then $s(mm')=s(m) s(m')$ because in this case the multiplication of matrices only involve multiplication of the coefficients.

We consider again in this subsection an admissible pseudo-deformation $(\Pi,\rhob,t,d)$ over $A$. Given a $(t,d)$-representation $\rho: \Pi \rightarrow R^\ast$, we recall that by definition $G=\rho(\Pi)$, and $\Gamma=G \cap SR^1$.

Our aim is to find elements of the image $G$ that are constant.
It is important to observe that {\bf the notion of constant element of $G$  depends on the chosen $(t,d)$-representation $\rho$.}
 Therefore, our aim is, more precisely stated, for a given admissible pseudo-deformation $(\Pi,\rhob,t,d)$ 
to find a suitable $(t,d)$-representation $\rho: \Pi \rightarrow R^\ast$ such that the associated group $G$ has enough constant elements.

\begin{theorem} \label{diagonallift} Let $g_0$ be such that $\rhob(g_0)$ has distinct eigenvalues $\lambda_0,\mu_0$ in $\F$, and
let $\rho: \Pi \rightarrow R^\ast$ be any $(t,d)$-representation adapted to $(g_0,\lambda_0,\mu_0)$.
Let $D$ be the subgroup of $\Gb$ generated by $\rhob(g_0)$ and by the scalar matrices in $\Gb$. Then $s(D) \subset G$. 

Furthermore, let $n \in N(D)-Z(D)$, where $N(D)$ is the normalizer and $Z(D)$ is the centralizer of $D$ in $\bar G$. 
Then, up to changing $\rho$ into another $(t,d)$-representation adapted to $(g_0,\lambda_0,\mu_0)$, one has $s(n) \in G$.
As a consequence, if $D=Z(D)$ then $s(N(D)) \subset G$. 
\end{theorem}
\begin{pf} By assumption $\rho(g_0)$ is diagonal and reduces modulo $\rad R$ to $\rhob(g_0)=\mat{\lambda_0 & 0 \\ 0 & \mu_0}$. Let us write $\rho(g_0) = s(\rhob(g_0))+m$ with $m \in \rad R$ a diagonal matrix.

 Since $s(\rhob(g_0))$ and $m$ commute, being two diagonal matrices, we get for every integer $n \geq 1$ (denoting by $q$ the cardinality of $\F$):
$$\rho(g_0^{q^n}) = s(\rhob(g_0)^{q^n}) + \sum_{k =1}^n {q^n \choose k}  s(\rhob(g_0)^{q^n-k}) m^{k}.$$
Denoting by $v_p$ the $p$-valuation of an integer, one has $v_p \left( {q^n \choose k} \right) = n v_p(q) - v_p(k)$ if $k \geq 1$, as is 
well-known.
The matrix $m^k$ is diagonal with coefficients in $\m^k$, and $s(\rhob(g_0)^{q^n-k})$ is diagonal with coefficients in $A$.
Therefore, since $p \in \m$, the term ${q^n \choose k}  s(\rhob(g_0)^{q^n-k}) m^{k}$ for $k \geq 1$ is a diagonal matrix whose coefficients belong to $\m^{n v_p(q)-v_p(k)+k}$, hence to $\m^{n v_p(q) + 1}$.
 
On the other hand, since $\rhob(g_0)$ is a diagonal matrix in $\GL_2(\F)$, its order divides $q-1$, hence $\rhob(g_0)^q=\rhob(g_0)$ and $\rhob(g_0)^{q^n}=\rhob(g_0)$. 

Therefore  $$\rho(g_0^{q^n})  \equiv s(\rhob(g_0))  \pmod{\m^{n v_p(q)+1}}$$
Since $\rho(g_0^{q^n})$ belongs to $G$ by definition, and $n v_p(q)+1$ tends to $+\infty$, we see that $s(\rhob(g_0))$ is the limit of a sequence of elements of $G$. Since $G$ is closed,  $$s(\rhob(g_0)) \in G.$$

Let $h \in \Pi$ such that $\rhob(h)$ is  a scalar matrix. Then we can write $\rho(h)=s(\rhob(h)) + m$ with $m \in \rad R$ a matrix commuting with $s(\rhob(h))$ (since $s(\rhob(h))$ is a scalar matrix in $R$). The same argument as above then shows that  $s(\rhob(h)) \in G.$ Since $D$ is generated by $\rhob(g_0)$ and the scalar matrices in $D$, and $s_{|D}$ is a morphism of groups, we have $s(D)\subset G$. This proves the first assertion of the theorem.

Now let $N$ be the normalizer of $D$ in $\bar G$, and $Z$ its centralizer. If $N=Z$ there is nothing else to prove. If $N \neq Z$, then there is an anti-diagonal element in $N$, which shows that $\rhob$ is irreducible and we are in the case $R=M_2(A)$.
It is easy to see that $|N|=2|Z|$. Since $Z$ consists of diagonal matrices, $|Z|$ divides $(q-1)^2$, and
the order $|N|$ is prime to $p$. Considering the exact sequence $1 \rightarrow \Gamma \rightarrow G \rightarrow \Gb \rightarrow 1$, and the fact that $\Gamma$ is a pro-$p$-group, we see by Zassenhaus' theorem that there is a map $s': N \rightarrow G$ which is a section of $G \rightarrow \Gb$ over $N \subset \Gb$. The restriction of $s'$ to $D$ is a section over $D$ of $G \rightarrow \Gb$. Since $|D|$ is prime to $p$, such a section is unique up to conjugation (again by Zassenhaus' theorem) by an element $g$ of $G$. Replacing $s'$ by $gs'g^{-1}$ we may assume that the section $s'$ on $N$ restricts to the section $s$ on $D$. 

Let us choose $n \in N-Z$. The element $n$ normalizes $D$ and therefore $s'(n)$ normalizes $s'(D)=s(D)$, which is a non-scalar diagonal subgroup of $R^\ast$. Therefore $s'(n)$ is either diagonal or anti-diagonal. If it was diagonal, then it would commute with $s(D)$, hence $n$ would commute with $D$ and be in $Z$, a contradiction. Therefore $s'(n)$ is anti-diagonal, say 
$s'(n)=\smat{0 & b \\ c & 0}$. Since $n^2 \in D$, $s'(n^2) = \smat{ bc & 0 \\ 0 & bc}$ is in $s(D)$ and therefore $bc \in s(\F)$.
By conjugating $\rho$ by the matrix $\smat{s(b) & 0 \\ 0 & 1}$, we may assume that $b=1$ (with $\rho$ still a $(t,d)$-representation adapted to $g_0$.) Thus $c \in s(\F)$, and therefore $s'(n) = s(n)$. It follows that $s(n) \in G$, as claimed.
\end{pf}

Let us note two important consequences:
\begin{cor} \label{Ldecomp} Let $\rho$ be adapted to an element $(g_0,\lambda_0,\mu_0)$ as above, and let $G$, $\Gamma$, $L$ defined using this $\rho$. Then $L$ is decomposable.
\end{cor}
\begin{pf} Let us denote by $u : R \rightarrow R$ the {\it conjugation by $s(\rhob(g_0))$}, that is the map $m \rightarrow s(\rhob(g_0)) m s(\rhob(g_0))^{-1}$. The map $u$ is a $W(\F)$-linear endomorphism of $R$. 
By the theorem $s(\rhob(g_0)) = \smat{s(\lambda_0) & 0 \\ 0 & s(\mu_0)}$ is in $G$, and therefore normalizes $\Gamma$, hence $L$.
In other words, $u$ stabilizes the additive subgroup $L$ of $R$. 

In order to simplify notation, let us set $r := s(\lambda_0/\mu_0) \in W(\F)$.  Clearly, $u$ fixes diagonal matrices in $R$, and acts by multiplication by $r$ (resp. $r^{-1}$) on matrices of the form $\smat{0 & b \\ 0 & 0}$ (resp. $\smat{0 & 0 \\ c & 0}$). It follows that $u$ is killed by the polynomial $X (X-r) (X-r^{-1})$. If $\Sigma = \Gal(\F/\F_p) = \Aut_{\Z_p} W(\F)$, then the polynomial
$X Q(X)$ also kills $u$, with $Q(X) = \prod_{\sigma \in \Sigma} (X-\sigma(r)) (X-\sigma(r)^{-1}) \in \Z_p[X]$.
Since by assumption, $r \neq 1$, the value $Q(1)$ is invertible in $\Z_p$ and the operator $Q(u)/Q(1)$ of $R$ is the projection
onto diagonal matrices relatively to antidiagonal matrices. This operator, being in $\Z_p[u]$, stabilizes $L$, which shows that if a matrix is in $L$, its diagonal part is also in $L$.
\end{pf}

\begin{cor} \label{Lstrongdecomp} Let $\rho$ be adapted to an element $(g_0,\lambda_0,\mu_0)$ as above, and let $G$, $\Gamma$, $L$ defined using this $\rho$. Let $\F_q$ be a subfield of $\F$, and assume that there exists an integer $n$ such that $\lambda_0^n/\mu_0^n \in \F_q^\ast - \{1,-1\}$. Then $W(\F_q) L$ is strongly decomposable. More precisely, $L$ is decomposable, and with $I_1$, $B_1$, $C_1$
as in \S\ref{decomposable}, one has $W(\F_q) L = \mato{W(\F_q) I_1 & W(\F_q) B_1 \\ W(\F_q) C_1 & W(\F_q) I_1}$
\end{cor}
\begin{pf} We already know that $L$, hence $W(\F_q)L$, is decomposable. Using the notation of the previous proof, the hypothesis becomes $r^n \in s(\F_q-\{1,-1\})$ and it follows that $r^n-r^{-n}$ is invertible in $W(\F_q)$. The operator 
$(u^n - r^n)/(r^{-n}-r^n)$ acts on anti-diagonal matrices of $R$ as the map $\smat{0 & b \\ c & 0} \mapsto \smat{0 & b \\ 0 & 0}$,
and this operator stabilizes $W(\F_q) \nabla$. The result follows.
\end{pf}

\subsection{Consequences of Theorem~\ref{diagonallift} in the cases of cyclic or dihedral projective image of $\rhob$}

\subsubsection{Well-adapted $(t,d)$-representations and splitting of the exact sequence~\ref{exactGammaG}}

We still consider an admissible pseudo-deformation $(\Pi,\rhob,t,d)$. In the cases $\rhob$ of abelian or dihedral projective image, we shall use the following terminology:
\begin{definition} \label{welladapted} A $(t,d)$-representation $\rho$ is said to be {\it well adapted} if 
\begin{itemize} \item[(i)] The representation $\rho$ is adapted to an element $g_0 \in \Pi$ such that $\rhob(g_0)$ together with the scalar matrices
in $\Gb$ generates $\Gb$ in the cyclic case, and a subgroup of index $2$ in $\Gb$ in the dihedral case.
\item[(ii)] $s(\Gb) \subset G$.
\item[(iii)] If $\Gb$ is non-abelian, then it contains a matrix of the form $\smat{ 0 & b \\ c & 0}$ with $bc^{-1} \in \F_p^\ast$. 
\end{itemize} 
\end{definition}
Note that in the abelian case, (ii) follows from (i) by Theorem~\ref{diagonallift} and (iii) is empty.

\begin{prop} \label{propsplit} Assume that the projective image of $\rhob$ is either cyclic or dihedral. 
Then there  exists a $(t,d)$-representation $\rho$ that is well adapted. Moreover, for such a $\rho$ the restriction of $s$ to $\Gb$ is a group-theoretic section of that exact sequence, and $G$ is therefore the semi-direct product of $\Gamma$ by $\Gb$, acting on $\Gamma$ by $g \cdot \gamma = s(g) \gamma s(g)^{-1}$.
\end{prop}
\begin{pf} Let $D$ be the group $\Gb$ if $\rhob$ is reducible, and $D$ be a subgroup of index $2$ in $\Gb$ containing all scalar matrices if $\rhob$ is dihedral.
In both case, one has $D=Z(D)$ and $D$ is diagonal in a certain basis, which implies that $D$ modulo its subgroup of scalar matrices is cyclic, say generated by $\rhob(g_0)$. By \ref{tmultfree}, $\rhob(g_0)$ is not scalar, and thus has two distinct eigenvalues $(\lambda_0,\mu_0)$. Let us choose for $\rho$ a $(t,d)$-representation adapted to $(g_0,\lambda_0,\mu_0)$ and, in the case $\rhob$ dihedral, chosen as to satisfy the second paragraph of Prop.~\ref{diagonallift}. Then by Prop.~\ref{diagonallift}, one has $s(N(D)) \subset G$ and since $N(D) = \Gb$, we see that $s$ is a section of $1 \rightarrow \Gamma \rightarrow G \rightarrow \Gb \rightarrow 1$. Moreover $\rho$ satisfies (i) and (ii) of the definition of a well adapted representation. Since $\Gb$ normalizes $D$ but is not abelian, it must contain a matrix of the form $\smat{0 & b \\ c & 0}$.
Up to conjugating $\Gb$ by $\smat{1 & 0 \\ 0 & x}$, it contains the matrix $\smat{0 & b x \\ c x^{-1} & 0}$. One can choose $x \in \F^\ast$ such that $(bx) (cx^{-1})^{-1} = bc^{-1} x^2$ be in $\F_p^\ast$. Thus, conjugating $\rho$ by $s\left (\smat{1 & 0 \\ 0 & x}\right)$ doesn't affect property (i) and (ii) and ensure property (iii).
\end{pf}

\begin{cor} \label{corsplit} Assume that the projective image of  $\rhob$ is either cyclic or dihedral. Then the exact sequence $1 \rightarrow \Gamma \rightarrow G \rightarrow \Gb \rightarrow 1$ is split.
\end{cor}

Note that a for well adapted $\rho$, the corresponding Lie Algebra $L$ is decomposable (Cor.~\ref{Ldecomp}) and can be written $L= I_1 J \oplus \nabla$.

\subsubsection{Consequences in the cyclic case}

\begin{cor}  \label{Amreducible} Assume that the projective image of $\rhob$ is cyclic. Let $\rho$ be a well adapted $(t,d)$-representation. Then one has (with the notation of \S\ref{decomposable})
\begin{num} $W(\F) 1 + W(\F) I_1+ W(\F) I_1^2 + W(\F) \tr(\nabla^2) = A$. \end{num}
\begin{num} The $A$-module generated by $B_1$ is $B$. \end{num}
\begin{num} The $A$-module generated by $C_1$ is $C$. \end{num}
\end{cor}
\begin{pf} By \ref{tgeneratemodule}, $W(\F) \tr(G)=A$. By Prop.~\ref{propsplit}, every element $g$ in $G$ can be written $g=\gamma \smat{\lambda_1 & 0 \\ 0 & \lambda_2}$ with $\lambda_1,\lambda_2 \in s(\F^\ast) \subset W(\F)$ and $\gamma \in \Gamma$. We can write $\gamma=\theta^{-1}\left(\smat{a & b \\ c & -a}\right)$, with $\smat{a & b \\ c & -a} \in L$.
 We have $\tr g = (\lambda_1-\lambda_2) a  + (\lambda_1 + \lambda_2) \sqrt{1+{a^2+bc}}$; the first term on the RHS is in $W(\F) I_1$, and the second in $W(\F) 1 + W(\F) P=W(\F) 1 + W(\F) I_1^2 + W(\F) \tr(\nabla^2)$ by Lemma~\ref{descPdec}. The first result follows.

For the second and third, if $g \in G$ is written $g=\gamma \smat{\lambda_1 & 0 \\ 0 & \lambda_2}$ as above,
then the anti-diagonal part of $g$ is $\smat{0 & \lambda_1 b \\ \lambda_2 c & 0}$ which belongs to $\smat{0 & W(\F) B_1 \\ W(\F) C_1 & 0}$.  Recalling that $G$ generates $R$ as an $A$-module, we get $A B_1=B$ and $A C_1 =C$.
 \end{pf}
 
\subsubsection{Consequences in the dihedral case}

\label{consequencesdihedral}

We now make some general observations concerning the case where the projective image of $\rhob$ is dihedral.
In this case, choosing a well adapted $(t,d)$-representation $\rho$ defines an abelian subgroup of index $2$ in $\Gb$, namely the subgroup $D$ generated by $\rhob(g_0)$ and the scalar matrices in $\Gb$. When the projective image of $\rhob$ has order $>4$, then this group $D$ is the unique abelian subgroup of index $2$ in $\Gb$, hence is independent of the choice of $\rho$, but when $\Gb = \Z/2\Z \times \Z/2\Z$, there are three possible index $2$ subgroups $D$ in bar $G$, and each of them is associated with a well-adapted $(t,d)$-representation $\rho$.

In any case, we fix a well-adapted $\rho$, which fixes a cyclic subgroup $D$ of index 2 in $\Gb$, and we define $\Pi'$ as the inverse image of $D$ by the map $G \rightarrow \Gb$. Hence $\Pi'$ is a subgroup of index $2$ of $\Pi$.
The image $G'=\rho(\Pi')$ lies in an exact sequence $1 \rightarrow \Gamma \rightarrow G' \rightarrow D \rightarrow 1$, and is exact sequence is split, a splitting being the restriction of $s$ to $D$. 

By Lemma~\ref{imagesubGma}, the sub-$A$-module $AG'$ of $R$ is a sub-$A$-GMA of $R=M_2(A)$, that is of the form $\mat{A & B \\ C & A}$ with $B,C$ ideals of $A$. Since $G$ contains anti-diagonal matrices with coefficients in $s(\F) \subset A^\ast$, and normalizes $AG'$,
one has $B=C$, and $R=\mat{A & B \\ B &A}$. It is not hard to see that the ideal $B$ depends only of the admissible pseudo-representation $(\Pi,\rhob,t,d)$ and the subgroup $D$ of $G'$, not of the choice of the well-adapted $(t,d)$-representation $\rho$: see e.g. Prop.~\ref{IBKab} below. 

\par\medskip

We write as usual $L=I_1 J \oplus \nabla$, and $B_1,C_1$ for the subgroups of upper-right and lower-left coefficients of $\nabla$; since
elements in $\Gamma$ have upper-right and lower-left coefficients in $B$, and $\Theta$ does not affect non-diagonal coefficients, we have $B_1 \subset B$, $C_1 \subset B$.

\begin{cor}  \label{Amdihedral} If $\rhob$ is dihedral, and $\rho$ is a well adapted $(t,d)$-representation, then:
\begin{num} \label{symnabla}There exists $\lambda \in s(\F_p^\ast)$ such that the subgroup $\nabla$ of $\mat{0 & B \\ B & 0}$ is stable by the map $\smat{0 & b \\ c & 0} \mapsto \smat{0 & \lambda c \\ \lambda^{-1} b & 0}$. In particular $B_1=C_1$. \end{num}
\begin{num} \label{WF1dihedral} One has $W(\F) 1 + W(\F) I_1+ W(\F) I_1^2 + W(\F) \tr(\nabla^2) + W(\F) B_1= A$. \end{num}
\begin{num} The $A$-module generated by $B_1$ is $B$. \end{num}
\end{cor}
\begin{pf}
By definition of a well adapted representation, the group $G$ contains a matrix $\smat{0 & s(\beta) \\ s(\gamma) & 0}$ with $s(\beta \gamma^{-1}) \in \F_p^\ast$.
The conjugation by that matrix stabilizes $\Gamma$, $L$, and $\nabla$, and is given by $\smat{0 & b \\ c & 0} \mapsto \smat{0 & \lambda c \\ \lambda^{-1} b & 0}$ with $\lambda = s(\beta \gamma^{-1})$. The first part of~\ref{symnabla} follows and we have $C_1 = \lambda^2 B_1$. Since $B_1$ is a $\Z_p$-module, and $\lambda \in \Z_p^\ast$, one gets  $C_1  = B_1$.

By \ref{tgenerate}, $W(\F) \tr(G)=A$. Every element $g$ in $G$ can be written either $g=\gamma \smat{\lambda_1 & 0 \\ 0 & \lambda_2}$ or $g=\gamma \smat{ 0 & \lambda_1 \\ \lambda_2 & 0}$ with $\lambda_1,\lambda_2 \in s(\F^\ast) \subset W(\F)$. We can write $\gamma=\theta^{-1}\smat{a & b \\ c & -a}$, with $\smat{a & b \\ c & -a} \in L$. 
In the first case, we have $\tr g = (\lambda_1-\lambda_2) a  + (\lambda_1 + \lambda_2) \sqrt{1+{a^2+bc}}$; the first term on the RHS is in $W(\F) I_1$, and the second in $W(\F) 1 + W(\F) P=W(\F) 1 + W(\F) I_1^2 + W(\F) \tr(\nabla^2)$.  In the second case, we have $\tr(g) = \lambda_1 b + \lambda_2 c \in W(\F) B_1$. Formula~\ref{WF1dihedral} follows.

Finally,  any $g \in G'$ can be written $g=\gamma \smat{\lambda_1 & 0 \\ 0 & \lambda_2}$ as above, and 
 the anti-diagonal part of $g$ is $\smat{0 & \lambda_1 b \\ \lambda_2 c & 0}$ which belongs to $\mat{0 & A B_1 \\ A B_1 & 0}$.  Recalling that by definition $G'$ generates $\mat{A & B \\ B & A}$ as an $A$-module, we get $A B_1 =B$.
\end{pf}

\subsection{The structure of $L$ when the projective image of $\rhob$ has order 2}

That is, we consider the case where $\rhob=\chi_1 \oplus \chi_2$ is reducible, with $\chi_1^2 = \chi_2^2$ (but still $\chi_1 \neq \chi_2$ by ~\ref{tmultfree}).
In this case, there is nothing more to say that what we have already said:

\begin{theorem} \label{structureLZ2Z} Let $(\Pi,\rhob,d,t)$ be an admissible pseudo-deformation such that the projective image of $\rhob$ has order 2 and let $\rho: \Pi \rightarrow R^\ast$, $R = \mat{A & B \\ C & A}$ a well adapted $(t,d)$-representation. Then there exists a closed subgroup $I_1$ of $\m$, and a closed subgroup $\nabla$ of $\mat{0 & B \\ C & 0}$ such that $$L = I_1 J \oplus \nabla$$ and
\begin{num}  \label{structureLZ2Z1} $[\nabla,\nabla] \subset I_1 J$,\end{num}
\begin{num} \label{structureLZ2Z2} $I_1[ J,\nabla] \subset \nabla$, \end{num}
\begin{num} \label{structureLZ2Z3}$\tr(\nabla^2) I_1 \subset I_1$, \end{num}
\begin{num} \label{structureLZ2Z4} $\tr(\nabla^2) \nabla \subset \nabla$, \end{num}
\begin{num} \label{structureLZ2Z5}$I_1^3 \subset I_1$, \end{num} 
\begin{num}\label{structureLZ2Z6} $W(\F) 1 + W(\F) I_1+ W(\F) I_1^2 + W(\F) \tr(\nabla^2) = A$ \end{num} 
\begin{num} \label{structureLZ2Z7} $AB_1=B$. \end{num}
\begin{num} \label{structureLZ2Z8}$AC_1=C$. \end{num}
\noindent

Conversely, if $R = \mat{A & B \\ C & D}$ is a faithful GMA over $A$, $I_1$ is any closed subgroup of $\m$, and $\nabla$ any closed subgroup  of $\mat{0 & B \\ C & 0}$ satisfying the eight conditions above, then there exists an admissible pseudo-deformation
$(\Pi,\rhob,t,d)$ with $\rhob$ of projective image of order $2$, and a $(t,d)$-representation $\rho : \Pi \rightarrow R^\ast$ such that the Lie algebra attached to $\rho$ is $L= I_1 J \oplus \nabla$.
\end{theorem}
\begin{pf} For the direct sense, if $\rho$ is well adapted, and $G$, $\Gamma$, $L$ attached to $\rho$, then $L$ is decomposable by Corollary~\ref{Ldecomp}, so $L = I_1 J \oplus \nabla$ and since $L$ is the Lie algebra of $\Gamma$, it satisfies the first five given conditions by
Prop.~\ref{structuredecomposable}. Moreover $L$ satisfies the last three conditions by Corollary~\ref{Amreducible}.

\par \bigskip

Conversely, if $L=I_1 J \oplus \nabla$ with $I_1$ and $\nabla$ satisfying the eights conditions above, then by \ref{structureLZ2Z1} to \ref{structureLZ2Z5} and Prop.~\ref{structuredecomposable}, $L$ is a Lie subring of $(\rad R)^0$ and $\Gamma := \Theta^{-1} (L)$ is a closed subgroup of $SR^1$ whose Lie algebra is $L$. Let $\Gb$ be any diagonal subgroup of $\GL_2(\F)$ containing at least one non-scalar matrix. It is clear that the conjugation by the subgroup $s(\Gb)$ of $R^\ast$ normalizes $L$, hence $\Gamma$. We can thus form the closed subgroup $G := \Gamma s(\Gb)$ of 
$SR^1$, a semi-direct product of $s(\Gb)$ by $\Gamma$. The composition $G \rightarrow s(\Gb) \simeq \Gb \subset \GL_2(\F)$ is
a representation $\rhob: G \rightarrow \GL_2(\F)$ which is the sum of two distinct characters and whose projective image has order $2$.  

The restriction $(t,d)$ to $G$ of the maps $(\tr,\det)$ on $R$ is a pseudo-representation over $G$. We claim that $(G,\rhob,t,d)$ is an admissible pseudo-deformation. The only condition that is not trivial to check is that the closed $W(\F)$-algebra generated by $\tr(G)$ is $A$.
Let us call this $W(\F)$-subalgebra by $\tilde A$. Since $\tr(1)=2$, $\tilde A$ contains $W(\F)1$.  Since $\tr(J \gamma) =  I_1$, $\tilde A$ contains $W(\F) I_1$. Also $\tilde A$ contains $\tr(\Gamma)$, hence it contains the closed sub-pseudoring generated by the element $\tr(\gamma)-2$, $\gamma \in \Gamma$, that is, it contains $P$ by Cor.~\ref{corPpseudo}. Therefore $\tilde A$ contains $W(\F)1 + W(\F) I_1 + W(\F) P = W(\F)1 + W(\F) I_1 + W(\F) I_1^2 + W(\F) \tr(\nabla^2)$, which is $A$ by \ref{structureLZ2Z6}. This concludes the proof of the claim that $(G,\rhob,t,d)$ is an admissible pseudo-deformation.

Let us define $\rho$ as the inclusion map $G \hookrightarrow R^\ast$. Then $\tr \rho=t$, $\det \rho=d$. We claim that $A \rho(G)=AG$ is the full algebra $R$. By Lemma~\ref{imagesubGma}, we know that $A \rho(G)=AG$ is a sub-$A$-GMA $\mat{A & B' \\ C' & A}$ of $R$, where $B'$ is a sub-$A$-module of $B$ and $C'$ a sub-$A$-module of $C$. By definition, $B'$ contains $B_1$ and $C'$ contains $C_1$, so  \ref{structureLZ2Z7} and \ref{structureLZ2Z8} imply that $B'=B$ and $C'=C$, so $A G = R$. It follows that $\rho: G \rightarrow R^\ast$ is a $(t,d)$-representation.
It is clear that the Lie algebra attached to $\rho$ is $L$, which proves the converse part of the theorem.
\end{pf}

\subsection{The structure of $L$ when the projective image of $\rhob$ is cyclic of order $>2$}

That is, $\rhob = \chi_1 \oplus \chi_2$ with $\chi_1^2 \neq \chi_2^2$.
In this case, we shall only determine the structure of the Lie algebra $W(\F_q)L$ where $\F_q$ is a large enough subfield
of $\F$.

\begin{theorem}  \label{structureLcyclic}
  Let $(\Pi,\rhob,d,t)$ be an admissible pseudo-deformation such that the projective image of $\rhob$ is cyclic of order $m>2$ and let $\rho: \Pi \rightarrow R^\ast$, $R = \mat{A & B \\ C & A}$ a  well adapted $(t,d)$-representation. Let $\F_q$ be any subfield of $\F$ such that 
  $\gcd(m,q-1)>2$ (a condition always satisfied when $\F_q=\F$).
  
  Then there exists a closed $W(\F_q)$-submodule $\tilde I_1$ of $\m$, and closed $W(\F_q)$-submodules $\tilde B_1$ of $B$ and $\tilde C_1$ of $C$ such that $W(\F_q) L =  \mato{\tilde I_1 & \tilde B_1 \\ \tilde C_1 & \tilde I_1}$
 and
\begin{num}  $\tilde B_1 \tilde C_1 \subset \tilde I_1$ \end{num}
\begin{num} $\tilde I_1^3 \subset \tilde I_1$, \end{num} 
\begin{num} $W(\F) 1 + W(\F) \tilde  I_1+ W(\F) \tilde I_1^2 = A$ \end{num} 
\begin{num} $W(\F) \tilde B_1=B$ and $W(\F) \tilde C_1=C$. \end{num}
 
Conversely, if $\tilde I_1$, $\tilde B_1$, $\tilde C_1$ are $W(\F)$-submodules of $\m$ satisfying those three conditions, and 
$L = \mato{\tilde I_1 & \tilde B_1 \\ \tilde C_1 & \tilde I_1}$, then there exists an admissible pseudo-deformation $(\Pi,t,d,\rho)$ such that the projective image of $\rhob$ is cyclic of order $>2$ and a $(t,d)$-representation $\rho: \Pi \rightarrow R^\ast$ such that the Lie algebra attached to $\rho$ is $L = W(\F) L$. 
\end{theorem}
\begin{pf} Let $g_0 \in \Pi$ be such that $\rhob(g_0)$ generates the group $\Gb$ modulo scalar matrices, and let $\lambda_0,\mu_0$ be the eigenvalues of $\rhob(g_0)$. Since the group $\Gb$ modulo scalar matrices has order $>2$, one has $\lambda_0/\mu_0 \neq \pm 1$.
By Cor. ~\ref{Ldecomp}, $L$ is decomposable, so we can write $L = I_1 J \oplus \nabla$ as usual, and by Cor. \ref{Lstrongdecomp} (applied with $n=1$), $W(\F) L$ is even strongly decomposable,  and we can  write $W(\F_q) L =  \mato{W(\F_q) I_1 & W(\F_q) B_1 \\ W(\F_q) C_1 & W(\F_q) I_1}$. Let us set $\tilde I_1 := W(\F_q) I_1$, $\tilde B_1= W(\F_q) B_1$, $\tilde C_1= W(\F_q) C_1$.
By Prop.~\ref{structuredecomposable}, one has $[\nabla,\nabla] \subset I_1 J$, which gives after taking the $W(\F_q)$-modules generated 
by the two terms of that inclusion, 
$\tilde B_1 \tilde C_1 \subset \tilde I_1$;
one has $I_1[J,\nabla] \subset \nabla$ which gives similarly $\tilde I_1 \tilde B_1 \subset \tilde B_1$, $\tilde I_1 \tilde C_1 \subset \tilde C_1$;
and $I_1^3 \subset I_1$, which gives $\tilde I_1^3 \subset \tilde I_1$. By Prop.~\ref{Amreducible}, $W(\F) 1 \oplus W(\F) I_1 \oplus W(\F) I_1^2 \oplus W(\F) \tilde B_1 \tilde C_1 = A$, and since $\tilde B_1 \tilde C_1 \subset \tilde I_1$, one has simply $W(\F) 1 \oplus W(\F) \tilde I_1 \oplus W(\F) \tilde I_1^2 = A$. Since $W(\F) B_1$ is stable by $W(\F) \tilde I_1$, it is stable by $A$, i.e. an $A$-module. But by Prop.~\ref{Amreducible}, the $A$-module generated by $B_1$, or by $W(\F) B_1$ is $B$. Therefore $W(\F) B_1=B$ and similarly $W(\F) C_1 =C$.
This completes the proof of the direct sense of the theorem.

Conversely, suppose $L = \mato{\tilde I_1 & \tilde B_1 \\ \tilde C_1 & \tilde I_1}$ satisfying the four given conditions. Then by Prop.~\ref{conditionstrongdec}, $\Gamma:= \Theta^{-1}(L)$ is a closed subgroup of $SR^1$ of Lie algebra $L$ (note that the condition~\ref{I1Jnablastrong}, i.e. $\tilde I_1 B \subset B$ and $\tilde I_1 C \subset C$, is automatically satisfied since $B$ and $C$ are $A$-modules). Let $\Gb$ be any group of diagonal matrices in $\GL_2(\F)$ whose quotient modulo scalar matrices is of order $>2$. It is clear that $\Gb$ normalizes $L$, hence $\Gamma$, and we can form a subgroup $G := \Gamma \Gb$ of $R^\ast$. Then the construction of $\rhob$, $t$, $d$, $\rho$ and the end of the proof of the converse is exactly as in the preceding theorem, so we leave details to the reader.
\end{pf}

\subsection{The structure of $L$ when the projective image of $\rhob$ is $\Z/2\Z \times \Z/2\Z$} 

\begin{theorem} \label{structureLZ2ZZ2Z} Let $(\Pi,\rhob,d,t)$ be an admissible pseudo-deformation such that the projective image of $\rhob$ is $\Z/2\Z \times \Z/2\Z$ and let $\rho: \Pi \rightarrow \GL_2(A)$ a well adapted $(t,d)$-representation.


There exists a closed subgroup $I_1$ of $\m$, and a closed subgroup $\nabla$ of $\mat{0 & \m \\ \m & 0}$ such that $$L = I_1 J \oplus \nabla$$ and
\begin{num}  \label{structureLZ2ZZ2Z1} $[\nabla,\nabla] \subset I_1 J$,\end{num}
\begin{num}  $I_1[ J,\nabla] \subset \nabla$, \end{num}
\begin{num} $\tr(\nabla^2) I_1 \subset I_1$, \end{num}
\begin{num}  $\tr(\nabla^2) \nabla \subset \nabla$, \end{num}
\begin{num}  \label{structureLZ2ZZ2Z5} $I_1^3 \subset I_1$, \end{num} 
\begin{num} \label{invlambda} There exists $\lambda \in s(\F_p^\ast)$ such that $\nabla$ is invariant by $\mat{ 0 & b \\ c & 0} \mapsto \mat{ 0 & \lambda c \\ b & 0}$ \end{num}
\begin{num}  \label{structureLZ2ZZ2Z7} $W(\F) 1 + W(\F) I_1+ W(\F) I_1^2 + W(\F) \tr(\nabla^2) + W(\F) B_1= A$ \end{num} 
\noindent

Conversely,  if $I_1$ is any closed subgroup of $\m$, and $\nabla$ any closed subgroup  of $\mat{0 & \m \\ \m & 0}$ satisfying the seven conditions above, then there exists an admissible pseudo-deformation
$(\Pi,\rhob,t,d)$ with $\rhob$ of projective image $\Z/2\Z \times \Z/2\Z$, and a $(t,d)$-representation $\rho : \Pi \rightarrow \GL_2(A)$ such that the Lie algebra attached to $\rho$ is $L= I_1 J \oplus \nabla$.
\end{theorem}
\begin{pf}
For the direct sense, if $\rho$ is well adapted, and $G$, $\Gamma$, $L$ attached to $\rho$, then $L$ is decomposable by Corollary~\ref{Ldecomp}, so $L = I_1 J \oplus \nabla$ and since $L$ is the Lie algebra of $\Gamma$, it satisfies  conditions  \ref{structureLZ2ZZ2Z1} to \ref{structureLZ2ZZ2Z5} by
Prop.~\ref{structuredecomposable}. Moreover $L$ satisfies conditions \ref{invlambda} and \ref{structureLZ2ZZ2Z7} by Corollary~\ref{Amdihedral}.

Conversely, if $L=I_1 J \oplus \nabla$ with $I_1$ and $\nabla$ satisfying the seven conditions above, then by Prop.~\ref{structuredecomposable} $L$ is a Lie subring of $\M_2(\m)$ and $\Gamma := \Theta^{-1} (L)$ is a closed subgroup of $SR^1$ whose Lie algebra is $L$. Let $\Gb$ be the subgroup of $\GL_2(\F)$ containing all matrices of the form $\smat{x & 0 \\ 0 & x}$, $\smat{x & 0 \\ 0 & -x}$, $\smat{0 & x \\ \lambda x & 0}$, $\smat{0 & x \\ -\lambda x & 0}$. This is a subgroup of order $4 |\F^\ast|$ which contains the subgroup of scalar matrices $\F^\ast$ of $\GL_2(\F)$, and the quotient $\bar G / \F^\ast$ is $\Z/2\Z \times \Z/2\Z$. The Lie algebra $L$ is stable by conjugation by $\Gb$ by~\ref{invlambda}. Therefore, so is $\Gamma$, and we can define a closed subgroup $G := \Gamma \Gb$ of 
$\GL_2(A)$. We thus have a split exact sequence $1 \rightarrow \Gamma \rightarrow G \rightarrow \Gb \rightarrow 1$.
We define a representation $\rhob: G \rightarrow \GL_2(\F)$ by composing the
natural map $G \rightarrow \Gb$ with the inclusion $\Gb \rightarrow \GL_2(\F)$. It is clear that $\rhob$ has projective image isomorphic to $\Z/2\Z \times \Z/2\Z$. We define a pseudo-repreresentation $(t,d)$ on $G$ by  restricting the trace and determinant map on $\GL_2(A)$ to $G$. 

We claim that $(G,\rhob,t,d)$ is an admissible pseudo-deformation. We just need to check that the closed $W(\F)$-algebra $\tilde A$ generated by $t(G)$ is $A$. Since $\tilde A$ contains $t(\Gamma)$ ad $t(J \Gamma)$, we see as in the proof of Theorem~\ref{structureLcyclic} that $\tilde A$ contains $W(\F) 1 + W(\F) I_1+ W(\F) I_1^2 + W(\F) \tr(\nabla^2)$. Moreover $\tilde A$ contains $\tr\left(\smat{0 & 1\\ \lambda & 0} \Gamma\right)$ and  $\tr(\smat{0 & 1\\ -\lambda & 0} \Gamma)$. When $\gamma= \smat{a & b \\ c & d}$ runs in $\Gamma=\Theta^{-1}(L)$, $\smat{0 & b \\ c & 0}$ runs in $\nabla$. Thus for any $\smat{0 & b \\ c & 0} \in \nabla$,
 $\tr \left(\smat{0 & \lambda\\ 1 & 0}\gamma \right) = \lambda c + b$ and $\tr\left( \smat{0 & -\lambda\\ 1& 0} \gamma \right) = -\lambda c + b$ are in $\tilde A$, and therefore $b \in \tilde A$. Thus $\tilde A$ contains $W(\F) B_1$ as well. By condition \ref{structureLZ2ZZ2Z7}, $\tilde A=A$, and this proves the claim.
 
 Let $\rho: G \rightarrow \GL_2(A)$ be the inclusion map. The representation $\rho$ is of trace $t$ and determinant $d$, and $A \rho(G)=M_2(A)$ by Lemma~\ref{imagesubGma}. So $\rho$ is a $(t,d)$-representation. The image of $\rho$ is $G$, its intersection with $SR^1$ is $\Gamma$, and the Lie algebra of $\Gamma$ is $L$. This proves the converse part of the theorem.
\end{pf}

\subsection{The structure of $L$ when the projective image of $\rhob$ is a non-abelian dihedral group}

\begin{theorem}  \label{structureLdihedral}
  Let $(\Pi,\rhob,d,t)$ be an admissible pseudo-deformation such that the projective image of $\rhob$ is a non-abelian dihedral group of order $2m>4$, $\rho: \Pi \rightarrow \GL_2(A)$ a  well adapted $(t,d)$-representation. Let $\F_q$ be any subfield of $\F$ such that 
  $\gcd(m,q-1)>2$ (a condition always satisfied when $\F_q=\F$).

   Then there exist  closed $W(\F_q)$-submodules $\tilde I_1$ and $\tilde B_1$ of $\m$ such that $W(\F_q) L = \mato{\tilde I_1 & \tilde B_1 \\ \tilde  B_1 & \tilde I_1}$
 and
\begin{num}  $\tilde B_1^2 \subset \tilde I_1$ \end{num}
\begin{num} $ \tilde I_1 \tilde B_1 \subset \tilde B_1$. \end{num}
\begin{num} $ \tilde I_1^3 \subset \tilde I_1$, \end{num} 
\begin{num} $W(\F) 1 + W(\F) \tilde I_1 + W(\F) \tilde I_1^2 + W(\F) \tilde B_1 = A$ \end{num} 

Conversely, if $\tilde I_1$ and $\tilde B_1$ are $W(\F)$-submodules of $\m$ satisfying those four conditions, and $L =  \mato{\tilde I_1 & \tilde B_1 \\ \tilde B_1 & \tilde I_1}$, then there exists an admissible pseudo-deformation $(\Pi,t,d,\rho)$ such that the projective image of $\rhob$ is dihedral of order $>4$ and a $(t,d)$-representation $\rho: \Pi \rightarrow R^\ast$ such that the Lie algebra attached to $\rho$ is $L = W(\F) L$. 
\end{theorem}
\begin{pf} We show as in the proof of Theorem~\ref{structureLcyclic} than $W(\F_q) L$ is strongly decomposable,  and we can thus write with the usual notations $W(\F_q) L = W(\F_q) I_1 \oplus \mat{0 & W(\F_q) B_1 \\ W(\F_q) C_1 & 0}$. By Prop.~\ref{Amdihedral}, $B_1 =  C_1$.
Thus, setting $\tilde I_1 = W(\F_q) I_1$, $\tilde B_1 = W(\F_q) B_1$, one has $W(\F_q) L = \mato{\tilde I_1 & \tilde B_1 \\ \tilde  B_1 & \tilde I_1}$.
By Prop.~\ref{structuredecomposable}, one has $[\nabla,\nabla] \subset I_1 J$, which gives after taking the $W(\F_q)$-modules generated 
by the two terms of that inclusion,  $\tilde B_1 ^2 \subset \tilde I_1$;
one has $I_1[J,\nabla] \subset \nabla$ which gives similarly $\tilde I_1 \tilde B_1 \subset \tilde B_1$,
and $ I_1^3 \subset I_1$, which gives $\tilde I_1^3 \subset \tilde I_1$. By Prop.~\ref{Amreducible}, $W(\F) 1 + W(\F) \tilde I_1 + W(\F) \tilde I_1^2 + W(\F) \tilde B_1^2  \tilde B_1= A$, and since $\tilde B_1^2 \subset \tilde I_1$, one has more simply $W(\F) 1 + W(\F) \tilde I_1 + W(\F) \tilde I_1^2 + W(\F) \tilde B_1 = A$.  This completes the proof of the direct sense of the theorem.

Conversely, suppose $L = \mato{\tilde I_1 & \tilde B_1 \\ \tilde B_1 & \tilde I_1}$ satisfying the four given conditions, then by Prop.~\ref{conditionstrongdec}, $\Gamma:= \Theta^{-1}(L)$ is a closed subgroup of $SR^1$ of Lie algebra $L$. Let $\Gb$ be for instance the group of all diagonal and anti-diagonal matrices in $\GL_2(\F)$. It is clear that $\Gb$ normalizes $L$, hence $\Gamma$, and we can form a subgroup $G := \Gamma \Gb$ of $R^\ast$. Then the construction of $\rhob$, $t$, $d$, $\rho$ and the end of the proof of the converse is exactly as in the preceding theorem, so we leave details to the reader.
\end{pf}

\subsection{Structure of $L$ in the large and exceptional image case}
 
\subsubsection{Results} 
 
This is the simplest case insofar as the description of $L$ is concerned, but the case where the proofs are  the hardest.

\begin{theorem} \label{structureLlargeimage} Let $(\Pi,\rhob,t,d)$ be an admissible pseudo-deformation. We assume that $\rhob$ is either of the large image type or of the exceptional type. 

Let $\F_q$ be a subfield of $\F$. If $\rhob$ is octahedral (resp. tetrahedral, resp. icosahedral), we assume that $\F_q$ contains cubic roots of unity (resp. either cubic or quartic roots of unity, resp. quintic roots of unity). We put no condition on $\F_q$ when $\rhob$ has large image (we can take $\F_q=\F_p$ in this case).

Then there exists a $(t,d)$-representation $\rho$, and a closed $W(\F_q)$-submodule $\tilde I_1$ of $\m$ such that
$$W(\F_q) L =  \mato{ \tilde I_1 & \tilde I_1 \\ \tilde I_1 & \tilde I_1},$$
and 
\begin{num} $\tilde I_1^2 \subset \tilde I_1$.\end{num}
\begin{num} $W(\F) \tilde I_1 = \m$. \end{num}
\end{theorem}

\begin{remark} Note that if we take $\F_q=\F$, and $\F$ large enough as we always do, the hypothesis of the theorem are obviously satisfied.
Thus, the theorem describes the structure of $W(\F)L$ for $\F$ large enough.
\end{remark}
\begin{cor}  \label{corstructureLlargeimage} With the same notation as in the above theorem, one has $\Gamma = \Theta^{-1}(L)$.
In the case where $\F_q=\F_p$, $\Gamma$ is precisely the group of matrices $\mat{a & b \\ c & d}$ in $\SL_2(A)$ such that $a \in 1  + \tilde I_1$,
$b \in \tilde I_1$, $c \in \tilde I_1$, $d \in 1+ \tilde I_1$.
\end{cor}
\begin{pf} This follows from the preceding theorem and Prop.~\ref{exampleliesimple}. \end{pf}

\begin{remark} In the appendix of \cite{MW}, Boston proves that if $G$ is a closed subgroup of $\SL_2(A)$ such that the image of $G$ in $\SL_2(A/\m^2)$ is $\SL_2(A/\m^2)$, then $G=\SL_2(A)$. This results follow easily from our classification result Theorem~\ref{structureLlargeimage} (indeed, we are in the large image case so we can take $\F_q=\F_p$, and the hypothesis implies that $\tilde I_1$ maps surjectively to $\m/\m^2$, thus is $\m$, which implies $\Gamma = \SL_2^1(A)$ by the corollary and $G=\SL^2(A)$). In does not seem that Boston's method generalizes to the other cases covered 
by Theorem~\ref{structureLlargeimage}.
\end{remark}

\subsubsection{Proof of Theorem~\ref{structureLlargeimage}}

We fix $(\Pi,\rhob,t,d)$ an admissible pseudo-deformation, and we assume that $\rhob$ is either of the large image type or of the exceptional type. We call $Z$ the subgroup of scalar matrices in $\GL_2(\F)$, isomorphic to $\F^\ast$.
 
 \begin{lemma} \label{exitencerholargeimage} There exists a $(t,d)$-representation $\rho$ such that:
\begin{itemize}
\item if $\rhob$ is octahedral, then $\Gb \subset Z \GL_2(\F_q)$ and there exists an element $g_0 \in \Pi$ 
such that $\rho(g_0) = \smat{s(\lambda_0) & 0 \\ 0 & s(\mu_0)}$
with $\lambda_0,\mu_0 \in \F^\ast$, $(\lambda_0/\mu_0)^3 = 1$, $\lambda_0 \neq \mu_0$;
\item if $\rhob$ is tetrahedral, then $\Gb \subset Z \GL_2(\F_q)$ and there exists an element $g_0 \in \Pi$ such that
$\rho(g_0) =  \smat{s(\lambda_0) & 0 \\ 0 & s(\mu_0)}$
with  $\lambda_0,\mu_0 \in \F^\ast$, $(\lambda_0/\mu_0)^3 = 1$ or $(\lambda_0/\mu_0)^4 =1$, and $\lambda_0^2 \neq \mu_0^2$;
\item if $\rhob$ is icosahedral, then $\Gb \subset Z \GL_2(\F_q)$ and there exists an element $g_0 \in \Pi$ 
such that
$\rho(g_0) = \smat{s(\lambda_0) & 0 \\ 0 & s(\mu_0)}$
with  $\lambda_0,\mu_0 \in \F^\ast$, $(\lambda_0/\mu_0)^5 = 1$, $\lambda_0 \neq \mu_0$;
\item if $\rhob$ has large image but the projective image of $\rhob$ is not isomorphic to $\PSL_2(\F_3)$ or $\PGL_2(\F_3)$, then $\SL_2(\F_p) \subset \Gb$ and there exists an element $g_0 \in \Pi$ such that 
$\rho(g_0) = \smat{ s(\lambda_0) & 0 \\ 0 & s(\mu_0)}$, $\lambda_0^2 \neq \mu_0^2$
\end{itemize}
\end{lemma}
\begin{pf}
The image of $\Gb$ on $\PGL_2(\F)$ contains an element of order $3$ in the octahedral case (a $3$-cycle in $A_4$),
an element of order 3 and of order 4 in the tetrahedral case (a $3$-cycle and a $4$-cycle in $S_4)$, an element of order $5$
in the icosahedral case (a $5$-cycle in $A_5$), and an element of order $>2$ in the large image case. Choosing an element $g_0$ such that $\rhob(g_0)$ maps to that element, we can diagonalize $\rhob(g_0)$ and write $\rhob(g_0) = \mat{\lambda_0 & 0 \\ 0 & \mu_0}$.
Choosing a $\rho$ adapted to $(g_0, \lambda_0,\mu_0)$ ensures that, in each case, the condition regarding $\rho(g_0)$.
For such a $\rho$, $\rhob(g_0) = \mat{\lambda_0 & 0 \\ 0 & \mu_0}$. On the other hand, we know that in the 
conjugacy class of $\rhob$ there is a representation $\rhob'$ satisfying $\rhob'(\Pi) \subset Z \GL_2(\F_q)$ in the exceptional cases and 
$\SL_2(\F_p) \subset \Gb$ in the large image case, and $\rhob'(g_0) = \mat{\lambda_0 & 0 \\ 0 & \mu_0}$. (This is because, in the exceptional case, there is a conjugate of $\rhob$ whose projective image is defined over $\F_q$, and after a base change over $\F_q$,
we may suppose that $\rhob'(g_0)$ is diagonal). The agreement of $\rhob$ and $\rhob'$ on $g_0$ implies that they are conjugate through a diagonal matrix. Conjugating $\rho$ by a diagonal lift of that diagonal matrix doesn't affect the condition on $\rho(g_0)$ but ensures that $\Gb = \rho'(\Pi)$ satisfies the required condition.
\end{pf}

Since in any case the eigenvalues $\lambda_0$ and $\mu_0$ of $\rhob(g_0)$ have distinct squares, Prop.~\ref{Lstrongdecomp} applies and ensures that $W(\F_q) L$ is strongly decomposable. That is, there exists three $W(\F_q)$-submodules of $A$, $\tilde I_1$, $\tilde B_1$ and $\tilde C_1$ such that 
$$W(\F_q) L = \tilde I_1 J \oplus \mat{0 & \tilde B_1 \\ \tilde C_1 & 0}.$$

\begin{lemma} \label{conjbarg} Let $\bar g = \smat{\bar \alpha & \bar \beta \\ \bar \gamma & \bar \delta} \in \GL_2(\F_q) \cap (\Gb Z)$. Then there exists a lift $g = \smat{\alpha & \beta \\ \gamma & \delta} \in \GL_2(A)$ of $\bar g$
such that
\begin{num} $\alpha^2 \tilde B_1 \subset \tilde B_1$, $\gamma^2 \tilde B_1 \subset \tilde C_1$, $\alpha \gamma \tilde B_1 \subset \tilde I_1$. \end{num}
\begin{num} $\beta^2 \tilde C_1 \subset \tilde B_1$, $\beta \delta \tilde C_1 \subset \tilde I_1$, $\delta^2 \tilde C \subset \tilde B_1$. \end{num}
\begin{num}  $\alpha \beta \tilde I_1 \subset \tilde B_1$,
$\gamma \delta \tilde I_1 \subset \tilde C_1$. \end{num}
\end{lemma}
\begin{pf} By hypothesis $\bar g = \bar g_1 \bar z$ with $\bar g_1 \in G$ and $\bar z$ a scalar matrix in $\F^\ast$. 
Let $g_1$ be any lift of $\bar g_1$ in $G$, $z = s(\bar z)$ which is a scalar matrix lifting $\bar z$, and set $g = g_1 z \in \GL_2(A)$ which is a lift of $\bar g$. Then $g  L g^{-1} = g_1  L g_1^{-1} =  L$ since $z$ is scalar and $g_1 \in G$. Moreover, if $g=\smat{\alpha & \beta \\ \gamma & \delta}$, set $g' = \smat{\delta & -\beta \\ - \gamma & \alpha} = g^{-1} \det(g)$. 
Note that $\det(g) = \det(g_1) \det(z) = s(\det(\bar g_1)) s(\det(\bar z))$ by \ref{dconstant}, so $\det(g) = s(\det \bar g) \in W(\F_q)^\ast.$ Then $g W(\F_q) L g' = g W(\F_q) L g^{-1} = W(\F_q) L$ since multiplication by $\det(g)^{-1}$ stabilizes $W(\F_q) L$.

The first line of the lemma then follows from the computation $g \smat{0 & 1 \\ 0 & 0} g' = \smat{-\alpha \gamma & \alpha^2 \\ - \gamma^2 & \alpha \gamma}$, the second line from $g \smat{0 & 0 \\ 1 & 0} g' = \smat{\beta \delta & -\beta^2 \\ \delta^2 & -\beta \delta}$, and the last line from $g J g^{-1} = \smat{* & -2 \alpha \beta \\ -2 \gamma \delta & *}$.
\end{pf}

\begin{lemma} \label{exbarg} There exists an element $\bar g = \smat{\bar \alpha & \bar \beta \\ \bar \gamma & \bar \delta} \in \GL_2(\F_q) \cap (\Gb Z)$
such that $\bar \alpha \bar \beta \neq 0$ (resp. $\bar \alpha \bar \gamma \neq 0$, resp. $\bar \beta \bar \delta \neq 0$, resp. $\bar \gamma \bar \delta \neq 0$.)
\end{lemma}
\begin{pf} In the large image case, we can take  for instance $\bar g = \smat{ 1 & 1 \\ 1 & 2}  \in \Gb \cap \GL_2(\F_p)$.

We assume that we are in some of the exceptional cases. It suffices to find one matrix $g_1$ in $\Gb$ satisfying $\bar \alpha  \bar \beta \neq 0$ for then since $\Gb \subset GL_2(\F_q) Z$, a suitable product of $g_1$ by a scalar matrix will belong to $\GL_2(\F_q)$ and obviously will still satisfies the required condition.

If all matrices in $\PGb$ had $\bar \beta = 0$, then the representation $\rhob$ would be reducible. Among the matrices such that $\bar \beta \neq 0$, if there is one with $\bar \alpha \neq 0$, we are done. Otherwise, all matrices with $\bar \beta \neq 0$ are of the form $\smat{ 0 & \bar \beta \\ \bar \gamma & \bar \delta}$ and their square is $\smat{ \bar \beta \bar \gamma & \bar \beta \bar \delta \\ * & *}$. Now $\bar \gamma$ is not $0$ because the matrix is invertible, and if $\bar \delta \neq 0$ either, we are done. Otherwise, this means that all matrices with $\bar \beta \neq 0$ have 
both $\bar \alpha$ and $\bar \delta$ equal zero, that is are antidiagonal. But then it is easy to see that $\PGb$ is contained in the normalizer of the diagonal torus, a contradiction with the hypothesis that $\PGb$ is exceptional.

\end{pf}

 \begin{lemma} \label{xX} Let $X$ be a closed $W(\F_q)$-submodule of $A$. Let $x \in A^\ast$ whose image in $A/\m=\F$ lies in $\F_q$.
 \begin{itemize}
 \item If $x X \subset X$, then $x X =X$.  
 \item If $x^2 X = X$, then $x X = X$.
 \end{itemize}
 \end{lemma} 
 \begin{pf} Replacing $x$ by $x^{-1}$,  the hypothesis becomes $X \subset x X$, and the contention is still that $X = x X$.
 Then by induction $X \subset x X \subset \cdots \subset x^n X$ for all $n > 0$.
 Writing $x=s(\bar x)+m$ with $\bar x\in \F_q$ the reduction of $x$ and $m \in \m$, we get
 $x^{q^n} \equiv s(\bar x) \pmod{\m^{nv_p(q)+1}}$ and so $X \subset x X \subset  X + \m^{n v_p(q)+1}X $. Since $X$ is a closed subgroup, the intersection of all $X + \m^{n v_p(q)+1}X$ when $n \geq 1$ is $X$, and we get $X \subset x X \subset X$, as desired. This proves the first point.
 
 For the second point, not that if $x^2 X = X$, then $x^{2n} X = X$. Choosing a sequence of positive integers $n$ which converges to $1/2$ $p$-adically gives the result.
 \end{pf}

 \begin{lemma} There exists $x,y \in A^\ast$ such that the image of $x$ and $y$ in $A/\m=\F$ are in $\F_q$,  $\tilde B_1 = x \tilde I_1$, and $\tilde C_1 = y \tilde I_1$. 
 \end{lemma}
\begin{pf} Pick a matrix $\bar g$ as in Lemma~\ref{exbarg} such that $\bar \alpha \bar \gamma \neq 0$. By Lemma~\ref{conjbarg}, there is a lift $u \in A^\ast$ of $\bar \alpha \bar \gamma$ such that $u \tilde B_1 \subset \tilde I_1$. Also pick a matrix $\bar g'$ as in Lemma~\ref{exbarg} such that $\bar \alpha' \bar \beta' \neq 0$. By Lemma~\ref{conjbarg}, there is a lift $v \in A^\ast$ of $\bar \alpha \bar \beta$ such that $v \tilde I_1 \subset \tilde B_1$. Thus $u v \tilde I_1 \subset u \tilde B_1 \subset \tilde I_1$. The inclusion $uv \tilde I_1 \subset \tilde I_1$ is an equality by Lemma~\ref{xX} (note that the image of $uv$ in $A/\m$ is $\bar \alpha \bar \gamma \bar \alpha' \bar \beta' \in \F_q^\ast$). Therefore, $u \tilde B_1 = \tilde I_1$ and the first result follows with $x=u^{-1}$. The second is similar.
 \end{pf}
 
 Now let $t = \sqrt{xy^{-1} / s(xy^{-1})} \in A^\ast$. We check easily that $\frac{x}{t} = yt\, s(xy^{-1})$. We conjugate $\rho$ by the diagonal matrix $\smat{1 & 0 \\ 0 & t}$. This doesn't affect any of the property of $\rhob$ already stated, and doesn't change $\tilde I_1$ but change $\tilde B_1$ into $\frac{1}{t} \tilde B_1=\frac{x}{t} \tilde I_1$ and $\tilde C_1$ into $t \tilde C_1=yt \tilde C_1 = yt\, s(xy^{-1}) \tilde C_1$.
Replacing $x$ by $\frac{x}{t}$, we get:
\begin{num} There exists $x \in A^\ast$ such that the image of $x$ in $A/\m=\F$ is in $\F_q$, such that $\tilde B_1 = \tilde C_1 = x \tilde I_1$.
 \end{num}
 
 Now we again pick  a matrix $\bar g$ as in Lemma~\ref{exbarg} such that $\bar \alpha \bar \gamma \neq 0$. By Lemma~\ref{conjbarg}, for some lift $\alpha$, $\gamma$ of $\bar \alpha, \bar \gamma$, one has $\alpha^2 \tilde B_1 \subset \tilde B_1$
 and $\gamma^2 \tilde B_1 \subset \tilde C_1=\tilde B_1$. By the first point of Lemma~\ref{xX}, this means  $\alpha^2 \tilde B_1 =\tilde B_1$ and $\gamma^2 \tilde B_1 = \tilde B_1$, and by the second point  $\alpha \tilde B_1 =\tilde B_1$ and $\gamma \tilde B_1 = \tilde B_1$. Therefore $\alpha \gamma \tilde B_1 = \tilde B_1$. On the other hand, by Lemma~\ref{conjbarg}, $\alpha \gamma \tilde B_1 \subset \tilde I_1$, and thus $\tilde B_1 \subset \tilde I_1$. The converse inclusion $\tilde I_1 \subset \tilde B_1$ is proved similarly using a matrix with  $\bar \alpha \bar \beta \neq 0$. We have therefore proved:
 
 \begin{num} \label{B1C1I1L}$\tilde B_1 = \tilde C_1 = \tilde I_1$, $L= \mato{ \tilde I_1 & \tilde I_1\\ \tilde I_1 & \tilde I_1}$. \end{num}
 \noindent
 From \ref{I1Jnablastrong}, one has $\tilde I_1 \tilde B_1 \subset \tilde B_1$, that is 
 \begin{num} \label{I1pseudo} $\tilde I_1^2 \subset \tilde I_1$. \end{num}
   
\begin{prop} One has $W(\F) \tilde I_1 = \m$. \end{prop}
\begin{pf}  By \cite[Theorem 7.16(b)]{eisenbud}, it is enough to prove that the natural composed map $f: W(\F) \tilde I_1 \hookrightarrow \m \rightarrow \m/\m^2$ is surjective. To prove this, it is enough to prove that for each non-zero linear form $l : \m/\m^2 \rightarrow \F$, the composition
$l \circ f: W(\F) \tilde I_1 \rightarrow \F$ is surjective, which is the same as being non-zero. Such a linear form $l$ (geometrically, a tangent vector to the unique closed point of $\spec A$) induces a surjective morphism of rings $A \rightarrow A/\m^2 \rightarrow \F[\epsilon]$
where the second map sends $m \in \m/\m^2$ to $l(m) \epsilon$. We need to prove that the image of $\tilde I_1$ in that map is non zero.
By functoriality (see \S\ref{functorialitypink}), the image of $\tilde I_1$ in $\F[\epsilon]$ is the same as the $\tilde I_1$ obtained for the admissible pseudo-deformation
$(\Pi,\rho,t',d')$ over $\F[\epsilon]$, where $t',d'$ are $t,d$ composed with the map $A \rightarrow \F[\epsilon]$. 

In other words, we have reduced the proof of the proposition to the case $A=\F[\epsilon]$, and in this case we just have to prove that $\tilde I_1 \neq 0$. We proceed by contradiction. Assume $\tilde I_1 = 0$. Then by~\ref{B1C1I1L}, $L=0$, so $\Gamma \subset \Theta^{-1}(L)$ is the trivial group and the reduction map $G \rightarrow \bar G$ is an isomorphism. The morphism $r: \bar G \simeq G \subset \GL_2(A)$ is thus a deformation to $A = \F[\epsilon]$ of the tautological representation $\bar G \subset \GL_2(\F)$. Such deformation are parametrized
by $H^1(\bar G, V)$, where $V$ is the trace-zero adjoint representation of the tautological representation of $\bar G$, and this cohomology group is trivial by Prop.~\ref{H1adj0}. Therefore, the trace of $r$ is constant, that is $\tr(G) \subset \F$, in contradiction with the hypothesis~\ref{tgenerate} that $\tr(G)$ generates $A = \F[\epsilon]$ as an $\F$-algebra.
\end{pf} 
Together, this proposition, \ref{B1C1I1L} and \ref{I1pseudo} complete the proof of Theorem~\ref{structureLlargeimage}

 \section{Congruence-large image}  \label{sectioncongruencelarge}
 
 This section is not used in the rest of the paper. Its aim is to establish a connection between our results on the structure of the image of pseudo-deformation and a series of recent results by Hida \cite{hida}, Lang \cite{lang} and Conti-Iovita-Tilouine \cite{CIT} concerning the image of the Galois representation carried by certain $p$-adic families of modular forms, ordinary in the work of first two named authors, of positive slope for the last group. Our setting is more general as we work with families of $2$-dimensional representations of arbitrary pro-finite groups, over arbitrary noetherian compact local domain. The aim of this section is to show that we can obtain, in this general setting, results that are quite close (and sometimes stronger) to those proved for families of modular forms.
 
In all this section, $A$ is a compact noetherian local ring with maximal ideal $\m$ and residue field $\F$ finite of characteristic $p>2$. We also assume that $A$ is a domain, of fraction field $K$. 
 
 \subsection{The notion of congruence-large image}

 \begin{definition} Let $R$ be a GMA over $A$. If $I$ is an ideal of $A$, the {\it principal congruence subgroup} $\Gamma_R(I)$ of $I$ is the subgroup of $R^\ast$ defined as the kernel of the map $SR \rightarrow (R/IR)^\ast$. A closed subgroup of $R^\ast$ is called {\it a congruence subgroup} if it contains $\Gamma_R(I)$ for some non-zero ideal $I$ of $A$. 
 \end{definition}

By definition,  $\Gamma_R(I)$ is the set of matrices $\smat{a & b \\ c & d}$ in $R$ such that $a,d \equiv 1 \pmod{I}$, $b \in IB$, $c \in IC$ and $ad-bc=1$. When $R=M_2(A)$, we retrieve the usual notion of the group of matrices congruent to the identity modulo $I$.
 
\begin{lemma} \label{lemmacongruencesubgroup} Let $R$ be a topological GMA over $A$, $\Gamma$ a closed subgroup of $SR^1$. Then $\Gamma$ is a congruence subgroup if and only if $L(\Gamma)$ contains $\mato{I & I \\ I & I}$ for some non-zero ideal $I$ of $A$.
\end{lemma}
\begin{pf} First, a trivial computation gives $\Theta(\Gamma_R(I)) = \smato{I & I \\ I & I}$, hence $L(\Gamma_r(I))= \smato{I & I \\ I & I}$.
By Prop~\ref{exampleliesimple}, $\Theta^{-1}\left( \smato{I & I \\ I & I}\right)$ is the unique closed subgroup of $SR^1$ whose Lie algebra is $ \smato{I & I \\ I & I}$, hence $\Theta^{-1}\left( \smato{I & I \\ I & I}\right) = \Gamma_R(I)$ (this can also be obtained by a direct computation).

Let $\Gamma$ be a closed subgroup of $SR^1$. If $\Gamma$ contains $\Gamma_R(I)$, then $L=L(\Gamma)$ contains $\Theta(\Gamma_R(I))=\smato{I & I \\ I &I}$. Conversely, assume that $L(\Gamma)$ contains $\smato{I & I \\ I &I}$.  Then $L_2$ contains $\smato{I^2 & I^2 \\ I^2 & I^2}$ and $\Gamma = \Theta^{-1}(L_2)$ by Theorem~\ref{thmpink1}, so $\Gamma$ contains $\Gamma_R(I^2)$.
\end{pf}

\begin{lemma} \label{lemmacongruence} Let $\Pi$ be a group, $(t,d)$ a $2$-dimensional pseudo-representation of $\Pi$ over $A$ which is not the sum of two 
characters. Let $R_1$ and $R_2$ be two faithful GMA over $A$, finite-type as $A$-modules, and let $\rho_1: \Pi \rightarrow R_1^\ast$,
$\rho_2: \Pi \rightarrow R_2^\ast$ be two representations both of trace $t$ and determinant $d$. Then $\rho_1(\Pi)$ is a congruence subgroup of $R_1$ if and only if $\rho_2(\Pi)$ is a congruence subgroup of $R_2$.
\end{lemma}
\begin{pf} By~\ref{basicGMA2}, we can assume that both $R_1$ and $R_2$ are sub-algebras of $M_2(K)$.  Seen as representations over $K$, $\rho_1$ and $\rho_2$ have the same trace and determinant, hence are conjugate. Let $g \in \GL_2(K)$ such that $\rho_2 = g \rho_1 g^{-1}$.  Since $R_1$ and $R_2$ are of finite type, there exists $z \in A-\{0\}$ such that $z g R_1 g^{-1} \subset R_2$.

If  $\gamma - 1 \in Iz R_1$, we have
 $ g (\gamma-1) g^{-1} \in Iz gR_1g^{-1} \subset I R_2$, and hence $g \gamma g^{-1} \in 1+ I R_2 \subset R_2$, so $g (\gamma-1) g^{-1} \in \Gamma_{R_2}(I)$. Therefore, $ \Gamma_{R_1}(Iz) \subset g^{-1} \Gamma_{R_2}(I) g $, and it follows that if $\rho_2(\Pi)$ contains a congruence subgroup of $R_2^\ast$, $\rho_1(\Pi)=g^{-1} \rho_2(\Pi) g$ contains a congruence subgroup of $R_1^\ast$.
\end{pf}

\begin{definition}  We say that an  two-dimensional pseudo-representation $(t,d)$ of a group $\Pi$ over $A$ has {\it congruence-large image} if for one (equivalently for any) representation $\rho: \Pi \rightarrow R^\ast$, with $R$ a faithful finite-type GMA over $A$, such that $\tr \rho = t $ and $\det \rho =d$, $\rho(\Pi)$ is a congruence subgroup of $R^\ast$.
\end{definition}

 \subsection{Sufficient conditions for a congruence-large image}
 
\begin{definition} \label{defregular} We say that a representation $\rhob: \Pi \rightarrow \GL_2(\F)$ is {\it regular} if there exists an element $g_0$ in $\Pi$ such that $\rhob(g_0)$ is diagonalizable of eigenvalues $\lambda$ and $\mu$ in $\F_p^\ast$, with $\lambda^2 \neq \mu^2$.
\end{definition} 
\begin{remark}
If $\rhob$ is regular, it has an element of order $>2$ in its projective image, which therefore cannot be cyclic of order 2, or dihedral of order 4.
In the other cases (cyclic of order $>2$, dihedral of order $>4$, large or exceptional), there exist many regular representations, for instance all that have $\F_p$ as field of definition.

The notion of regularity is related to the notion of an {\it $H$-regular representation} of Lang (\cite{lang}) and of an {\it $(H,\Z_p)$-regular representation} of Conti-Iovita-Tilouine of \cite{CIT}. Let us recall
that $H$-regular means that $H$ is a subgroup such that $\rhob_{|H}$ is reducible, and there is an element $g_0 \in H$ which is diagonalizable with distinct eigenvalues $\lambda,\mu$, while $(H,\Z_p)$-regular requires in addition that $\lambda^2 \neq \mu^2$ and $\lambda,\mu \in \F_p$. It is obvious that $(H,\Z_p)$-regular (for any $H$) implies regular in our sense, while regular implies $\Pi$-regular, but not in general $H$-regular for a proper subgroup $H$ of $\Pi$.
\end{remark}

 \begin{theorem} \label{thmconglarge} Assume that $A$ is a domain. Let $(\Pi,\rhob,t,d)$ be an admissible pseudo-deformation such that $\rhob$ is regular. Moreover, we assume that
 \begin{itemize}
 \item If $\rhob$ is reducible, $t$ is not the sum of two continuous characters $\Pi \rightarrow A^\ast$.
 \item If $\rhob$ is dihedral, then if $\Pi'$ is the unique subgroup of index $2$ of $\Pi$ such that $\rhob(\Pi')$ is abelian,
 $t_{|\Pi'}$ is not the sum of two characters.
 \end{itemize}
 
 Then there exists a subring $A_0$ of $A$, which is a complact noetherian local ring of maximal ideal $\m \cap A_0$, and an open subgroup $\Pi_0$  of $\Pi$, containing $\ker \rhob$, such that
 \begin{itemize}
\item $t(\Pi_0) \subset A_0$, $d(\Pi_0) \subset A_0^\ast$.
\item  $(\Pi_0,\rhob_{|\Pi_0},t_{|\Pi_0},d_{|\Pi_0})$ is an admissible pseudo-deformation over $A_0$, and has congruence-large image.
\end{itemize}
 \end{theorem}
 \begin{pf}
We choose a $g_0 \in \Pi$ as in the definition~\ref{defregular} and a $(t,d)$-representation $\rho : \Pi \rightarrow R^\ast$ with $R= \mat{A & B \\ C & D}$ adapted to $g_0$. In particular, if $D_0$ denotes the subgroup of $\rhob(\Pi)$ generated by $\rhob(g_0)$, 
then $D_0$ is a group of diagonal matrices and $s(D_0) \subset G$ by Theorem~\ref{diagonallift}. We write $\rhob(g_0)=\smat{\lambda & 0 \\ 0 & \mu}$.

By Cor.~\ref{Lstrongdecomp}, $L$ is strongly decomposable. We write $ L \mato{I_1 &  B_1 \\  C_1 & I_1}$ with $I_1$, $ B_1$ and $ C_1$ closed subgroups of $A$, $B$ and $C$ respectively.

We define $$A_0 := \Z_p +   I_1 + I_1^2.$$
 By \ref{I13strong}, $A_0$ is a subring of $A$, and it is clearly a compact local ring of maximal ideal $\m_0 = p \Z_p + I_1 + I_1^2 = \m \cap A_0$. By \ref{I1Jnablastrong}, both $\tilde B_1$ and $\tilde C_1$ are $A_0$-modules. 
 
We define $$\Pi_0 = \rhob^{-1}(D_0).$$  
This is obviously a subgroup of finite index in $\Pi$, containing $\ker \rhob$. The restriction of $\rhob$ to $\Pi_0$ is a reducible representation, sum of two distinct characters. 

We claim that the closed $W(\F_0)$-subring of $A$ generated by $t(\Pi_0)$ is $A_0$. Indeed, let us call $A'_0$ that subring. Any element of $\Pi_0$ can be written $s(d) \gamma$, with $d \in D_0 $ and $\gamma \in \Gamma \subset \Theta^{-1}(L)$, and thus has trace in $\Z_p + I_1 + P = \Z_p + I_1 + I_1^2$ (by~\ref{Pstrongdec} and~\ref{nablanablastrong}). Thus we see that $t(\Pi_0) \subset A_0$, hence $A'_0 \subset A_0$. On the other hand, $A'_0$ contains $\Z_p$ by definition. It therefore contains $\tr(\gamma)-2$ for every $\gamma \in \Gamma$, hence it contains $P$ by Cor.~\ref{corPpseudo}. And it contains $\tr(s(g_0)^n \Gamma)$ for any $n$, hence $I_1$. Thus $A_0=A'_0$.

It follows easily that $(\Pi_0,\rhob_{|\Pi_0},t_{|\Pi_0},d_{|\Pi_0})$ is an admissible pseudo-deformation over $A_0$.
By Cor.~\ref{admissiblenoetherian}, $A_0$ is a noetherian ring.

We define $R_0$ as the $A_0$-sub-GMA $\mat{A_0 & B_1 \\ C_1 & A_0}$ of $R$. Since this is a sub-GMA of $M_2(K)$, $R_0$ is faithfull provided that $B_1 \neq 0$ and $C_1 \neq 0$, and this follows from the hypothesis made on $(t,d)$. One has clearly $\rho(\Pi_0) \subset R_0^\ast$. Moreover $\rho(\Pi_0)$ generates $R_0$ as an $A_0$-module, since clearly the $s(g_0)^n$ generates the subring of diagonal matrices $\mat{A_0 & 0 \\ 0 & A_0}$ of $R_0$, and $\rho(\Pi_0)$ contains $\Gamma$, whose projection on anti-diagonal matrices topologically generates as a group, hence as an $A$-module, $\nabla = \mat{0 & B_1 \\ C_1 & A_0}$. Thus, the restriction $\rho_{|\Pi_0}$ of $\rho$ to $\Pi_0$ is a $(t_{|\Pi_0},d_{|\Pi_0})$-representation. 

Its image $\rho(\Pi_0)$ contains $\Gamma$, hence also $\Gamma_2 = \Theta^{-1}(L_2)$. From the description of $L$, it follows that $L_2 = \mat{ B_1 C_1 & I_1 B_1 \\ I_1 C_1 & B_1 C_1}$. Since $I_1$ contains $B_1 C_1$, $L_2$ contains $B_1 C_1 J \oplus \mat{ 0 &  B_1^2 C_1 \\ B_1 C_1^2 & 0} \supset B_1 C_1 R_0^0$, and it follows that the image of $\rho_{|\Pi_0}$ contains the congruence subgroup 
$\Gamma_{R_0}(B_1C_1)$.
 \end{pf}
 
 \begin{remark} With the notation of the preceding theorem and its proof, let $K$ be the fraction field of $A$, and $K_0$ be the fraction field of $A_0$. The representation $\rho: \Pi \rightarrow R^\ast$ induces a representation $\rho_K: \Pi \rightarrow \GL_2(K)$ since $R \otimes_A K = M_2(K)$. Similarly, $\rho_{\Pi_0}: \Pi_0 \rightarrow R_0^\ast$ induces a representation $\rho_{K_0}: \Pi_0 \rightarrow \GL_2(K_0)$. 
 The representations $\rho_{K_0}$ and $\rho_K$ have the same trace and determinant on $\Pi_0$. Therefore there exists $g \in \GL_2(K)$ such that $g \rho_K g^{-1} = \rho_{K_0}$ on $\Pi_0$. The conclusion of our theorem implies that $\rho_K(\Pi_0)$ contains $1+JR_0$
 for some non-zero ideal $J$. It follows from  Lemma~\ref{lemmacongruence} that $\rho_K(\Pi_0)$ contains $1+ J' M_2(A_0)$ for some ideal $J'$. Hence $g \rho(\Pi_0) g^{-1}$ contains the congruence subgroup $\Gamma_{M_2(A_0)}(J')$.
 
 This is the way the conclusion of the main theorem of Lang \cite[Theorem 2.4]{lang} is stated, as well as the main theorem of \cite{CIT}. 
 
 On the other hand, the hypotheses of Lang are that $\Pi=G_\Q$, $A$ a local domain finite over the Iwasawa algebra $\Z_p[[T]]$,
 $(t,d)$ the pseudo-representation carried by a Hida's family which is residually absolutely irreducible -- a very special case of the situation we are studying. She assumes in addition that the family is not CM, an hypothesis which is equivalent (under other running assumptions) to our assumption that $t_{|\Pi'}$ is not the sum of two characters. Finally she is assuming that $(t,d)$ is $\Pi_0$-regular, an hypothesis which imply our regularity assumption (it allows, it seems, for some $\rhob$ with projective image dihedral of order 4), nor it implied by ours.
 
To summarize Theorem~\ref{thmconglarge} implies the congruence-large image result of \cite[Theorem 2.4]{lang} in many casesthough not in all cases, and it implies  the congruence-large image result of \cite{CIT} in all cases.

In the references \cite{lang} and \cite{CIT}, the congruence-large image result are made more precise by an explicit description of the subring $A_0$ of $A$ and the subgroup $\Pi_0$ of $\Pi$, in terms of the {\it conjugate self-twist} of $(t,d)$ (see \cite[definition 2.1]{lang}). Our method also gives an explicit description of $\Pi_0$ and $A_0$, though a different one. It would be interesting to compare these descriptions.
 \end{remark}

 \section{The essential submodule attached to an admissible pseudo-deformation}
 
 \label{sectionessential}
 
In this section, we assume that $A$ satisfies the condition~\ref{baseringA2}. We also assume throughout that $p>2$.
   
 \subsection{Definition of the essential submodule}

\begin{definition} Let $(\Pi,\rhob,t,d)$ be an admissible pseudo-deformation over $A$. Let $\rho: \Pi \rightarrow R^\ast$ be a $(t,d)$-representation, and define $G$, $\Gamma$, $L$ accordingly, with $L_2$ the derived Lie algebra of $L$. 
We call $S$ the set of elements $g \in G$ such that $\tr(g)=0$ and $-\det(g)$ is a square in $A^\ast$. We shall say that $(\Pi,\rhob,t,d)$ is {\it weakly odd} if the set $S$ is non empty.
\end{definition} 
\begin{definition} With the same notation as in the preceding definition, 
define
 $$A_\ess = \sum_{g \in S} W(\F) \tr(g L_2) \subset A.$$
 We call this $W(\F)$-submodule $A_\ess$ of $A$ the
 {\it essential submodule} of $A$ attached to $(\Pi,\rhob,t,d)$.
 \end{definition}
 Note that the condition of being weakly odd, and the $W(\F)$-submodule $A_\ess$ of $A$ depend only on $(\Pi,t,d)$, and not on the $(t,d)$-representation $\rho: \Pi \rightarrow R^\ast$, for  if $\rho': \Pi \rightarrow R'^\ast$ is another $(t,d)$-representation, then there exists an isomorphism $f: R \rightarrow R'$ preserving trace and determinant such that $\rho'=f \circ \rho$; the group $G' = \rho'(\Pi)$ is  the image $f(G)$,
and $\Gamma'=f(\Gamma)$, $L'=f(L)$, $L'_2=f(L_2)$. It is clear that $f$ realizes a bijection between $S$ and $S'$ and for every $g \in S$ a bijection between
the subgroups $\sum_{g \in S} gL_2$ of $R$ and $\sum_{g' \in S'} g' L'_2$ of $R'$. Since $f$ preserves traces, it follows that  $$A_\ess = \sum_{g \in S} W(\F) \tr(gL_2) = \sum_{g' \in S'} W(\F) \tr(g' L'_2).$$

The real motivation for introducing the submodule $A_\ess$ is its essential rôle in analyzing the density of modular forms modulo $p$, see section~\ref{sectiondensity} below. Meanwhile, $A_\ess$ can be considered as a very rough measure of how big the image $G$ if the pseudo-deformation is:
the bigger $G$, the bigger $\Gamma$, $L$ and $L_2$, and the more numerous the $g \in G$ such that $g^2=1$, hence the bigger $A_\ess$.
In this sense, most of the results below can be seen as big image theorems, though of a different type than the big image theorem of the previous section.

\begin{lemma} Let $(\Pi,t,d)$ be an admissible pseudo-deformation and $\rho: \Pi \rightarrow R^\ast$ a $(t,d)$-representation such that $J = \smat{1 & 0 \\ 0 & -1}\in G$ (and therefore $L=I_1 J \oplus \nabla$ is decomposable). Let $\Delta_2$ and $\nabla_2$ be the subgroups of diagonal and anti-diagonal matrices in $L_2$.
\begin{num} \label{decL20} One has $L_2=\Delta_2 \oplus \nabla_2,\ \Delta_2=[\nabla,\nabla]\text{ and }\nabla_2=[\Delta, \nabla] .$ \end{num}
One can write $\Delta_2 = I_2 J$ for some closed subgroup $I_2$ of $I_1$, and one has a decomposition:
\begin{num} \label{decL2} $L_2 = I_2 J \oplus \nabla_2.$\end{num}
\begin{num} \label{decL2inclusion} One has $I_2 \subset I_1$, $\Delta_2 \subset \Delta$, $\nabla_2 \subset \nabla$. \end{num}
 \begin{num}  \label{trJgammaL} For every $\gamma \in \Gamma$, one has $\tr(J \gamma L_1)=I_1$ and $\tr(J \gamma L_2)=I_2$. \end{num}
 \end{lemma} 
\begin{pf}
One has $L_2= [L,L] = [\Delta \oplus \nabla, \Delta \oplus \nabla] = [\nabla, \nabla] + [\Delta, \nabla]$ since $[\Delta,\Delta]=0$ (two diagonal matrices commute). 
But $[\nabla,\nabla]$ consists of diagonal matrices, and $[\Delta,\nabla]$ of antidiagonal ones. This proves~\ref{decL20} and~\ref{decL2}. Since $L_2 \subset L$, \ref{decL2inclusion} is clear

Let us prove \ref{trJgammaL}. By decomposition~\ref{decL2}, one has $\tr(J \gamma L_1) = \tr(\gamma) I_1 + \tr(J \gamma \nabla)$, and by Lemma~\ref{multtrgamma}, $\tr(\gamma) I_1 = I_1$. It therefore suffices to prove that $\tr (J \gamma \nabla) \subset I_1$. For this, let us denote by $\epsilon  \in \nabla$  the anti-diagonal part of $\gamma$ or of $\Theta(\gamma)$, and by $\eta$ any matrix in $\nabla$. One needs to show that $\tr (J \gamma \eta) = \tr (J \epsilon \eta) \in I_1$. Since $\epsilon$ and $\eta$ are anti-diagonal, one has $\tr(J \epsilon \eta) = - \tr(J \eta \epsilon)$, and thus $\tr (J \epsilon \eta) = \tr (J [\epsilon,\eta])/2$. Since  $[\nabla,\nabla]=\Delta_2 =  I_2 J$, one has $[\epsilon,\eta] \in I_2 J$ and one gets $\tr(J \epsilon \eta) \in \tr(J J I_2) = I_2 \subset I_1$, which completes the proof of \ref{trJgammaL} for $L_1$. The proof for $L_2$ is exactly the same. 
\end{pf}

\begin{lemma} If $(\Pi,t,d)$ is an admissible pseudo-deformation over $A$, and $f : A \rightarrow A'$ a surjective morphism of rings,
then $A'$ is a again a compact semi-local ring (for the quotient topology), and setting $t'=f \circ t$, $d'=f \circ d$, $(\Pi,t',d')$ is an admissible pseudo-deformation over $A'$. Moreover, if $A_\ess$ (resp. $A'_\ess$) is the essential submodule of $(\Pi,t,d)$ (resp. of $(\Pi,t',d')$), then $f(A_\ess)=A'_\ess$.
\end{lemma}
This is clear. 

In particular, if $A_{i,\ess}$ is the essential sub-module of $(\Pi,\rhob_i,t_i,d_i)$, then the projection $A \rightarrow A_i$ sends $A_\ess$ onto $A_{i,\ess}$. Note however that the map $A_\ess \rightarrow \prod_{i=1}^r A_{i, \ess}$ is not in general surjective.

Fix  an admissible pseudo-deformation $(\Pi,\rhob,t,d)$ and $\rho: G \rightarrow R^\ast$ a $(t,d)$-representation.
\begin{definition} Let $\bar S$ be the set of elements $\bar g \in \Gb$ such that $\tr(\bar g)=0$ and $-\det(g)$ is a square in $(A/\m)^\ast$.
\end{definition}
Note that the reduction map $G \rightarrow \Gb$ obviously induces a map $S \rightarrow \bar S$.
\begin{prop} \label{AessI2} The natural reduction map $S \rightarrow \bar S$ is surjective. For $g \in S$, the subgroup $\tr(g L_2)$ of $A$ only depends
 on the image $\bar g$ of $g$ in $\bar S$.
Moreover, for every $g \in S$, there exists a GMA $R'$, an isomorphism of $A$-algebra $f: R \rightarrow R'$ preserving traces and determinants, such that $f(g)=J$, and such that if $\rho'$ denotes the $R'$-valued $(t,d)$-representation $\rho'=f \circ \rho$, and  $L'_2$, $I'_2$ are defined using $\rho'$, then  one has $W(\F) \tr(g L_2) = W(\F) \tr(J L'_2)=W(\F) I'_2$.
\end{prop}
\begin{pf}
Let $\bar g \in \bar S$. By~\ref{dconstant}, there exists $\lambda \in \F^\ast$ such that $\det(\bar g) = - \lambda^2$ with $\lambda \in \F^\ast$. Denotes by $\bar g_i$ the image of the element $\bar g$ of $(R/\rad R)^\ast$ in $(R_i/\rad R_i)^\ast$. By definition of $\bar S$, there exists an element $g_0 \in \Pi$ such that $\rhob_i(g_0)=\bar g_i$ for $i=1,\dots,r$.

Since $\tr(\rhob_i(g_0))=0$, the eigenvalues of $\rhob_i(g_0)$ in $R_i/(\rad R_i)$ are $\pm \lambda$, two distinct elements of $\F^\ast$. Let us choose a $(t,d)$-representation $\rho'_i:  \Pi \rightarrow {R'_i}^\ast$ adapted to $(g_0,\lambda,-\lambda)$ (Prop.~\ref{exGMA}(iii)); let us set $R' = \prod_{i=1}^r R'_i$ and  $\rho' = \prod_{i=1}^r \rho'_i$, 
and let denote by $G'$, $\Gamma'$, $L'$, etc the group-theoretic and Lie theoretic data attached to $\rho'$.
Then $\rhob'(g_0) =\lambda J$ and by Theorem~\ref{diagonallift}, $\rho'(g_0)=s(\lambda) J \in G$.

Moreover, any lift $g' \in G'$ of $\rhob'(g_0)=\lambda J$ is of the form $s(\lambda) J \gamma$ with $\gamma \in \Gamma$, so by~\ref{trJgammaL}, $W(\F) \tr(g' L_2)= W(\F) \tr(J \gamma L_2)=W(\F) I'_2$, which is independent of $g'$.
There exists (Prop.~\ref{exGMA}(ii)) an isomorphism of $A$-algebras $f: R \rightarrow R'$ such that $f \circ \rho = \rho'$, preserving trace and determinant. By definition, if $g$ is a lift of $\bar g=\rhob(g_0)$ in $S$, then $g':= f(g)$ is a lift of $\rhob'(g_0)$ in $S'$, and $f(L_2)=L'_2$, so that $\tr(g L_2) = \tr(g' L'_2)=I'_2$, which is independent of $g$. 
\end{pf}

\begin{cor} One has $A_\ess = \sum_{\bar g \in \bar S} W(\F) \tr(gL_2)$ where in the summand $\tr(g L_2)$,
$g$ is an arbitrarily chosen lift of $\bar g$ in $S$. In particular, $A_\ess$ is a closed submodule of $A$.
\end{cor}
\begin{pf} The first assertion follows from the definition of $A_\ess$ and the proposition. Since the $\Z_p$-module $\tr(gL_2)$ is compact, so is $W(\F) \tr(g L_2)$. Since the set $\bar S$ is finite, it follows that $A_\ess$ is compact, hence closed in $A$. 
\end{pf}

\subsection{The key measure computation} \label{keymeasure}

In this subsection, we assume in addition to the preceding hypotheses that $A$ is an $\F$-algebra (equivalently, that $pA=0$). Therefore, in the results stated above,
each time there is a $W(\F) X$ where $X$ is an additive subgroup of $A$ or of $R$, it can just be replaced by $\F X$.

For any compact group $X$, we denote by $\mu_X$ the Haar measure on $X$ of total mass 1. We fix an admissible weakly odd pseudo-deformation $(\Pi,\rhob,t,d)$.

\begin{theorem} \label{compmeasure} Let $l : A \rightarrow \F$ be a linear form that is not identically $0$ on $A_\ess$.
Then 
\begin{num}\label{muP} $\mu_\Pi( (l \circ t)^{-1}(\F^\ast)) \geq  \frac{p-1}{pn},$ \end{num}
\noindent
 where  $n = | \Gb |$.
\end{theorem}

Since $l$ does not vanish on $A_\ess$, then by Prop.~\ref{AessI2}, for some $(t,d)$-representation $\rho$, $l$ does not vanish on $I_2$.
For the rest of this proof, we fix such a representation $\rho$ and the attached groups $G$, $\Gamma$, $L$, $L_2$, $I_2$.

Since $\rho$ is a surjective morphism of groups, the Haar measure $\mu_G$ is the direct image of  the measure $\mu_{\Pi}$ by $\rho$. Since $t = \tr_G \circ \rho$,  \ref{muP} is equivalent to:
\begin{num} \label{toprove0}
$ \mu_{G}  \left( (l\circ \tr_G)^{-1}(\F^\ast) \right) \geq \frac{p-1}{pn},$
\end{num}
\noindent which is the same thing as 
\begin{num} \label{toprove1}
$ \mu_{G}  \left( (l \circ \tr_G)^{-1}(0) \right) \leq \frac{1}{pn} + \frac{n-1}{n}$.
\end{num}
\noindent
To prove this, it is clearly enough to prove that
\begin{num}  \label{toprove2}
$ \mu_{G} \left( (l \circ \tr_G)^{-1}(0) \cap J \Gamma \right) \leq \frac{1}{pn}$,
\end{num}
\noindent
since $\mu_G(G-J \Gamma) = \frac{n-1}{n}$, $G-J\Gamma$ being the union of $n-1$ 
$\Gamma$-cosets each of measure $1/n$.
Let $m_J$ be the injective map $\Gamma \rightarrow G$, $\gamma \mapsto J \gamma$, whose image is the coset $J \Gamma$, and let $\mu_\Gamma$ be the Haar measure of total measure $1$ on $\Gamma$.
Clearly, \ref{toprove2} is equivalent to
\begin{num} \label{toprove3}
$ \mu_{\Gamma} \left( (l \circ \tr_G \circ m_J)^{-1}(0) \right) \leq \frac{1}{p}$
\end{num}
\noindent
Now consider the exact sequence $1 \rightarrow \Gamma_2 \rightarrow \Gamma \rightarrow \Gamma/\Gamma_2 \rightarrow 1$.
By Fubini's theorem, to prove~\ref{toprove3} it is enough to prove that for all $\gamma \in \Gamma/\Gamma_2$,
\begin{num}  \label{toprove4}
$\mu_{\Gamma_2} \left( (l \circ \tr_G \circ m_{J\gamma})^{-1}(0) \right) \leq \frac{1}{p}$
\end{num} \noindent where $m_{J \gamma}$ is 
the map $\Gamma_2 \rightarrow G$, $\gamma_2 \mapsto J \gamma \gamma_2$ and $\mu_{\Gamma_2}$ the Haar measure on $\Gamma_2$ of total measure 1.
Since $\Theta^{-1} : L_2 \rightarrow \Gamma_2$ is a measure-preserving homeomorphism (Prop.~\ref{ThetapreservesHaar}), it suffices to prove
\begin{num}  \label{toprove5}
$\mu_{L_2} \left( (l \circ \tr_G \circ m_{J\gamma} \circ \Theta^{-1})^{-1}(0) \right) \leq \frac{1}{p}$
\end{num}
\noindent
To simplify notation let us define the map $$h_\gamma = \tr_G \circ m_{J \gamma} \circ \Theta^{-1}  : L_2 \stackrel{\Theta^{-1}}{\longrightarrow} \Gamma_2 \stackrel{m_{J\gamma}}{\longrightarrow} R^\ast \stackrel{\tr}{\longrightarrow}A,$$
so that \ref{toprove5} becomes
\begin{num}  \label{toprove6}
$\mu_{L_2} \left( (l \circ h_\gamma)^{-1}(0) \right) \leq \frac{1}{p}$
\end{num}
\noindent To prove~\ref{toprove6}, we shall use the following result:
\begin{prop}  
\begin{itemize} Fix $\gamma \in \Gamma$.
\item[(i)] There exists a measure preserving homeomorphism $\Psi = \Psi_\gamma: L_2 \rightarrow L_2$ such that $h_\gamma \circ \Psi^{-1}: L_2 \rightarrow A$ is $\F_p$-affine.
\item[(ii)] The image of $h_\gamma$ is the $\F_p$-affine subspace $\tr(J\gamma) + I_2$ of $A$.
\end{itemize}
\end{prop}
\begin{pf}

Let us define a map $\Psi : L_2 \rightarrow L_2$ by setting
$$\Psi(m) = m + \sigma(m)\text{ with }\sigma(m)= (\sqrt{1+\tr(m^2)/2}-1) \frac{\tr(J \gamma)}{\tr( \gamma )}  J .$$
Let us check that $\Psi$ is well-defined. Write $\gamma=\smat{a & b \\ c & d}$.
First, one has $\tr(\gamma) = a + d \equiv 2 \pmod{\m}$, hence $\tr(\gamma)$ is invertible in $A$,
and the formula defining $\Psi(m)$ makes sense as an element of $R$. We need to check that it is indeed in $L_2$.  
By definition, $\Theta(\gamma) = \smat{ (a -  d)/2 & b \\ c & (d-a)/2 }$ is in $L$, and since $L$ is decomposable,  $\smat{ (a -  d)/2 & 0\\ 0 & (d-a)/2 }$ is in $L$; on the other hand and one computes
$\tr(J \gamma) J = \smat{ a -  d & 0\\ 0 & d-a }$, so  $\tr(J \gamma) J \in L$.
One has $\tr(\gamma)^{-1} L = L$ by Lemma~\ref{multtrgamma},
hence $\frac{\tr(J \gamma)}{\tr( \gamma )}  J$ is in $L$. On the other hand  $\sqrt{1+\tr(m^2)/2}-1=\sum_{n=1}^\infty {n \choose 1/2}\frac{tr(m^2)^n}{2^n}$ is in $P(\Gamma)$ so sends $L$ into $L_3\subset L_2$. Hence $\sigma(m)$ is in $L_2$ and $\Psi$ is well-defined. 

If $m,m'$ are in $L_2 \subset I^2 R$, and $m-m' \in \m^n R$ then one sees that $$\sqrt{1+\tr(m^2)/2}-\sqrt{1+\tr(m'^2)/2} = \sum_{n=1}^\infty {n \choose 1/2}\frac{\tr(m^2)^n-\tr(m'^2)^n}{2^n} \in \m^{n+2},$$ hence $\sigma(m)-\sigma(m') \in \m^{n+2} L \subset \m^{n+3} R$. Therefore, by Lemma~\ref{lemmahomeo}, $\Psi: L_2 \rightarrow L_2$ is a measure-preserving homeomorphism.

For $m \in L_2$, one has
\begin{eqnarray*} h_\gamma(m) &=& \tr ( J \gamma \Theta^{-1}(m) ) \\ &=& \tr(J \gamma m) +  \tr (J \gamma) \sqrt{1+\tr(m^2)/2} \\
&=& \tr(J \gamma) + \tr (J \gamma \Psi(m)).\end{eqnarray*}
Therefore $h_\gamma(\Psi^{-1}(m)) =  \tr(J \gamma) + \tr (J \gamma m)$,
which shows that $h_\gamma \circ \Psi^{-1}$ is an affine map as stated in (i), whose image is the affine space $\tr(J \gamma) + \tr(J \gamma L_2)=\tr(J\gamma)+ I_2$ by~\ref{trJgammaL}.
 
\end{pf}

Using the proposition and the map $\Theta$ it introduces, we see that to prove ~\ref{toprove6}, it is enough to prove
that
\begin{num}  \label{toprove7}
$\mu_{L_2} \left( (l \circ h_\gamma \circ \Theta^{-1})^{-1}(0) \right) \leq \frac{1}{p}$
\end{num}
But $h_\gamma \circ \Theta^{-1}$ is a $\F_p$-affine map. So $l \circ h_\gamma \circ \Theta^{-1}$ is 
an $\F_p$-affine map on $L_2$ with values in $\F$, and with image the $\F_p$-affine subspace $l(\tr(J\gamma)) + l(I_2)$.
 Since $l(I_2) \neq 0$, the image $S$ of our map $l \circ h_\gamma \circ \Theta^{-1}$ is an affine $\F_p$-subspace of positive dimension of $\F$.  The measure $\mu_{L_2} \left( (l \circ h_\gamma \circ \Theta^{-1})^{-1}(0) \right)$ is 0 if $S$ does not contain 0, and $1/|S|$ otherwise. In any case, it is less than $\frac{1}{p}$ which proves \ref{toprove7}
 and the theorem.

\begin{remark} \label{LL2} If we assume that $\Theta(\Gamma)=L$, then we can prove that the inequality  $\mu_\Pi( (l \circ t)^{-1}(\F^\ast)) \geq  \frac{p-1}{pn},$ holds only when $l(I_2) \neq 0$, but more generally when $l(I_1) \neq 0$.
Indeed, to prove \ref{toprove3} for such an $l$, that is that $ \mu_{\Gamma} \left( (l \circ \tr_G \circ m_J)^{-1}(0) \right) \leq \frac{1}{p}$, it is enough to prove that $ \mu_{L} \left( (l \circ \tr_G \circ m_J \circ \Theta^{-1})^{-1}(0) \right) \leq \frac{1}{p}$.
But the map $\tr_G \circ m_J \circ \Theta^{-1}$ is very simple: it sends a matrix $m=\smat{a & b \\ c &-a}$ to $\tr(J\Theta^{-1} m)=2a$.
In particular, this map is linear, and its image is the group $I_1$. Thus if $l$ is non-zero on $I_1$, the map   $(l \circ \tr_G \circ m_J \circ \Theta^{-1})^{-1}(0) $ is a $\F_p$-affine map from $L$ to $\F$ whose image has positive dimension, and we conclude easily.
\end{remark}
\subsection{A sufficient condition for the largeness of $A_\ess$}
 
 In this subsection (and for the rest of this section) we assume that $A$ is local.
 
\begin{definition} An admissible pseudo-deformation $(\Pi,\rhob,t,d)$ is said to be {\it virtually abelian} if there exists an open subgroup $\Pi_0$ of $\Pi$ such that the restriction $(t_{|\Pi_0},d_{|\Pi_0})$ factors trough an abelian quotient of $\Pi_0$.
\end{definition}

\begin{lemma} \label{Aessnotzero} Let $(\Pi,\rhob,t,d)$ be a weakly odd admissible pseudo-deformation. Assume that $A$ is a domain. If $A_\ess = 0$, then $(\Pi,\rhob,t,d)$ is virtually abelian.
\end{lemma}
 \begin{pf} Let us pick $g_0 \in \bar S$  and choose $\rho: \Pi \rightarrow R^\ast$  a $(t,d)$-representation adapted to $g_0$. Thus $L$ is decomposable and $\tr(g L_2) = I_2 \subset A_\ess$ so by hypothesis $I_2=0$. By~\ref{trnabla2I1},
 $I_1 \tr(\nabla^2) \subset I_2 = 0$. Since $A$  is a domain, it follows that either $I_1=0$ or $\tr(\nabla^2)=0$. We claim that in both case,
 $(\Pi,\rhob,t,d)$ is virtually abelian.
 
 If $I_1=0$, then $L=\nabla$, and $L_2 = [\nabla,\nabla] =\Delta_2 \subset I_1 J = 0$. Thus $L$ is commutative. It follows that $\Theta^{-1}(L)$ is commutative and $\Gamma$ is commutative. Let $\Pi_0 = \ker \rhob$. Then $\rho(\Pi_0) = \Gamma$ is commutative, which proves that $(t_{|\Pi_0},d_{|\Pi_0})$ factors trough an abelian quotient of $\Pi_0$.
 
  Remember (\cite[\S1.3]{BC}) that since $A$ is a domain, we may assume that if $R = \mat{A & B \\ C & A}$, for $b \in B$ and $c \in C$, $bc = 0 \Rightarrow b=0 \text{ or }c=0$. If $\tr(\nabla^2)=0$, and if $\epsilon = \smat{0 & b \\ c & 0} \in \nabla$, then $\tr(\epsilon^2)=2bc=0$.
  If $\epsilon \neq 0$, by symmetry we may assume that $c = 0$, $b \neq 0$. Then for $\epsilon'=\mat{0 & b' \\ c' &0} \in \nabla$, we  
 have $bc'=\tr(\epsilon \epsilon')=0$, so $c'=0$. Thus all matrices in $\nabla$ are upper-triangular. It follows that $L$ itself, hence $\Gamma \subset \Theta^{-1}(L)$
is contained in the set of triangular matrices. If again we set $\Pi_0 = \ker \rhob$, we see that $t_{| \Pi_0}=\tr \rho_{|\Pi_0}$ is the sum of two characters, hence factors through an abelian quotient.
\end{pf} 

\begin{prop} \label{propAesslarge} Let $(\Pi,\rhob,t,d)$ be a weakly odd admissible pseudo-deformation, which is not virtually abelian. Assume that $A$ is a domain. Let $g_0 \in \bar S$, and $\rho: \Pi \rightarrow R^\ast$  a $(t,d)$-representation adapted to $g_0$. Assume that for this representation, either $W(\F) B_1$ or $W(\F) C_1$ is not a finite-type $W(\F)$-module. Then $A_\ess$ is not a finite type $W(\F)$-module either.
\end{prop}
\begin{pf} 
Assume by contradiction that $A_\ess$ is a finite type $W(\F)$-module. Since $W(\F) I_2 \subset A_\ess$, so is $W(\F) I_2$. By the preceding lemma, $I_2 \neq 0$, and since $P I_2 \subset I_2$ and $A$ is a domain, it follows that $W(\F) P$ is a finite type $W(\F)$-module. Therefore $W(\F) P + W(\F) I_2$ is a finite type $W(\F)$-module. But if $\smat{0 & b \\ c & 0}$ and $\smat{0 & b' \\ c' & 0}$ are any elements in $\nabla$, then $bc'-b'c \in I_2$ and $bc'+b'c \in P$,
so $bc'$ is in $I_2 + P$, and $I_2 + P$ contains $B_1 C_1$, so $W(\F) I_2 + W(\F) P$ contains $W(\F) B_1 C_1$. But $W(\F) B_1$ and $W(\F) C_1$ are non-zero (otherwise the deformation would be virtually abelian) and by assumption one of them is not a finite $W(\F)$-module,
so it follows that $W(\F) B_1 C_1$ is not a finite $W(\F)$-module, a contradiction.
\end{pf}

 \subsection{The essential subgroup in the reducible case}

In this subsection we keep assuming that $A$ is local and we fix an admissible weakly odd pseudo-deformation $(\Pi,\rhob,t,d)$, and we assume throughout that $\rhob$ is reducible.
\begin{num} There exists two continuous characters $\chi_1, \chi_2 : \Pi \rightarrow \F^\ast$, such that $
\rhob \simeq \chi_1 \oplus \chi_2$.
\end{num}
Let us chose  a $(t,d)$-representation which is well-adapted in the sense of Definition~\ref{welladapted}. Thus the group $\bar G$ is a diagonal subgroup of $\GL_2(\F)$, and $s(\bar G) \subset G$. Since $\rhob$ is weakly odd there exists in $\bar G$ an element of order $2$ other than $\pm 1$, and since this element is diagonal, it is either $J$ or $-J$. There is no loss of generality in supposing that $J \in \bar G$, hence $J=s(J) \in G$.

\begin{prop} \label{Aessredcase} One has $A_\ess = W(\F) I_2$.
\end{prop}
\begin{pf} Indeed, $A_\ess = \sum_g W(\F) \tr(\bar g L_2)$ where the sum runs over all 
$\bar g \in \Gb$ except $\pm 1$.
But the only such $\bar g$ in the diagonal subgroup $\Gb$ are $J$ and possibly $-J$. In any case, $s(\bar g)=\pm J$ 
is a lift of $\bar g$ in $G$ whose square is $1$, so $\tr(\bar g L_2)=\tr(\pm J IL_2)=I_2$.
\end{pf}

\begin{prop} \label{Aesslargered} Assume that $A$ is a domain and is not a finite $W(\F)$-module. Then if $(\Pi,\rhob,t,d)$ is not virtually abelian, $A_\ess$ is not a finite $W(\F)$-module. 
\end{prop}
\begin{pf}  Since $A= W(\F) \oplus W(\F) I_1 + W(\F) P$ (Prop.~\ref{Amreducible}), either $W(\F) I_1$ or $W(\F) P$ is not finite as a $W(\F)$-module. If $W(\F) I_1$ is not finite, then neither is $W(\F) I_1^2$ since $A$ is a domain, and since $I_1^2 \subset P$, neither is $W(\F) P$. So in any case $W(\F) P$ is not finite as a $W(\F)$-module.

Under our hypotheses $A_\ess = W(\F) I_2$ is not zero by Prop.~\ref{Aessnotzero}. Since $P I_2 \subset I_2$, and $A$ is a domain, $W(\F) I_2$ is not  finite as a $W(\F)$-module.
\end{pf}

\begin{theorem} \label{thmAessred} Assume that the character $\chi_1/\chi_2$ is not of order $2$, or in other words that the projective image of $\rhob$ is not $\Z/2\Z$. Then $A_\ess$ is an ideal of $A$, and more precisely it is the {\it reducibility ideal} of the pseudo-representation $(t,d)$ (see \cite[\S1.5]{BC}).
\end{theorem}
\begin{pf} By theorem~\ref{structureLdihedral}, one has $W(\F) L = \mato{\tilde I_1 & B \\ C & \tilde I_1}$ for some $W(\F)$-module $\tilde I_1$. 
 It follows that $W(\F) I_2=BC$.
\end{pf}
 
 \subsection{The essential subgroup in the dihedral case}
 
 In this subsection we still assume that $A$ is local and we fix an admissible odd pseudo-deformation $(\Pi,\rhob,t,d)$, and we assume throughout that 
  \begin{num} \label{Aessdih} The projective image of $\rhob$ is dihedral. \end{num}
  
As in~\S\ref{consequencesdihedral}, we choose a well-adapted $(t,d)$-representation $\rho: G \rightarrow \GL_2(A)$ which encompasses
the choice of a subgroup $D$ of index $2$ in $\bar G$ consisting of diagonal matrices.

\subsubsection{Largeness of $A_\ess$}

\begin{prop} \label{Aesslargedih} Assume that $A$ is a domain and is not a finite $W(\F)$-module, that $(\Pi,\rhob,t,d)$ is not virtually abelian, and~\ref{Aessdih}.
Then $A_\ess$ is not a finite $W(\F)$-module.
\end{prop}
\begin{pf} We first claim that $W(\F) B_1$ is not a finite $W(\F)$-module. Indeed, it is non-zero otherwise  $(\Pi,\rhob,t,d)$ would be virtually abelian. Moreover,  $W(\F) I_1 + W(\F) P + W(\F) B_1$ is not a finite $W(\F)$-module. Therefore at least one of the three terms is not a finite $W(\F)$-module. If it is the third, then we are done, and if it is one of the two first, we are also done since $I_1 B_1 \subset B_1$ and $P B_1 \subset B_1$. 

The proposition then follows from Prop.~\ref{propAesslarge}
\end{pf}

\subsubsection{Description of $A_\ess$ in the case $4 \mid n$, $n>4$}

Let $n$ be the order of the projective image of $\rhob$. Since $\rhob$ is dihedral, $n \geq 4$ and $n$ is even.
\begin{num}  
\label{Aessn4}
We assume that $n>4$, and that $4 \mid n$. 
\end{num}
Under this assumption, the image of the diagonal group $D$ in $\PGL_2(\F)$ has even order, and thus contains an element of order $2$.
Fix a lift $\bar g$ of that element in $D$. This element $\bar g$ has trace zero, hence is of the form $\lambda J$ for some $\lambda \in \F^\ast$, and is an element of $\bar S$.
An element of $\Gb - D$ also has trace 0. We shall make the following supplementary assumption:
\begin{num}
\label{AessGb} There exists an element $\bar g'$ of $\Gb-D$ such that $-\det \bar g'$ is a square in $\F^\ast$, or in other words, such that $\bar g' \in \bar S$. \end{num}
This assumption will be harmless in the applications (see \S\ref{sectiondensity} below), since if not true, we can always choose an element $\bar g'$ in $\Gb -D$ and extend the scalars from $\F$ to the quadratic extension $\F'$ of $\F$ generated by $\sqrt{-\det \bar g'}$. 

\begin{theorem} \label{thmAessdih} Assume~\ref{Aessdih}, \ref{Aessn4} and \ref{AessGb}. 
Then $A_\ess = \m B$. In particular $A_\ess$ is an ideal of $A$.
\end{theorem} 
\begin{pf} Since $n>4$, $W(\F) L =\mato{ W(\F) I_1 & W(\F) B_1 \\ W(\F) B_1 & W(\F)}$ by Theorem~\ref{structureLdihedral}. It follows that
$W(\F) L_2 =  \mato{ W(\F) B_1^2 & W(\F) I_1 B_1 \\ W(\F) I_1 B_1 & W(\F) B_1^2 }$

We claim that
$$A_\ess = W(\F) B_1^2 + W(\F) I_1 B_1.$$
Indeed, $W(\F) \tr(s(\bar g) L_2) = W(\F) I_2 = W(\F) B_1^2 \subset A_\ess$, and $W(\F) \tr(s(\bar g') L_2) = W(\F) I_1 B_1 \subset A_\ess$,
and if there are other elements $\bar g''$ in $S$, they are either diagonal or anti-diagonal, contributing the same summand $W(\F) B_1^2$ or
$W(\F) I_1 B_1$.

To prove that $A_\ess$ is an ideal, we recall that $A=W(\F) + W(\F) I_1 + W(\F) I_1^2 + W(\F) B_1$, so we only need to check that $A_\ess$ is stable by multiplication by $I_1$ and $B_1$. We have $I_1 A_\ess = W(\F) I_1 B_1^2 + W(\F) I_1^2 B_1$, and since $I_1 B_1 \subset B_1$, we see  that $I_1 A_\ess \subset A_\ess$. We have $B_1 A_\ess = W(\F) B_1^3 +  W(\F) I_1 B_1^2$, and since $W(\F) B_1^2 \subset W(\F) I_1$ and $I_1 B_1 \subset B_1$, we see that $B_1 A_\ess \subset A_\ess$.

Since $B_1 \subset B$ and $B$ is an $A$-ideal, it is clear that $A_\ess \subset B$.
We claim that the ideal (of $A$) generated by $B_1$ is $B$.  

Since $A_\ess$ is an ideal, we get $A_\ess = B^2 + I_1 B$. Since $\m = W(\F) I_1 + W(\F) I_1^2 + B$, we have $\m B = I_1 B + I_1^2 B + B^2 = A_\ess + I_1^2 B$. But since $B$ is an ideal, $I_1 B \subset B$ and $I_1^2 B \subset I_1 B \subset A_\ess$, so $\m B=A_\ess$.
\end{pf}

 \subsection{The essential subgroup in the large image or exceptional case}

We assume that $A$ is local, and we assume that $\rhob$ has large and exceptional projective image. In this case, things are pretty simple:
\begin{theorem} \label{Aesslargeirr} If $\rhob$ has large or exceptional projective image, then $A_\ess = \m^2$.
\end{theorem}
\begin{pf} By Theorem~\ref{structureLlargeimage}, one has for a suitable $(t,d)$-representation $\rho$, $W(\F) L = \mato{\m & \m \\ \m & \m}$. Since $\mato{\m & \m \\ \m & \m}$ is invariant by conjugation by any trace-preserving automoprhism of $R$, it follows  that $W(\F) L = \mato{\m & \m \\ \m & \m}$ for {\it any} $(t,d)$-representation $\rho$. For any $g \in S$ we therefore have $W(\F) \tr (g L_2) = W(\F) I_2 = \m^2$, and $A_\ess = \m^2$.
\end{pf}

\section{An example}
\label{sectionexample}
 The aim of this section is to provide an example of an admissible pseudo-representation whose image is `complicated', and which violates the conclusions (and of course, the hypotheses) of certain theorem we have proved earlier. It can be safely skipped.
 
Let $\F$ be a finite field of characteristic $p>2$.
Let $A = \F[[X]]$, with maximal ideal $\m = X \F[X]]$.

\subsection{A two-generator closed subgroup $\Gamma$ of $SL_2^1(A)$ and its Lie algebra}

Define $$g = \smat{X + \sqrt{1+X^2} & 0 \\ 0 & -X + \sqrt{1+X^2}} \text{ and }h = \smat{ \sqrt{1-X^2} & X \\ -X & \sqrt{1- X^2}}.$$ Note that those two matrices belongs to $\SL_2^1(A)$. Let $\Gamma$ be the topological closure of the subgroup of $SL_2^1(A)$ generated by $g$ and $h$. 

\begin{lemma} \label{lemmaGammaJ} With $J=\smat{1 & 0 \\ 0 &-1}$, one has $J g J = g$ and $J h J = h^{-1}$. One has $J \Gamma J = \Gamma$. With $J'=\smat{0 & 1 \\ 1 & 0}$, one has $J' g J' = g^{-1}$ and $J' h J'= h^{-1}$. One has $J'\Gamma J'=\Gamma$.
\end{lemma}
\begin{pf} The first and third sentences consist of two trivial computations each and the second and fourth sentences follow. \end{pf}

\begin{lemma} \label{lemmaGammaabcd} Suppose that $\gamma = \smat{a(X) & b(X) \\ c(X) & d(X)}$ is in $\Gamma$. Then $a(X)=d(-X)$ and $b(X) = c(-X)$.
\end{lemma} 
\begin{pf} The equalities $a(X)=d(-X)$ and $b(X)=c(-X)$ are clearly true for the matrices $g$ and $h$, and also $g^{-1}$ and $h^{-1}$. If these equalities are true for $\gamma = \smat{a(X) & b(X) \\ c(X) & d(X)}$ and $\gamma' = \smat{a'(X) & b'(X) \\ c'(X) & d'(X)}$, then $\gamma \gamma' = \smat{a a' + b c' & ab'+bd' \\ a'c+c'd & d d' + b' c}$ and one see
that $(aa'+bc')(X)=a(X) a'(X) + b(X) c'(X) = d(-X) d'(-X) + c(-X)b'(-X) = (dd'+b'c) (-X)$, and $(ab'+bc')(X) = a(X) b'(X) + b(X) d'(X) = d(-X) c'(-X) + c(-X) a'(-X) =  (a'c+c'd) (-X)$.  Therefore they are true for any element of the subgroup generated by $g$ and $h$, and of its closure, hence the lemma.
\end{pf}

Define a subspace $L$ of $R$ as follows:
$$L = \{ \smat{a & b \\ c & -a}, a,b,c \in \m = X \F[[X]],\  a(X)=-a(-X),\  b(X)=c(-X)\}.$$
In other words, $L = X \F[[X^2]] J \oplus \nabla$, with $\nabla = \{ \smat{0 & b \\ c & 0 }, b,c \in \m = X \F[[X]],\   b(X)=c(-X)\}$. In particular, $L$ is decomposable, but not strictly decomposable.

\begin{lemma} The Pink's Lie algebra $L(\Gamma)$ of $\Gamma$ is $L$.
\end{lemma}
\begin{pf} First we prove that $L(\Gamma) \subset L$. It suffices to prove that $\Theta(\gamma) \in L$ for every $\gamma \in \Gamma$.
If $\gamma=\smat{a & b \\ c & d}$, then by Lemma~\ref{lemmaGammaabcd}, $a(X)=d(-X)$ and $b(X) =  c(-X)$, and $\Theta(\gamma)=\smat {(a(X) - a(-X))/2 & b(X) \\ c(X) & (a(X) - a(-X))/2}$ which is clearly in $L$.

Next, observe that by Lemma~\ref{lemmaGammaJ}, $L(\Gamma)$ is decomposable. We write $L(\Gamma)=I_1 J \oplus \nabla_1$, with $\nabla_1$ anti-diagonal. Also, $\Theta(g) = \smat{X/2 & 0 \\ 0 & -X/2}$ belongs to $L(\Gamma)$, so $X \in I_1$ and $2 \tr(\Theta(g)^2) = X^2$ belongs to the closed sub-pseudoring $P(\Gamma)$ of $A$. It follows that $X^2 \F[[X^2]] \subset P(\Gamma)$.  Since $I_1$ is stable by $P(\Gamma)$, we get $I_1=X \F[[X^2]]$.
From $\Theta(h) \in L(\Gamma)$ and $L(\Gamma)$ decomposable, we get   $\smat{0 & X \\ -X & 0} \in \nabla_1$. Since $\nabla_1$ is stable by 
taking the Lie bracket with $X J \in L(\Gamma)$, we see that $\smat{0 & X^2 \\ X^2 & 0}$, $\smat{0 & X^3 \\ -X^3 & 0}$, etc. belong to $\nabla_1$, and finally $\nabla_1 =  \{ \smat{0 & b \\ c & 0 }, b,c \in \m = X \F[[X]],\   b(X)=-c(-X)\}$. Hence $L(\Gamma)=L$.
 \end{pf}

\subsection{Construction of two admissible pseudo-deformations}

We define $$G = \Gamma \coprod J \Gamma.$$ It follows from the first part of Lemma~\ref{lemmaGammaJ} that $G$ is a closed subgroup of $\GL_2(A)$,
containing $\Gamma$ as a subgroup of order $2$, and that $G$ is the semi-direct product of  $\{1,J\}$ by $\Gamma$.

Let $\Pi$ be any pro-finite group with a continuous surjective morphism onto $G$ (for example $\Pi=G$ with the identity). Let $\rho$ be the composition $\Pi \rightarrow G \rightarrow \GL_2(A)$. Let $t = \tr \rho$, $d = \det \rho$. Let $\rhob: \Pi \rightarrow \GL_2(\F)$ be the reduction modulo $\m$ of $\rho$. Then $\rhob$ is a continuous semi-simple representation of $\Pi$ with image (and projective image) isomorphic to $\Z/2\Z$.

\begin{lemma} $(\Pi,\rhob,t,d)$ is an admissible pseudo-deformation. The projective image of $\rhob$ is cyclic of order $2$. \end{lemma}
\begin{pf} The representation $\rhob$ is the sum of the trivial character and a character of order $2$ of $\Pi$, so $\rhob$ 
satisfies~\ref{tmultfree}. The property~\ref{deform} is obvious. One has $d(\Gamma)=1$, $d(G-\Gamma)=d(J\Gamma)=-1$, which makes clear that~\ref{dconstant} holds. For~\ref{tgenerate}, one has $\tr(Jg) = 2 X$, hence the smallest closed subring of $A$ containing $\tr(G)$ contains $\F[X]$, hence is $A$.
\end{pf}

Let $H$ be the subgroup of order $8$ of $\GL_2(A)$ generated by $J$ and $J'$. By Lemma~\ref{lemmaGammaJ}, $H$ normalizes $\Gamma'$. We define $G'=\Gamma H$, a semi-direct product of $H$ by $\Gamma$. Let $\Pi'$ be any pro-finite group with a continuous surjective morphism onto $G'$ (for example $\Pi'=G'$ with the identity). Let $\rho'$ be the composition $\Pi' \rightarrow G' \rightarrow \GL_2(A)$. Let $t' = \tr \rho'$, $d' = \det \rho'$. Let $\rhob': \Pi' \rightarrow \GL_2(\F)$ be the reduction modulo $\m$ of $\rho'$. Then $\rhob'$ is a continuous semi-simple representation of $\Pi$ with image  isomorphic to $H$.

\begin{lemma} $(\Pi',\rhob',t',d')$ is an admissible pseudo-deformation. The projective image of $\rhob'$ is dihedral of order $4$. \end{lemma}
\begin{pf}
The proof if the same as above, except for the projective image, which is the image of $H$ in $\PGL_2(\F)$. This image is generated by the image of $J$ and $J'$, elements of order $2$ that commute in $\PGL_2(\F)$ since in $\GL_2(\F)$ one has $J' J J' = -J$.
\end{pf}

\subsection{Counter-examples to over-optimistic statements}

We now use the admissible pseudo-deformations $(\Pi,\rhob,t,d)$ and  $(\Pi',\rhob',t',d')$ to construct counter-examples.

First, we show that Theorem~\ref{thmconglarge} is false if we do not assume that $\rhob$ is regular. More precisely, we show that it does not assume in one case where $\rhob$ has projective image cyclic of order $2$, and in one case where its has projective image dihedral of order $4$.
\begin{prop} Let $(\Pi,\rhob,t,d)$ be the admissible pseudo-deformation constructed in the above subsection. 
There is no subgroup $\Pi_0$ of $\Pi$ containing $\ker \rhob$, and subring $A_0$ of $A$ such that
the pseudo-representation $(t,d)$ of $\Pi_0$ takes value in $A_0$, is admissible, and $(t_{|\Pi_0},d_{\Pi_0})$ has congruence-large image.
The same holds with $(\Pi,\rhob,t,d)$ replaced by $(\Pi',\rhob',t',d')$ constructed in the above subsection.
 \end{prop}
 \begin{pf} If $\Pi_0$ is a subgroup as in the statement, then either $\Pi_0=\Pi$ or $\Pi_0 = \ker \rhob$ has index $2$ in $\Pi$.  The second case is excluded since $\rhob_{|\ker \rhob}$ is the trivial representation of dimension 2, which is not multiplicity free. Thus $\Pi_0=\Pi$ and $A_0=A$. 
 We just have to show that for the unique $(t,d)$-representation $\rho$, $\rho(\Pi_0)=G$ does not contain any congruence subgroup.
 But if it did, $\Gamma$ would contain a congruence subgroup and $L$ would contain a sub-module of the form $\mat{I & I \\ I & I}^0$ for some non-zero ideal $I$ of $A$. Since up-left coefficients of $L$ are odd elements of $\F[[X]]$, $I$ would contain only odd functions, but this is absurd since $I$ is stable by multiplication by $X$. 
 
 The same result for $(\Pi',\rhob',t',d')$ is proved similarly.
 \end{pf}

 Second, we show that it may be false, when $n=2$, that $A_\ess$ is an ideal of $A$.
 \begin{prop} The $\Z_p$-submodule $A_\ess$ of $A$ attached to the admissible pseudo-deformation $(\Pi_0,\rhob,t,d)$
 is not an ideal of $A$.
 \end{prop} 
 \begin{pf} By Prop.~\ref{Aessredcase} we have $A_\ess=I_2\subset I_1$. Since $I_1$ consists of odd element of $A=\F_p[[X]]$, so does $I_2$,
 but no non-zero ideal of $A$ consists only of odd elements.
 \end{pf}
 
\subsection{The group $G$ as a Galois group}

Lest the reader think that the pathological example $(\Pi,\rhob,t,d)$ is allowed only by our too lenient definition of an admissible representation, and does not happen in the concrete   applications to number theory, we show that when $p=3$ (to fix ideas) one can take in the above example for $\Pi$ the absolute Galois group $G_{\Q,3}$ and for $(t,d)$ the quotient by a prime ideal of height one of the canonical pseudo-representation of $G_{\Q,3}$ over the Hecke algebra of modular forms modulo $3$.

Let $G_{\Q(\mu_3),3}$ be the Galois group of the maximal algebraic extension of $\Q(\mu_3)=\Q(\sqrt{-3})$ unramified above the unique place above$3$. This is a subgroup of order $2$ of $G_{\Q,3}$, and $G_{\Q,3}$ is a semi-direct product of $\{1,c\}$, where $c$ is any complex conjugation, by $G_{\Q(\mu_3),3}$.

Let $G_{\Q(\mu_3),3}^3$ be the largest quotient of $G_{\Q(\mu_3),3}$ which is a pro-$3$-group. The structure of that group is known. Let $c$ be a complex conjugation in $G_{\Q,3}$.
\begin{lemma} There exists an element $g \in G_{\Q(\mu_3),3}^3$ such that $G_{\Q(\mu_3),3}^3$ is a free pro-3-group with $g$ and $cgc$ as pro-generators. \end{lemma}
The freeness of $G_{\Q,3}^3$ is due to Shafarevich, see \cite[page 82, example after theorem 5]{sha}. The rest of the lemma is proven in \cite{annap3}. 

Consider the unique continuous morphism of groups $f: G_{\Q(\mu_3),3}^3 \rightarrow \Gamma$ sending $g$ to $(xy)^{1/2}$ and $cgc$ to 
$(xy^{-1})^{1/2}$ (the square root $z^{1/2}$ for $z$ an element of the pro $3$-group $\Gamma$ is defined as usual as the limit $z^{a_n}$ where $a_n$ is a sequence of natural integers converging $3$-adically to $1/2$). Since the group generated by $(xy)^{1/2}$ and $(xy^{-1})^{1/2}$ contains $x$ and $y$, $f$ is surjective. Using the structural surjective map $ G_{\Q(\mu_3),3} \rightarrow  G_{\Q(\mu_3),3}^3$, we see $f$ as a surjective morphism $G_{\Q(\mu_3),3} \rightarrow \Gamma$.
Since $f(cgc) = J f(g) J^{-1}$ in $G$ by Lemma~\ref{lemmaGammaJ}, and $J^2=c^2=1$, we can extend $f$ into a surjective morphism $f: G_{\Q,3} \rightarrow G$ sending $c$ onto $J$.
We thus get a pseudo-character $(t=\tr \circ f,d = \det \circ f)$ on the Galois group $\Pi=G_{\Q,3}$ which is an admissible pseudo-deformation of $\rhob=1 \oplus \omega_3$ and whose image is $G$. As seen above, this Galois pseudo-deformation is a counter-example to the assertion that $A_\ess$ is an ideal and that $(t,d)$ has congruence-large image.

Finally, note that if $R_\rhob$ denotes the universal deformation of $\rhob$ as a pseudo-representation in characteristic $p$ and with constant determinant, and $A_\rhob$ denotes the Hecke algebra of modular forms modulo $3$ and level $1$, the natural map $R_\rhob \rightarrow A_\rhob$ is an isomorphism by \cite{annap3}, and both rings are isomorphic to $\F_3[[Y,Z]]$. Thus, the pseudo-deformation $(t,d)$ induces a surjective map $R_\rhob= A_\rhob \rightarrow A=\F_3[[X]]$, such that $(t,d)$ is the composition of the natural pseudo-character $(t_\rhob,d_\rhob)$ with this map.

\section{Density of modular forms}

\label{sectiondensity}

In this section we prove the main results of our works, the ones regarding the density of modular forms, namely Theorems~\ref{densnonzero}, \ref{densbounded} and \ref{ithmspecialinfinite}.

We revert to the notation of the introduction: $p$ is prime, $N \geq 1$ an integer, $k \in \Z/(p-1)\Z$ and $\F$ a (large enough) finite extension of $\F_p$. The space of modular forms on $\F$ of weight $k$, level $N$, and coefficients null at indices not prime to $Np$ is denoted by $\Fc$.
We note that to prove Theorems~\ref{densnonzero}, \ref{densbounded} and \ref{ithmspecialinfinite}, we can without loss of generality replace $\F$ by a finite extension. We shall always assume that the finite field $\F$ is large enough below.

\subsection{The Hecke algebra of mod $p$ modular forms}

The space $\Fc$ is endowed with an action of the Hecke operators $T_\ell$  for $\ell \nmid Np$. Let $A = A_k(N,\F)$  
be the topological closure\footnote{The topology on $\Fc$ is the discrete topology and the topology on $\End_{\F}(\Fc)$ is the compact-open topology} of the $\F$-subalgebra of $\End_{\F}(\Fc)$ generated by the Hecke operators $T_\ell$ for $\ell$ not dividing $Np$. 

For every $k \in \Z/(p-1)\Z$, the $\F$-algebra $A=A_k(N,\F)$ is {\it semi-local}. More precisely, if $\F$ is large enough,
its maximal ideals are in bijection with a certain set $\Rc=\Rc(k,N,\F)$ of semi-simple continuous Galois representations $\rhob: G_{\Q,Np} \rightarrow \Gl_2(\F)$ up to $\F$-isomorphism: the correspondence is given by $\lambda_\ell = \tr \rhob(\Frob_\ell)$. This set $\Rc(k,N,\F)$ can   be described as the set of all semi-simple representations $\rhob: G_{\Q,Np} \rightarrow \Gl_2(\F)$ of determinant $\omega_p^{k-1}$ and Serre's level $N$. This is the content of Serre's conjecture, now a theorem of Khare and Wintenberger.

Still assuming that $\F$ is large enough, and $\rhob \in \Rc(k,N,\F)$, we shall denote
by $A_\rhob$ the corresponding local component of $A=A_k(N,\F)$, that is the localization of $A_k(N,\F)$ at the maximal ideal corresponding to $\rhob$. The generalized eigenspace $\Fc_\rhob = \Fc_\rhob(N,\F)$ for the $T_\ell$, $\ell \nmid Np$, with generalized eigenvalues $\lambda_\ell$
(already considered defined in the introduction) is equivalently the localization of the $A=A_k(N,\F)$-module $\Fc=\Fc_k(N,\F)$ at that maximal ideal $\m_\rhob$ corresponding to $\rhob$. 

Then, $A_\rhob(\F)$ is a compact local $\F$-algebra with residue field $\F$. The image of the elements $T_\ell$  of $A$ in that localization $A_\rhob$
shall also by denoted by $T_\ell$. The image of $T_\ell \in A_\rhob$ in the residue field $\F$ is $\tr (\rhob(\Frob_\ell))=\lambda_\ell$. 
Equivalently, the $A_\rhob$-module  $\Fc_\rhob$ can be described as the generalized eigenspace in $\Fc_k(\F)$ for the $T_\ell$, $\ell \nmid Np$, with generalized eigenvalues $\lambda_\ell$. To summarize, we  have decompositions
\begin{eqnarray} \label{decFck} A= \prod_{\rhob \in \Rc} A_\rhob,\ \ \  \Fc = \bigoplus_{\rhob \in \Rc} \Fc_\rhob. \end{eqnarray}

Recall that we have a perfect pairing $\Fc \times A \rightarrow \F,\ (f,t) \mapsto a_1(tf)$, which induces a perfect pairing $\Fc_\rhob \times A_\rhob \rightarrow \F$.

We note that the ring $A$ thus satisfies all hypotheses made in Section~\ref{sectionessential}. Moreover we have the following results on the structure of $A$:

\begin{prop} \label{dimArhob} The ring $A_\rhob$ are always infinite, and have Krull dimension $\geq 1$. 
If $p>3$, or if $p=3$ and $\rhob$ is a twist of $1\oplus \omega_3$ ($\omega_3$ the cyclotomic character), or if $p=2$ and $\rhob$ is a twist of $1 \oplus 1$, the Krull dimension of $A_\rhob$ is at least $2$.
\end{prop}
\begin{pf} See \cite{jochnowitz} for the first assertion, \cite{BK} and \cite{Deo} for the case $p>3$ and \cite{annathesis} in the case $p=3$, \cite{NS2} in the case $p=2$.
\end{pf}
It is expected that $A_\rhob$ always has dimension exactly $2$, and this is known in many cases, see the references above.

\subsection{The canonical Galois pseudo-representation over $A$}

\begin{prop} \label{expseudo} There exists a unique continuous pseudo-representation  $(t,d)$ of dimension $2$ of $G_{\Q,Np}$ with values in $A$ such that $t(\Frob_\ell)=T_\ell$ for all $\ell \nmid Np$. One has $d = \omega_p^{k-1}$ and $t(c)=0$.
\end{prop}
For a proof of the proposition, which is well-known to specialists, see \cite{Bcras} where the case $p=2$ is dealt with -- the case $p>2$ is exactly the same. We denote by $(t_\rhob,d_\rhob)$ the composition of $(t,d)$ with the map $A \rightarrow A_\rhob$, and observe that by definition, $t_\rhob = \tr \rhob \pmod{\m_\rhob}$ and $d_\rhob=\det \rhob \pmod{\m_\rhob}$.

\begin{cor} The pseudo-deformation $(G_{\Q,Np},(\rhob_i)_{i=1,\dots,r},t,d)$ is admissible. \end{cor}
\begin{pf}
By \cite{mazur}, $G_{\Q,Np}$ satisfied the $p$-finiteness condition, hence~\ref{Pipfinite}.
The hypothesis~\ref{tmultfree} is satisfied because the representation $\rhob$ are odd, hypotheses~\ref{deform} and ~\ref{dconstant} are clear, and 
~\ref{tgenerate} follows from the fact that $t(G_{\Q,Np})$ contains $T_\ell$ for all prime $\ell$ not dividing $Np$ and those operators, by construction, generates $A$ as an $\F$-algebra.
\end{pf}

\begin{cor} The ring $A$ is noetherian. \end{cor}
\begin{pf} This follows from the preceding corollary and Cor.~\ref{admissiblenoetherian}. \end{pf}

We observe that if $p=2$, the ideal generated by all the $T_\ell$, $\ell \nmid Np$, in $A$ is the maximal ideal is the orthogonal of the eigenform  $\Delta$, which is up to a scalar the only form in $\Fc_k(\F)$ killed by all Hecke operators. We shall denote that ideal by $\m_1$
since it is the maximal ideal of $A$ corresponding to the trivial representation $\rhob=1 \oplus 1$.

\begin{lemma} \label{vecdense}   The closed $\F$-subspace generated by $t(G_{\Q,Np})$ is $A$ when $p>2$ and
$\m_1$ when $p=2$.
\end{lemma}
\begin{pf} When $p>2$, the  lemma is just~\ref{tgeneratemodule}. When $p=2$, the same argument gives that the closed $\F$-subspace generated by $t(G_{\Q,Np})$ is an ideal, and contains all the $T_\ell$, $\ell \nmid Np$. Thus it is $\m_1$ or $A$. But $t(G_{\Q_{Np}}) \subset \m_1$ because $t \pmod{\m_1} = \tr (1+1)=0$.
\end{pf}

\subsection{Proof of Theorem~\ref{densnonzero}}

\label{proof1}

We now give the proof of Theorem~\ref{densnonzero}. 
Let $f \in \Fc_k(\F)$, $f \neq 0$. If $p=2$ we assume in addition that $f \not \in \F \Delta$. We want to show that $\delta(f)>0$.

 Let $l_f$ be the $\F$-linear form on $A_k(\F)$ defined by $l_f(T)=a_1(Tf)$. In other words, $l_f$ 
is the linear form on $A_k(\F)$ corresponding to $f \in \Fc_k(\F)$ through the perfect duality $A_k(\F) \times \Fc_k(\F) \rightarrow \F$, $(T,f) \mapsto a_1(Tf)$, 
and in particular, $l_f$ is non-zero. Let $H_f$ be the closed hyperplane $\ker l_f$ of $A_f(\F)$. If $p=2$, our supplementary assumption means that $H_f$ is not the maximal ideal $\m_1$.

If $\mu$ denotes the Haar measure of total mass $1$ on the compact group $G_{\Q,Np}$, we claim that
\begin{num} \label{deltaf} $\delta(f) = 1- \mu( t^{-1}(H_f) ).$ \end{num}
To prove the claim, note that for $\ell$ a prime not dividing $Np$, one has $a_\ell(f)=0 \Leftrightarrow a_1(T_\ell f)=0 \Leftrightarrow a_1(t(\Frob\ell) f)=0 \Leftrightarrow t(\Frob_\ell) \in H_f \Leftrightarrow \Frob_\ell \in t^{-1}(H_f)$. Observe that $H_f$, being closed and of finite index, is open in $A_f$, and therefore $t^{-1}(H_f)$ is open in $G_{\Q,Np}$. Thus Chebotarev's density theorem implies that the density of primes $\ell$ such that $\Frob_\ell$ is not in $t^{-1}(H_f)$ is $1- \mu( t^{-1}(H_f) )$, and the claim follows

To finish the proof, we therefore just have to prove that $t^{-1}(H_f)$ is a proper subset of  $G_{\Q,Np}$. We do not have  $t^{-1}(H_f) = G_{\Q,Np}$, because that would mean $t(G_{\Q,Np}) \subset H_f$, contradicting Lemma~\ref{vecdense}. This completes the proof of Theorem~\ref{densnonzero}.

\subsection{Definition of special modular forms}

\label{defformespeciales} 

{\it From now on, we assume $p>2$.} The admissible pseudo-deformation $(G_{\Q,Np},(\rhob_i),t,d)$ over $A$ defines a closed $\F$-subspace $A_\ess$ of $A$ (cf. \S\ref{sectionessential}). We say that a modular form $f \in \Fc$ is {\it special} if $a_1(tf)=0$ for all $t \in A_\ess$. Thus, special modular forms in $\Fc$ form  a $\F$-sub-vector space $\Fc_\spe$, which is the orthogonal of $A_\ess$ for the perfect pairing $A \times \Fc \rightarrow \F$.

For $\rhob \in \Rc$, we set as in the introduction $\Fc_{\rhob,\spe} = \Fc_\rhob \cap \Fc_\spe$. The admissible pseudo-deformation
$(G_{\Q,Np},\rhob,t_\rhob,d_\rhob)$ over $A_\rhob$ defines a closed $\F$-subspace $A_{\rhob,\spe}$ of $A_\rhob$, which is the image of $A_\spe$ by the projection mao $A \rightarrow A_\rhob$. Thus, $\Fc_{\rhob,\spe}$ is the orthogonal of $A_{\rhob,\ess}$ for the perfect pairing $A_\rhob \times \Fc_\rhob \rightarrow \F$. 

\subsection{Proof of Theorem~\ref{ithmspecialinfinite}}

\label{Proofthmspecialinfinite}
\label{proofithmspecialinfinite}
Given a representation $\rhob \in \Rc$ (which in the case $p=3$ is a twist of $1 \oplus \omega_3$), we need to show that $\Fc_{\rhob,\spe}$ is of infinite codimension in $\Fc_\rhob$, or equivalently, that $A_{\rhob,\ess}$ is infinite-dimensional.

\begin{prop} If $(G_{\Q,Np},\rhob,t,d)$ is a virtually abelian admissible pseudo-deformation over a noetherian local compact domain $A$ such that $pA=0$ for some odd prime $p$, then the Krull dimension of $A$ is at most 1.
\end{prop}
\begin{pf} Let $K$ the fraction field of $A$. Let $A'$ be the integral closure of $A$ in $K$. Since $A$ is a complete noetherian local ring, then by a theorem of Nagata (see \cite[Corollary 2 page 234]{matsumura}), $A'$ is a finite type module over $A$ and is a complete noetherian local ring as well.

Let $\rho:\G_{\Q,Np} \rightarrow R^\ast$ be a $(t,d)$-representation. By Lemma~\ref{basicGMA3}, $\rho$ can be seen as a representation $G_{\Q,Np} \rightarrow \Gl_2(K)$. 

Let $M$ be a finite Galois extension of $\Q$ such that $(t,d)$ factors through an abelian quotient of $G_{M,Np}$. 
The representation $\rho: G_{M,Np} \rightarrow \GL_2(K)$ becomes reducible over a quadratic extension $K'$ of $K$, so there are two 
characters $\chi_1,\chi_2: G_{M,Np} \rightarrow (K')^\ast$ such that $\rho=\chi_1 \oplus \chi_2$ as a representation over $K'$.
Since $\chi_i(g)$ for $i=1,2$ are the roots of the polynomials $X^2-t(g) X + d(g) \in A[X]$, $\chi_i(g)$ belongs to the integral 
closure $A'$ of $A$ in $K'$. Since $A$ is a complete noetherian local ring, then by a theorem of Nagata, $A'$ is a finite type module over $A$ and is a complete noetherian local ring as well. 

We claim that the characters $\chi_i: G_{M,Np} \rightarrow (A')^\ast$ for $i=1,2$ are continuous. Indeed, if they are equal  they are continuous  since $2 \chi_1(g)=t$.
If not,  there is a $g_0$ such that $\chi_1(g_0)\neq\chi_2(g_0)$. By the continuity of the roots of polynomial, there exists a neighborhood 
$U$ of $1$ and two continuous functions $\psi_1,\psi_2$ on $g_0 U$ (with values in $(A')^\ast$) such that $X^2-t(g)X +d(g)=(X-\psi_1(g))(X-\psi_2(g))$ on $g_0 U$ $\psi_i(g_0) =\chi_i(g_0)$ for $i=1,2$. Shrinking $U$ if necessary, we may assume that $U$ is an open 
subgroup of $G_{M,Np}$ and that $g \mapsto \psi_i(gg_0) \chi_i(g_0)^{-1}$ is a character on $U$. By uniqueness of the decomposition of a representation into sum of characters over a field ($K'$),  it follows that for $i=1,2$, there exists $j=1,2$ such that $\psi_i(gg_0) \chi_i(g_0)^{-1} = \chi_j(g)$ on $U$. It follows that the $\chi_i$ are continuous on $U$, hence everywhere.

Let $\Gamma = \Gal(M/\Q)$. Since the functions $t$ and $d=\chi_1\chi_2=\det \rhob$ on $G_{M,Np}$ are invariant by conjugation of the argument by any element of $\Gal(M/\Q)$, there exists a subgroup $\Gamma'$ of $\Gamma$ of index $1$ or $2$ such that
\begin{num} \label{condchi1}
for every $\gamma \in \Gamma'$, $\chi_i^\gamma = \chi_i$ and for every $\gamma \in \Gamma-\Gamma'$, $\chi_i^{\gamma} = \det \rhob \chi_i^{-1}$.
\end{num} 
\noindent Let $R_{\univ}$ be the universal deformation ring in characteristic $p$ of the character $\bar \chi_1: G_{M,Np} \rightarrow \F^\ast$ satisfying condition~\ref{condchi1}. The character $\chi_1: G_{M,Np} \rightarrow (A')^\ast$ defines a morphism of $\F$-algebras $R_\univ \rightarrow A'$ whose image $A_0$ is the closed $\F$-subalgebra of $A'$ generated by $\chi_1(G_{M,Np})$. For $g \in G_{M,Np}$, we can
write $\chi_1(g) = \bar \chi_1(g) + x$ with $\bar \chi_1(g) \in \F^\ast$ and $x$ in the maximal ideal of $A'$, and $\chi_2(g)=\det \rhob(g) (\bar \chi_1(g) + x)^{-1} = \bar \chi_2(g) (1 - \bar \chi_1(g)^{-1} x + \bar \chi_2(g)^{-2} x^2 - \dots)$. Thus $\chi_2(g)$ is in $A_0$, 
and so is $t(g)$. Since $A$ is the closed $W(\F)$-subalgebra generated by the image of $t$, we see that $A \subset A_0 \subset A'$.
Since $A'$ is finite as an $A$-module, Cohen-Seidenberg's theorem ensures that $A$, $A'$ and $A_0$ have the same Krull's dimension.
Thus to prove the proposition it suffices to prove that $A_0$ has dimension at most $1$, and for this it is enough to prove that $R_\univ$ has dimension at most 1. This follows easily from Class Field Theory.
 \end{pf}

By Prop.~\ref{dimArhob}, the ring $A_\rhob$ has Krull dimension at least $2$ under the hypothesis on $\rhob$ of Theorem~\ref{ithmspecialinfinite}. Let $B$ be the reduced ring of the ring of a $2$-dimensional irreducible component of $\spec A_\rhob$. Then
$B$ is a quotient of $A_\rhob$, which is domain of dimension $2$. To prove that $A_{\rhob,\ess}$ is infinite (as a set or $\F$-vector space),
it is enough to prove that the image $B_\ess$ in $B$ is infinite. The subspace $B_\ess$ is the essential subspace of the admissible pseudo-deformation $(\G_{\Q,Np},\rhob,t,d)$ over $B$, which is not virtually abelian by the above proposition. Therefore, $B_\ess$ is infinite by Propositions~\ref{Aesslargered}, \ref{Aesslargedih}, \ref{Aesslargeirr}, and Theorem~\ref{ithmspecialinfinite} is proved.

\subsection{Proof of Theorem~\ref{densbounded}} 

\label{proofdensbounded}

Let $f \in \Fc$ a modular form which is not in $\Fc_\spe$. This means that the modular form 
$l:A \rightarrow \F$, $t \mapsto a_1(tf)$ is not zero on the subspace $A_\ess$ of $A$.
By Theorem~\ref{compmeasure}
$$\mu_{G_{\Q,Np}}( (l \circ t)^{-1}(\F^\ast)) \geq  \frac{p-1}{pn},$$
that is by \ref{deltaf}
$$\delta(f) \geq  \frac{p-1}{pn},$$
where $n = |\Gb|$. This proves the main part of Theorem~\ref{densbounded}. This theorem also states that $\Fc_\spe$ is of infinite codimension in $\Fc$. To prove this, it is sufficient that for one $\rhob \in \Rc$, $\Fc_{\spe,\rhob}=\Fc_\spe \cap \Fc_\rhob$ is of infinite codimension
in $\Fc_\rhob$. The results from Theorem~\ref{ithmspecialinfinite} for any $\rhob \in \Rc$ if $p > 3$, and also for $p=3$ if we choose
for $\rhob$ the representation $1 \oplus \omega_3$, which always belong to $\Rc(N,3,\F)$ since it is the representation attached to the eigenform $\Delta \pmod{3}$.

\section{Cyclotomic and $K$-abelian modular forms}

\label{sectionspecial}

We keep the notation of the preceding section. We do not assume $p>2$ unless explicitly mentioned. We fix a representation $\rhob \in \Rc$.

For $f \in \Fc_\rhob$, we denote by $I_f$ the annihilator ideal of $f$ in $A_\rhob$, and
by $A_f$ the quotient $A_\rhob/I_f$. The perfect duality $A_\rhob \times \Fc_\rhob \rightarrow \F$ induces a perfect duality
$A_f \times A_\rhob f \rightarrow \F$. The space $A_\rhob f$ is finite, because the action of the Hecke operators is locally finite; it follows that 
the ring $A_f$ is finite, and it is therefore a local artinian $\F$-algebra.
We obtain an admissible pseudo-deformation $(G_{\Q,Np},\rhob,t_f,d_f)$ on $A_f$  by post-composing 
$t_\rhob$ and $d_\rhob$ with the surjective map  $A_\rhob \rightarrow A_f$. 

\subsection{Fields of determination of a modular form $f \in \Fc_\rhob$}

For $S$ a finite set of primes, let us denote by $\Q_S$ the maximal algebraic extension of $\Q$ unramified outside $S$ and $\infty$, and by $G_{\Q,S}$ the group $\Gal(\Q_S/\Q)$. If $S$ is the set of primes dividing an integer $N$,
we also use $N$ instead of $S$ in these notations.

Let us denote by $L_f$ the subfield of $\Q_{Np}$ fixed by $\ker t_f$. Note that $L_f$ is a Galois extension of $\Q$, unramified outside $Np$ and $\infty$, such that 
$\Gal(L_f/\Q)=G_{\Q,Np}/\ker(t_f,d_f)$. 

\begin{lemma} The field $L_f$ is  finite over $\Q$. \end{lemma}
\begin{pf} Let $\rho_f: G_{\Q,Np} \rightarrow R_f^\ast$ be a $(t_f,d_f)$-representation. 
By assertion (vii) of Proposition~\ref{exGMA}, $R_f$ is of finite type as a module over $A_f$, hence is finite as a set, and by assertion (vi) of the same, the kernel of $\ker (\rho_f)_{|G_{\Q,Np}}$ is $\ker (t_f,d_f)$.
Therefore, $\Gal(L_f/\Q)=\rho_f(G_{\Q,N,p})$ and since the later is a subset of the finite set $R_f$, it is finite.
\end{pf}

\begin{theorem} \label{theoremdetermination} Let $L$ be a Galois extension of $\Q$ contained in $\bar \Q$, unramified outside a finite set $S$ of primes dividing $Np$ and $\infty$.
The following properties are equivalent:
\begin{itemize}
\item[(i)] For every prime $\ell  \not \in S$, the form $T_\ell f$ depends on $\ell$ only through the conjugacy class
$\Frob_{\ell,L/\Q} \in\Gal(L/\Q)$.
\item[(i')] For almost every prime $\ell$, the form $T_\ell f$ depends on $\ell$ only through the conjugacy class
$\Frob_{\ell,L/\Q} \in \Gal(L/\Q)$.
\item[(ii)]  For every prime $\ell  \not \in S$, the coefficient $a_\ell(f)$ depends on $\ell$ only through the conjugacy class
$\Frob_{\ell,L/\Q} \in\Gal(L/\Q)$.
\item[(ii')]  For almost every prime $\ell$, the coefficient $a_\ell(f)$ depends on $\ell$ only through the conjugacy class
$\Frob_{\ell,L/\Q} \in\Gal(L/\Q)$.
\item[(iii)] One has $L_f \subset L$.
\end{itemize}
\end{theorem}
\begin{definition} If $L$ satisfies the conditions of the above theorem, we shall say that $L$ is a {\it determination field} of $f$.
\end{definition}
Obviously, there is always a smallest determination field, namely $L_f$, and it is finite over $\Q$ and unramified outside $Np$. However, it is sometimes convenient to consider also other determination fields.  

\begin{proof} 
We see the pseudo-representation $(t_f,d_f)$ of $\Gal(\Q_{Np}/\Q)$ as a pseudo-representation
of $\Gal(\Q_S/\Q)$ by inflation. Let us call $\pi$ the surjective map $\Gal(\Q_S/\Q) \rightarrow \Gal(L/\Q)$. To ease notations,
let us denote by $\Frob_\ell$ the element $\Frob_{\ell,\Q_s/\Q}$. Thus $\pi(\Frob_\ell)=\Frob_{\ell,L/\Q}$. 

Since $t_{f}(\Frob_{\ell} f) = T_\ell f$, the assertion (i) (resp. (i')), is equivalent to 
\begin{num}  $t_{f}(\Frob_{\ell})$ depends only on $\pi(\Frob_{\ell}) = \Frob_{\ell,L/\Q}$  for all $\ell$ not in $L$ (resp. for almost all $\ell$) 
\end{num}
By Chebotarev's density theorem, both these assertions are equivalent to:
\begin{num}
The map $t_{f}$ factors through $\pi$,
\end{num}
\noindent
which amounts to $\ker \pi \subset \ker t_{f}$, that is $L_f \subset L$. We thus have proved the equivalence between (i), (i') and (iii).

Since the coefficient $a_\ell$ of $f$ is the coefficient $a_1$ of $T_\ell(f)$, it is obvious that (i) implies (ii). Since (ii) obviously implies  (ii'), it just remains to prove  that (ii') implies (i'). 
For every prime $\ell$ not in $S$, one has
$$a_1(T_f(\Frob_{\ell}) f) = a_1(T_\ell f) = a_\ell(f),$$
so (ii') means that for almost all $\ell$, $a_1(t_f(\Frob_{\ell}) f)$ depends only on $\Frob_{\ell,L/\Q}=
\pi(\Frob_{\ell})$. Using Chebotarev, this means that there exists  a continuous map $\beta: \Gal(L/\Q) \rightarrow \F$ such that
\begin{num} \label{facalpha} for all $\gamma \in \Gal(\Q_S/\Q)$, $a_1(t_f(\gamma)f) = \beta(\pi(\gamma))$. \end{num}
\noindent
Let $q$ be a prime number not in $S$. 
\begin{eqnarray*} a_q(T_\ell f) &=& a_1 (T_\ell T_q f) \\ &=& a_1( t_{f}(\Frob_{\ell,\Q_S/\Q}) t_f(\Frob_{q,\Q_S/\Q}) f) \\
&=& a_1(t_{f}(\Frob_{\ell,\Q_S/\Q} \Frob_{q,\Q_S/\Q} f) + q^{k-1} a_1( t_{f}(\Frob_{\ell,\Q_S/\Q} \Frob^{-1}_{q,\Q_S/\Q}) f)\\
&=& 
\beta (\Frob_{\ell,L/\Q} \Frob_{q,L/\Q}) + q^{k-1} \beta( \Frob_{\ell, L/\Q} \Frob^{-1}_{q,L/\Q})
\end{eqnarray*}
Thus the coefficient $a_q$ (for $q$ any prime not in $S$), as well as the coefficient $a_1$ of the form $T_\ell f$
depends on $\ell$ only through $\Frob_{\ell,L/\Q}$. Since by the corollary of Theorem~\ref{densnonzero} a modular form is determined by its coefficient at primes (excluding a finite set) and at $1$ , it follows that the form $T_\ell f$ itself depends on $\ell$ only through $\Frob_{\ell,L/\Q}$. In other words, we have proved (i').
\end{proof}

\subsection{Cyclotomic modular forms}

\begin{prop} Let $f=\sum_n a_n q^n \in \Fc_\rhob(\F)$. The following are equivalent:
\begin{itemize}
\item[(i)] $f$ has a determination field which is abelian over $\Q$.
\item[(ii)] There exists an integer $M \geq 1$ such that for all prime $\ell$ not dividing $Np$, $a_\ell$ depends on $\ell$ only trough $\ell \pmod{M}$.
\item[(iii)] There exists an integer $M \geq 1$ such that for all prime $\ell$ not dividing $Np$ $T_\ell f$,  depends on $\ell$ only trough $\ell \pmod{M}$.
\end{itemize}
If they hold, we can take $M$ in (ii) and (iii) such that all prime factors of $M$ divide $Np$.
\end{prop}
\begin{pf} This is a special case of Theorem~\ref{theoremdetermination}, taking into account the Kronecker-Weber theorem that every number field abelian over $\Q$ is a subfield of a cyclotomic field $\Q(\zeta_M)$. \end{pf}
\begin{definition} We say that $f$ is {\it cyclotomic}  if it satisfies the conditions of the above proposition.
\end{definition}
 
\begin{definition} Let us denote by $I_\cycl$ the ideal generated by the elements $t_\rhob(xyx^{-1}y^{-1}s)-t_\rhob(s)$
for $x,y,s \in G_{\Q,Np}$. 
\end{definition}
Since $A_\rhob(\F)$ is noetherian the ideal $I_\cycl$ is finitely generated and closed.  Clearly, $I_\cycl$ is the smallest ideal $I$ of $A_\rhob$
such that $G/\ker (t_I, d_I)$ is abelian, where $t_I$ is the composition $t : G \rightarrow A \rightarrow A/I$ and similarly for $d$.

\begin{example} In the case $p=2$, $\rhob=1 \oplus 1$, the ideal $I_\cycl$ is principal, and generated by the element $T_5+T_3+T_3^3+T_3^5+T_3^9+T_3^{11}+T_3^{129}+\dots$: see \cite{Bcras}. 
\end{example}

\begin{prop} A form $f$ is cyclotomic if and only if it is annihilated by $I_\cycl$.
\end{prop} 
\begin{pf} A form $f$ is killed by $I_\cycl$ if and only if $I_\cycl \subset I_f$ which is visibly equivalent to $G_{\Q,Np}/\ker t_f$ being abelian, or $L_f$ being an abelian extension of $\Q$.
\end{pf}
\begin{prop} If $\rhob$ is irreducible, the only cyclotomic form in $\Fc_\rhob(\F)$ is $0$.
\end{prop} 

\begin{pf} Recall that if $\rhob$ is irreducible, it is absolutely irreducible, hence its image $\rhob(G_{\Q,Np})$ is not abelian. 
If there is a non-zero cyclotomic form $f$ in $\Fc_\rhob(\F)$, then the pseudo-representation $(t_f,d_f) : G_{\Q,Np} \rightarrow A_f$ reduces modulo the maximal ideal of $A_f$ to the pseudo-representation $(\tr \rhob, \det \rhob) : G_{\Q,Np} \rightarrow \F$, and it follows that the group 
$G_{\Q,Np} / \ker t_\rhob$ is a quotient of $G_{\Q,Np}/\ker t_f$, hence is abelian. Since $\rhob$ is semi-simple, $\ker \rhob = \ker \tr \rhob$,
hence $G_{\Q,Np} / \ker t_\rhob \simeq \rhob(G_{\Q,Np})$ is abelian, a contradiction.
\end{pf}

For the rest of this subsection we assume that $p>2$ (for similar but more complicated results in the case $p=2$, $N=1$, see \cite{Bcras}), and that the projective image of $\rhob$ is cyclic, in other words that $\rhob$ is reducible. 
Let $\rho: G_{\Q,Np} \rightarrow R^\ast$ be a $(t,d)$-representation with $R=\mat{A & B \\ C & D}$. 

\begin{prop} One has $I_\cycl = BC$. In other terms, $I_\cycl$ is just the reducibility ideal of the pseudo-representation $t_\rhob$ (see \cite[\S1.5]{BC}).
\end{prop}
\begin{pf} Let $I$ be any ideal of $A_\rhob$, and let $(G_{\Q,Np},\rhob,t_I,d_I)$ be the admissible pseudo-deformation 
obtained by reducing the pseudo-deformation over $A_\rhob$ modulo $I$. Let
$\rho_I : G_{\Q,Np} \rightarrow R_I^\ast$ be a $(t_I,d_I)$-representation
attached to the admissible pseudo-deformation $(G_{\Q,Np},\rhob,t_I,d_I)$ adapted to an element of $G_{\Q,Np}$ for which 
$\rho$ is also adapted.
Then $R_I = \mat{A/I & B_I \\ C_I &A/I}$ and there is a natural surjective morphism of algebras $R \otimes_A A/I = R/IR \rightarrow R_I$
inducing identity maps $A/I \rightarrow A/I$ on the diagonal components, and maps $B/IB \rightarrow B_I$, $C/IB \rightarrow C_I$ on the non-diagonal components. Note that the map $R/IR \rightarrow R_I$, as well as the maps $B/IB \rightarrow B_I$ and $C/IC \rightarrow C_I$
needs not be injective (this is because $R/IR$ may not be faithful.) The ideal $B_I C_I$ of $A/I$ is nevertheless the image in $A/I$ of the ideal $BC$ of $A$, because the map $R \mapsto R_I$ preserves multiplication of matrices (see~\cite[\S1.5]{BC} for more detailed proofs of the assertion of this paragraph).

By construction $I_\cycl$ is the smallest ideal $I$ of $A$ such that $G_{\Q,Np}/\ker t_I$ is abelian. One has $\ker (\rho_I)_{|G} = \ker t_I$ because $R_I$ is faithful. Hence $G/\ker t_I \simeq \rho_I(G)$, and $I_\cycl$ is the smallest ideal $I$ of $A$ such that $\rho_I(G)$ is abelian, or again, since $R_I$ is generated by $\rho_I(G)$ as an $A/I$-module, the smallest ideal $I$ such that $R_I$ is commutative. It is easy to see that the GMA $R_I = \mat{A/I & B_I \\ C_I & A/I}$ is commutative if and only if $B_I=C_I=0$. Since the product $B_I \times C_I \rightarrow A/I$ is a non-degenerate pairing, this is equivalent to $B_I C_I = 0$, that is by the above paragraph, to $BC \subset I$. Thus $BC=I_\cycl$.
\end{pf}
 
\begin{cor}\label{corspered} Assume as above that the projective image of $\rhob$ is cyclic, but also that it is not of order 2. Then $I_\cycl=A_{\rhob,\ess}$.
In other words, a form $f \in \Fc_\rhob$ is cyclotomic if and only if it is special.
\end{cor}
\begin{pf} This follows from the preceding proposition and Theorem~\ref{thmAessred}.
\end{pf}

\subsection{$K$-abelian forms}

In this subsection, we assume $p>2$. For $K$-abelian forms in the case $p=2$, see \cite{NS2} and an article in preparation by J. Bellaïche, J.-L. Nicolas, and Jean-Pierre Serre.
Let $K$ be a quadratic extension of $\Q$.

\begin{definition} A form $f \in \Fc_\rho$ is {\it $K$-abelian} if it has a field of determination $L$ which is an abelian extension of $K$.\end{definition}

Note that the composition of two Galois extensions of $\Q$ which contains $K$ and are abelian over $K$ is also a Galois extension of $\Q$ which contains $K$ and is abelian over $K$. It follows that if $f$ and $f'$ are $K$-abelian, $f+f'$ is $K$-abelian as well:
if $L$ and $L'$ are fields of determination of $f$ and $f'$, then $LL'$ is a field of determination of $f+f'$. Thus the set of $K$-abelian forms is a vector space.
It is also obviously stable by the Hecke operators $T_\ell$. Hence its orthogonal for the duality $A_\rhob \times \F_\rhob \rightarrow \F$ is an ideal $I_{\Kab}$. 

From now on, {\bf we assume that the projective image of $\rhob$ is dihedral of order $>4$}. Thus the projective image of $\rhob$ has a unique quotient of order $2$, which corresponds to a quadratic extension $K$ of $\Q$. Thus $G_{K,Np}$ is a subgroup of index $2$ in $G_{\Q,Np}$ and the projective image of $\rhob(G_{K,Np})$ is cyclic. We choose a well-adapted $(t_\rhob,d_\rhob)$-representation $\rho: G_{\Q,Np} \rightarrow \GL_2(A_\rhob)$. 
By \S\ref{consequencesdihedral}, the $A_\rhob$ algebra generated by $\rho(G_{K,Np})$ is a sub-GMA of $M_2(A_\rhob)$, of the form 
$R=\mat{A_\rhob & B \\ B & A_\rhob}$ for some proper ideal $B$ of $A_\rhob$. 

\begin{prop} \label{IBKab} One has $B=I_{\Kab}$.
\end{prop}
\begin{pf} By definition, $I_{\Kab}$ is the smallest ideal $I$ of $A_\rhob$ such that the image of $G_{K,Np}$ in the quotient $G_{\Q,Np}/\ker(t_I,d_I)$ is abelian. Since the representation $\rho_I : G_{\Q,Np} \rightarrow \GL_2(A_\rhob/I)$ obtained by reducing $\rho$ modulo $I$ has trace $t_I$ and determinant $d_I$, and since the GMA $M_2(A/I)$ is faithful, $\rho_I$ realizes an isomorphism $G_{\Q,Np}/\ker(t_I,d_I) \rightarrow \rho_I(G_{\Q,Np}) \subset \GL_2(A_\rhob/I)$, and the image of $G_{K,Np}$ into $G_{\Q,Np}/\ker (t_I,d_I)$ is $\rho_I(G_{K,Np})$.
Thus, $I_\Kab \subset I$ if and only if the group $\rho(G_{K,Np})$ is abelian, if and only if the $A_\rhob/I$-subalgebra of $M_2(A_\rhob/I)$ generated by $\rho_I(G_{K,Np})$is commutative, if and only if the image of $R=\mat{A_\rhob & B \\ B & A_\rhob}$ in $M_2(A_\rhob/I)$.
Clearly, the latter condition is equivalent to $B \subset I$. Thus $B=I_{\Kab}$.
\end{pf}

\begin{cor} Assume that the the projective image of $\rhob$ is dihedral of order $>4$, and divisible by $4$.
The one has $A_{\rhob,\ess} = \m_\rhob I_\Kab$.
\end{cor}
\begin{pf} By Theorem~\ref{thmAessdih}, one has $A_{\rhob,\ess}= \m_\rhob B$. The corollary follows.
\end{pf}

\begin{cor} \label{corspedih} Assume that the the projective image of $\rhob$ is dihedral of order $>4$, and divisible by $4$. 
A form $f \in \Fc_\rhob$ is special if and only if $(T_\ell - \lambda_\ell)f=0$ is $K$-abelian for every prime $\ell$
not dividing $Np$ (here $\lambda_\ell=\tr(\rhob(\Frob_\ell))$). In particular, the space $\Fc_{\rhob,\spe}$ contains the space of $K$-abelian forms as a finite dimensional subspace.
\end{cor}
\begin{pf} This is just a translation of the preceding corollary, using that $\m_\rhob$ is finitely generated since $A_\rhob$ is noetherian.
\end{pf}

\subsection{The case of a large or exceptional $\rhob$} 

\label{spelarge}

In this case it is easy to see that there is no cyclotomic forms in $\Fc_\rhob$, nor $K$-abelian forms for any $K$.
The space $A_{\rhob,\ess}$ is the ideal $\m_\rhob^2$, hence $\Fc_{\rhob,\spe}$ is the space of forms which are killed by 
$(T_\ell-\lambda_\ell)^2$ for all $\ell$ not dividing $Np$. This space is finite-dimensional since $A_\rhob/\m_\rhob^2$ is finite-dimensional.

 \end{document}